\documentclass[12pt,twoside,leqno]{amsbook}
\usepackage{amssymb}

\input xy
\xyoption{all}
\usepackage[all]{xy}
\usepackage{graphicx}

\def\R{{\mathbb R}}
\def\C{{\mathbb C}}
\def\H{{\mathbb H}}

\newcommand\mc[1]{{\mathcal{#1}}}

\newcommand\Lie[1]{{\mathfrak{#1}}}
\newcommand\group[1]{{\mathrm{#1}}}

\def\Aut{\mathop{\rm Aut}}
\def\Int{\mathop{\rm Int}}
\def\Ad{\mathop{\rm Ad}}
\def\ad{\mathop{\rm ad}}
\def\id{\mathop{\rm id}}

\def\eps{\varepsilon}

\def\sfrac#1#2{{\scriptstyle{\frac{#1}{#2}}}}

\newlength{\pictht}

\newtheorem{theorem}{\sc Theorem}
\newtheorem{cor}[theorem]{\sc Corollary}
\newtheorem{lemma}[theorem]{\sc Lemma}
\newtheorem{prop}[theorem]{\sc Proposition}

\begin{document}

\pagenumbering{}

\phantom{.}

\vskip1truecm

\centerline{DOTTORATO DI RICERCA IN MATEMATICA E APPLICAZIONI}

\centerline{XV
CICLO}

\vskip1truecm

\centerline{Sede Amministrativa: Universit\`a di Genova}

\centerline{Sedi consorziate:  Politecnico
di Torino, Universit\`a di Torino}

\vskip1truecm

\centerline{PhD Thesis}

\centerline{May 2004}

\vskip2truecm

\centerline{Alessandro Ottazzi}

\vskip2truecm

\centerline{ \LARGE\bf Multicontact mappings }
\vskip0.3cm
\centerline{ \LARGE\bf on Hessenberg manifolds}
\vskip 3truecm
\begin{align*}
&\text{Supervisor:}  &\hskip6cm  &\text{PhD program coordinator:}\\
&\text{Prof. Filippo De Mari} &\hskip6cm &\text{Prof. Claudio Pedrini}\\
&\text{(University of Genova)} &\hskip6cm &\text{(University of Genova)}
\end{align*}

\pagenumbering{}

\pagenumbering{roman}
\tableofcontents
\chapter*{Introduction}
\vskip0.5cm
This thesis is concerned with the study of multicontact structures on Hessenberg
manifolds and of the mappings that preserve them. The setting in which our
considerations take place is that of parabolic geometry, namely a homogeneous space
G/P where G is a semisimple Lie group  and P is a parabolic subgroup
of G. Furthermore, it is assumed that G has real rank greater than one
and that P is minimal. We show that it is possible
to define a notion of multicontact mapping, hence of multicontact vector
field, on every
Hessenberg submanifold ${\rm Hess}_{\mc{R}}(H)$ of G/P associated to a
regular element $H$ in the Cartan subspace $\Lie{a}$ of the Lie algebra
$\Lie{g}$ of G. The Hessenberg combinatorial data, namely  the subset
$\mc{R}$ of the positive restricted roots $\Sigma_+$ relative to
$(\Lie{g},\Lie{a})$ that defines the type of the manifold, single out an
ideal
$\Lie{n}_{\mc C}$ in the nilpotent Iwasawa subalgebra of $\Lie{g}$, labeled by the complement
$\mc C=\Sigma_+\setminus\mc{R}$. By means of a reduction theorem, it is
shown that without loss of generality one can work under the assumption
that $\mc{R}$ contains all the simple restricted roots. In order to avoid
certain degeneracies, we assume further that $\mc{R}$ contains all
height-two restricted roots as well. We prove that the normalizer of
$\Lie{n}_{\mc C}$ in $\Lie{g}$ modulo $\Lie{n}_{\mc C}$ is naturally
embedded in the Lie algebra of multicontact vector fields on ${\rm
Hess}_{\mc{R}}(H)$. If the data $\mc{R}$ satify the property of encoding
a finite number of positive root systems, each corresponding to an
Iwasawa nilpotent algebra, then the above quotient actually coincides
with the Lie algebra of multicontact vector fields on ${\rm
Hess}_{\mc{R}}(H)$. This situation covers a wide variety of cases (for
example all Hessenberg data in a root system of type $A_\ell$) but not
all of them. Explicit exceptions are given in the $C_\ell$ case. One of
the main motivations for the present study is the observation that
${\rm Hess}_{\mc{R}}(H)$ can be realized locally as a stratified
nilpotent group that is {\it not} always of Iwasawa type. Hence our work
is an extension of the theories of multicontact maps developed thus far.
\vskip0.5truecm
The notion of multicontact structure was introduced in \cite{CDKR1} and \cite{CDKR2}
in the context of the homogeneous spaces G/P. Roughly speaking, it
refers to a collection of special sub-bundles of the tangent bundle with
the property that their sections generate the whole tangent space by
repeated brackets. The selection of the special directions is not only
required to satisfy this H\"ormander-type condition, but it is also
dictated by the stratification of the tangent space $T_x$ at each point
$x\in$ G/P in terms of restricted root spaces. If for example P is
minimal, then  $T_x$ can be identified with a nilpotent Iwasawa Lie
algebra and therefore it may be viewed as  the direct sum of all the root
spaces associated to the positive restricted roots. Since a positive root
is a sum of simple roots, it is natural to expect that the tangent
directions along the simple roots will play a special r\^ole. Indeed, it
is proved in
 \cite{CDKR2} that, at least in rank greater than one,  G acts on G/P
by maps whose differential preserves each sub-bundle corresponding to a  simple restricted root,
or, at worst, it permutes them amongst themselves.
It is thus natural to say that $g\in \group{G}$ induces a {\it
multicontact mapping}. The main result in \cite{CDKR2} is that the
converse statement is also true: a locally defined $C^2$
multicontact mapping on G/P is the restriction of the action of a
uniquely determined element
$g\in \group{G}$. Hence the boundaries G/P are (in most cases) {\it
rigid}.
\par
This type of theorem is one in a long standing history of rigidity results, dating back to Liouville.
Around 1850, he proved that any $C^4$ conformal map between domains
in~$\R^3$ is necessarily a composition of translations, dilations
and inversions in spheres. This amounts to saying that the group $\group{O}(1,4)$ acts
on the sphere~$S^3$ by conformal transformations (and hence locally
on~$\R^3$, by stereographic projection), and then proving that any
conformal map between two domains arises as the restriction of the
action of some element of $\group{O}(1,4)$. The same result also holds
in~$\R^n$ when $n > 3$ (see, for instance,~\cite{Nv60}),
and with metric rather than smoothness assumptions (see~\cite{Gg62}).
\par
A cornerstone in the extension process of Liouville's result is
certainly the paper \cite{KR85} by A. Kor\'anyi and H.M. Reimann, where
the Heisenberg group $\H^n$ substitutes
the Euclidean space and the sphere in $\C^n$ with its Cauchy-Riemann
structure substitutes the real sphere. The authors study
smooth maps whose differential preserves the contact (``{\it horizontal}'')
plane $\R^{2n}\subset\H^n$ and is in fact given by a multiple of a unitary map.
These maps are called conformal by Kor\'anyi and Reimann. Their rigidity theorem
states that all conformal maps belong to the
group~$\group{SU}(1,n)$.
\par
A second step was taken by P.~Pansu~\cite{Pu90}, who
proved  that in the quaternionic and octonionic
analogues of this set-up Liouville's theorem
holds under the sole assumption that the map in question preserves a suitable
contact structure of codimension greater than one. Similar phenomena
have been studied in more general situations: see,
for example, \cite{Bm01}, \cite{BH01}, \cite{GK98}, \cite{Gv87}.
\par
A remarkable piece of work concerning this circle of ideas is \cite{Ym93}, by
K.~Yamaguchi. His approach is at the infinitesimal level and is  based on
the theory of G structures, as developed by N.~Tanaka~\cite{Ta70}. The
crucial step in his analysis uses heavily Kostant's Lie algebra
cohomology and classification arguments.
 It is perhaps fair to say that the latest
important contribution in this area is the point of view adopted by Cowling,
De~Mari, Kor\'anyi and Reimann in \cite{CDKR1} and \cite{CDKR2}. As
mentioned earlier, they introduce the notion of multicontact mapping.
Their results have a non-trivial overlap with those by Yamaguchi, but are
independent of classification, rely on entirely elementary techniques
and focus on a very important issue: the main step in proving a rigidity
result at the Lie algebra level consists in showing that the appropriate
system of differential equations has polynomial solutions.
\vskip0.5truecm
One may reverse the point of view presented above and argue that
 a rigidity theorem exhibits G as a group
of geometric transformations (of some natural sort) of a homogeneous space of
G. It is then very natural  to ask if rigidity phenomena occur in a wider
variety of circumstances, for other Lie groups or, rather, for
submanifolds of a rigid manifold. In this thesis we  consider the case of
a large class of Hessenberg  submanifolds of G/P. Hessenberg manifolds
were introduced in the mid 80's by G.~Ammar and
C.~Martin~\cite{A},~\cite{AM} in connection with the study of the QR
algorithm as a dynamical system on flag manifolds. Let us briefly recall
the definition of the simplest Hessenberg manifolds.
\par
If $\group{G=SL}(n,\R)$ and P is the minimal parabolic subgroup of
G consisting of the unipotent upper-triangular matrices, then G/P
is identified with the flag manifold ${\rm Flag}(\R^n)$. The
elements of
 ${\rm Flag}(\R^n)$ are the nested sequences $S_1\subset S_2\subset\dots\subset S_{n-1}$
of linear subspaces of $\R^n$, with ${\rm dim}(S_i)=i$, and the
identification takes place by viewing the first $i$ columns of the matrix
$g\in \group{G}$ as a spanning set for $S_i$. Clearly, two matrices
will identify the same flag if and only if they differ by right
multiplication by an upper-triangular matrix with ones along the
main diagonal. This shows that
 ${\rm Flag}(\R^n)\simeq \group{G/P}$.
 \par
 Fix now a  matrix $A$ and a positive integer $p\in\{1,2,\dots,n-1\}$.
 Ammar and Martin say that the flag $S_1\subset S_2\subset\dots\subset
S_{n-1}$ is a Hessenberg flag of type $p$ relative to $A$ if
 $AS_i\subset S_{i+p}$ for all $i=1,\dots,n-p-1$. Thus, a Hessenberg flag of type $p$ is one for which $A$
shifts  a space of the flag into a larger space (one with $p$ additional
dimensions), within the same flag. The set of all Hessenberg flags of
type $p$ is denoted by
 ${\rm Hess}_p(A)$ and referred to as a Hessenberg manifold. It was proved in \cite{FDMth} that
if $A$ is diagonal and has distinct non-zero eigenvalues, then ${\rm
Hess}_p(A)$ is indeed a smooth manifold.
 \par
Notice that the
 defining condition may be formulated by the single matrix
equation $Ag=gR$, where $g\in \group{G}$ represents
 the flag and $R$ is any $n\times n$ matrix that has no more than $p$
non-zero sub-diagonals. Indeed, the first $i$ columns of $gR$  are in
this case a linear combination of the first $i+p$ columns of $g$. A
matrix like $R$ is known in the Numerical Analysis literature as a
Hessenberg matrix of type $p$. This clarifies the terminology.
 \par
Asking that the identity $Ag=gR$ is satisfied for some Hessenberg matrix of type $p$ is equivalent to saying
$g^{-1}Ag\in{\mathcal H}_p$, where ${\mathcal H}_p$ is the space of all Hessenberg matrices of type $p$, a space
that is stable under conjugation by elements in P. A simple but
far-reaching observation is that the defining equation can be written as
$\Ad(g^{-1})A\in{\mathcal H}_p$ and then interpreted in Lie-theoretical terms.
This leads to a general notion of
Hessenberg manifold. The definition makes sense whenever
${\mathcal H}_p$ is replaced by a vector space ${\mathcal H}$ in the Lie
algebra of G  that contains the Lie algebra of P and is stable under
$\Ad(\group{P})$.
 In the context of  semisimple Lie algebras, these P-modules ${\mathcal
H}$ can be described in terms of root spaces and turn out to be labeled
by those subsets ${\mathcal R}$ of the set of positive restricted roots
that satisfy the so-called Hessenbeg condition: if $\alpha\in{\mathcal
R}$ and $\beta$ is another positive restricted root such that
$\alpha-\beta$ is again a positive restricted root, then
$\alpha-\beta\in{\mathcal R}$. Thus, one defines a whole class of
manifolds. Each  element in the class depends on two choices: the Lie
algebra element $A$ and the combinatorial structure $\mc{R}$, whence the
standard notation ${\rm Hess}_{\mc{R}} (A)$. In particular, one recovers
the notion of ``type $p$''  described above by choosing ${\mathcal R}$ to
be the set of positive roots of height less than or equal to $p$.
\par
The Hessenberg manifolds were studied primarily by F.~De~Mari in a series of papers with different collaborators.
The basic topological and geometric features in the complex setting were studied with C.Procesi and M.~A.~Shayman in  \cite{DPS} and \cite{DS}, whereas the real Hessenberg manifolds and their close connection with the generalized Toda flow were considered in collaboration with M.~Pedroni \cite{DP}.
Recently, they have been investigated by J.~Tymoczko~\cite{Ty} from the point of view of combinatorial topology.
\vskip0.5truecm
The first problem addressed in this thesis is an appropriate
definition of multicontact structure on the Hessenberg manifold,
and it is considered in Chapter~\ref{multi3}. Here, as
in most of the existing literature, we consider ${\rm Hess}_{\mc{R}}
(H)$ when $H$ is a regular element  in the Cartan subspace of the Lie
algebra of G. The construction of the special sub-bundles requires a
careful local description of the manifold, and this is relatively
staightforward  once the  Bruhat decomposition is taken into play. The
point is that by means of the Bruhat decomposition, an open and dense
subset of G/P can be identified with the nilpotent group N occuring
in the Iwasawa decomposition G=KAN. Then the tangent space at the base
point, namely the identity $e\in$ N, is naturally identified with the Lie
algebra
$\Lie{n}$ of N and the exponential coordinates enable  to endow G/P
with a (local) structure governed by the restricted roots. By writing
down the equations that define ${\rm Hess}_{\mc{R}} (H)$, one realizes
that in these coordinates it is the graph of a polynomial mapping. Thus
one looks at the independent variables as a natural model.  In other
words, a coordinate subspace of the euclidean space N is selected as a
substitute of the Hessenberg manifold, that is a ``slice''
$\group{S}\subset{\rm Hess}_{\mc{R}} (H)$. Inside $\group{S}$, the
coordinates that correspond to simple roots are very well visible, and
this calls for the correct identification of the special bundles. In this
way the multicontact structure is proved to exist and to satisfy all the
reasonable properties that it should satisfy.
\par
Next, we present a case-study, namely we take $\group{G=SL}(4,\R)$ and
we take P to be the minimal parabolic subgroup of unipotent
lower-triangular matrices. We consider the structure $\mc{R}$ consisting
of all the positive roots except the highest one. Equivalently, we look
at the standard Hessenberg manifold with $p=2$. The usual realization of
a Cartan subspace $\Lie{a}$ is the space of traceless diagonal matrices,
and an element $H\in\Lie{a}$ is regular if and only if it has distinct
eigenvalues. The specific choice of $H$ is irrelevant, but we do make a
choice for simplicity. In the attempt of understanding what the
multicontact  mappings are in this case, it is very natural to follow the
method developed in \cite{CDKR2}, namely to look for the vector fields
whose flow consists of multicontact local diffeomorphisms.
\par
The study of multicontact vector fields leads immediately to a system
of differential equations (for the components of the vector field) that
reveals some of the basic principles that appear in the general case, but
also some special feature. First of all, the system can be seen as a
bunch of systems (in this case two), each associated to a maximal root in
$\mc{R}$, namely a root $\mu\in\mc{R}$ to which no simple root can be
added to give another root in $\mc{R}$. This is a general feature.
Secondly, each of the systems exhibits a hierarchic structure: once the
component along the maximal root $\mu$ is known, then all the components
labeled by the roots in the ``cone'' below $\mu$ are also known, by
suitable differentiation. This is also a general feature. Thirdly, each
subsystem is identical to the system whose solutions identify the
multicontact vector fields on some other 
$\tilde{\rm G}/ \tilde{\rm P}$ (here
$\tilde{\rm G}={\rm SL}(3,\R)$ and $\tilde{\rm P}$ is its minimal
 parabolic subgroup of unipotent lower-triangular matrices). Thus, by
\cite{CDKR1} each of them can be solved. This is a special feature, and
is what we refer to as an instance of several Iwasawa models paired
together. Finally, the systems overlap, and the core of the analysis
consists in understanding what are the consequences of the overlapping.
This is again a general problem.
\par
In the case at hand, that is $\group{G}=\group{SL}(4,\R)$, the analysis
can be carried out without difficulties and leads to a very interesting
answer. The Lie algebra of multicontact vector fields has dimension $9$
and is naturally isomorphic to a quotient Lie algebra. Let us describe it.
Denote by $\Sigma_+$ the set of  positive restricted roots and by ${\mc
C}$ the complement
${\mc C}=\Sigma_+\setminus\mc{R}$. Consider the vector space direct sums
$$
\Lie{n}=\sum_{\alpha\in\Sigma_+}\Lie{g}_\alpha,
\qquad
\Lie{n}_{\mc C}=\sum_{\alpha\in\mc C}\Lie{g}_\alpha,
$$
where evidently $\Lie{g}_\alpha$ is the root space associated to
$\alpha$ and $\Lie{n}$ is the usual Iwasawa nilpotent Lie subalgebra of
$\Lie{g}$. Because of the Hessenberg condition, $\Lie{n}_{\mc C}$ is an
ideal in $\Lie{n}$.
 Now, let $\Lie{q}$ denote the normalizer of $\Lie{n}_{\mc C}$ in $\Lie{g}$,
that is the largest subalgebra of $\Lie{g}$ in which $\Lie{n}_{\mc C}$
 is an ideal. In our case, $\mc{C}$ consists of the highest root
alone, and the Lie algebra of multicontact vector fields is  isomorphic
to the quotient $\Lie{q}/\Lie{n}_{\mc C}$.
\par
Motivated both by \cite{CDKR1}, \cite{CDKR2} and by the
previous example, we define the notion of multicontact vector field on
${\rm Hess}_{\mc{R}}(H)$. Since the nature of our investigations is
local, we focus on the local model of ${\rm Hess}_{\mc{R}}(H)$, that is
the slice $\group{S}$.  We ask ourselves the following basic question:
what is the structure of the Lie algebra  $MC\group{(S)}$ of
multicontact vector fields on $\group{S}$ in terms of the
combinatorial data $\mc{R}$? The remaining part of the thesis is
devoted to giving a partial answer to this question. We  explain
below the main steps and keep in mind the case
$\group{G=SL}(n,\R)$.
\par
A simple look at the combinatorics of $\mc{R}$ suggests to partition
it in what one is naturally inclined to think of as a connected component.
This is what we shall call ``dark zones''. The reason for this
terminology is that each of these components may be viewed as the
(overlapping) union of shadows, each of which stems from a maximal root
$\mu$. In the $\Lie{sl}(n,\R)$ case, the  picture explains the
wording.
\vskip6.4truecm
\setlength{\unitlength}{1cm}
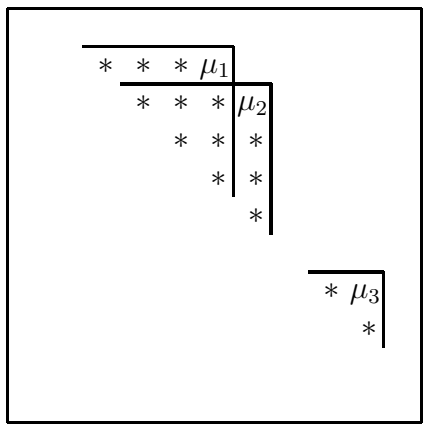
\begin{picture}(0,0)\thicklines
\put(3,1){\line(0,1){5.5}}
\put(3,1){\line(1,0){5.5}}
\put(3,6.5){\line(1,0){5.5}}
\put(8.5,1){\line(0,1){5.5}}
\put(4,6){\line(1,0){2}}
\put(6,4){\line(0,1){2}}
\put(4.5,5.5){\line(1,0){2}}%
\put(6.5,3.5){\line(0,1){2}}%
\put(7.0,3){\line(1,0){1}}%
\put(8.0,2){\line(0,1){1}}%
\put(4.2,5.65){$*$}
\put(4.7,5.65){$*$}
\put(5.2,5.65){$*$}
\put(5.55,5.65){$\mu_1$}
\put(4.7,5.15){$*$}
\put(5.2,5.15){$*$}
\put(5.7,5.15){$*$}
\put(6.05,5.15){$\mu_2$}
\put(5.2,4.65){$*$}
\put(5.7,4.65){$*$}
\put(6.2,4.65){$*$}
\put(5.7,4.15){$*$}
\put(6.2,4.15){$*$}
\put(6.2,3.65){$*$}
\put(7.2,2.65){$*$}
\put(7.55,2.65){$\mu_3$}
\put(7.7,2.15){$*$}
\end{picture}
\vskip-0.5truecm
\par
\noindent
These partitions of course reflect the nature of the differential
equations that the components of a multicontact vector field must satisfy.
First of all, a vector field has components labeled by $\mc{R}$;
secondly, the dark zones correspond to
independent subsystems and, thirdly, each shadow exhibits a
hierarchic structure.
As a consequence of the mutual independence of
dark zones, we may prove a reduction theorem, Theorem~\ref{reduction},
that unables us to safely assume that $\mc{R}$ is a single dark zone.
This really means that all the simple roots are in
$\mc{R}$. From now on we thus work under this assumption.
\vskip6truecm
\setlength{\unitlength}{1cm}
\hskip1truecm
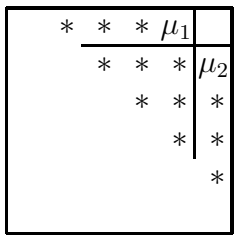
\begin{picture}(0,0)\thicklines
\put(3.5,6){\line(1,0){3}}
\put(6,4){\line(0,1){2}}
\put(4.5,5.5){\line(1,0){2.0}}%
\put(6.5,3){\line(0,1){3}}%
\put(3.5,3){\line(0,1){3}}%
\put(3.5,3.0){\line(1,0){3}}%
\put(4.2,5.65){$*$}
\put(4.7,5.65){$*$}
\put(5.2,5.65){$*$}
\put(5.55,5.65){$\mu_1$}
\put(4.7,5.15){$*$}
\put(5.2,5.15){$*$}
\put(5.7,5.15){$*$}
\put(6.05,5.15){$\mu_2$}
\put(5.2,4.65){$*$}
\put(5.7,4.65){$*$}
\put(6.2,4.65){$*$}
\put(5.7,4.15){$*$}
\put(6.2,4.15){$*$}
\put(6.2,3.65){$*$}
\end{picture}
\vskip-2.5truecm
\noindent
At this point, the hierarchy plays a key role.
If the components along the
maximal roots satisfy
the appropriate differential equations, then all the components in the
shadow below them are obtained by differentiation.
Therefore, the problem reduces to analyzing the differential
equations for the maximal roots  on one hand, and of understanding the
implications of the overlapping shadows on the other hand.
By doing so, we obtain a first result on  the Lie algebra 
 $MC\group{(S)}$
of the multicontact vector fields, namely that $MC\group{(S)}$
contains
$\Lie{q}/\Lie{n}_{\mc{C}}$. This latter Lie algebra corresponds to filling
the previous picture with the stars  that label the normalizer  $\Lie{q}$,
and then taking the quotient modulo the ideal
$\Lie{n}_{\mc{C}}$, that is the black dot.
\vskip6.5truecm
\setlength{\unitlength}{1cm}
\hskip-1truecm
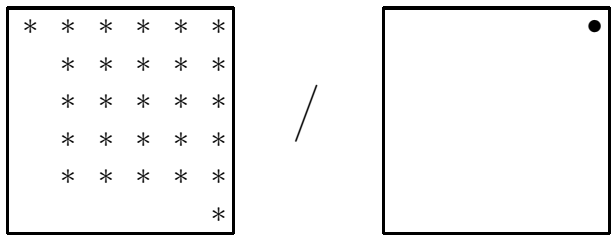
\begin{picture}(0,0)\thicklines
\put(3.5,6){\line(1,0){3}}
\put(6.5,3){\line(0,1){3}}%
\put(3.5,3){\line(0,1){3}}%
\put(3.5,3.0){\line(1,0){3}}%
\put(3.7,5.65){$*$}
\put(4.2,5.65){$*$}
\put(4.7,5.65){$*$}
\put(5.2,5.65){$*$}
\put(5.7,5.65){$*$}
\put(6.2,5.65){$*$}
\put(4.2,5.15){$*$}
\put(4.7,5.15){$*$}
\put(5.2,5.15){$*$}
\put(5.7,5.15){$*$}
\put(6.2,5.15){$*$}
\put(4.2,4.65){$*$}
\put(4.7,4.65){$*$}
\put(5.2,4.65){$*$}
\put(5.7,4.65){$*$}
\put(6.2,4.65){$*$}
\put(4.2,4.15){$*$}
\put(4.7,4.15){$*$}
\put(5.2,4.15){$*$}
\put(5.7,4.15){$*$}
\put(6.2,4.15){$*$}
\put(4.2,3.65){$*$}
\put(4.7,3.65){$*$}
\put(5.2,3.65){$*$}
\put(5.7,3.65){$*$}
\put(6.2,3.65){$*$}
\put(6.2,3.15){$*$}
\put(7.3,4.5){$\Big/$}
\put(8.5,6){\line(1,0){3}}
\put(11.5,3){\line(0,1){3}}%
\put(8.5,3){\line(0,1){3}}%
\put(8.5,3.0){\line(1,0){3}}%
\put(11.2,5.65){$\bullet$}
\end{picture}
\vskip-2truecm
\noindent
In the general case, however, the converse inclusion
$MC\group{(S)}\subseteq\Lie{q}/\Lie{n}_{\mc{C}}$ does not hold, as
explained in detail in Section~\ref{counterex}.
\par
A natural assumption
under which the converse does hold, is that each shadow defines a
subalgebra, necessarily an Iwasawa algebra in its own right. One may
use Theorem~4.1 in~\cite{CDKR2} for each single shadow. This is done in
the final step, Theorem~\ref{al}, where we also glue all the different
pieces together. This requires a technical description
(Lemma~\ref{norma}) of the normalizer $\Lie{q}$ in terms of roots.
We finally draw some conclusions concerning the group of multicontact
diffeomorphisms. In Proposition~\ref{mcg} we show that it contains
canonically the quotient
$\group{Q}/\group{N}_{\mc{C}}$, where $\group{Q}=\Int(\Lie{q})$ and
$\group{N}_{\mc{C}}=\exp\Lie{n}_{\mc{C}}$.
\vskip0.5truecm
This thesis also contains a chapter devoted to some decomposition results
for the polynomials that generate the Lie algebra of multicontact vector
fields on G/P, under the further assumption that $\Lie{g}$ is a split
form. We believe that this is of some independent interest and that it
indicates another possible area of investigation, namely the explicit
description of all the (special) polynomial algebras that the theory of
multicontact vector fields seems to produce. For clarity and internal
consistency, this part  actually precedes the study of Hessenberg
manifolds and is developed in Chapter~\ref{po}. Finally,
Chapters~\ref{pr} and~\ref{hesse} collect some prerequisites.

\chapter{Preliminaries}~\label{pr}
\setcounter{page}{1}
\pagenumbering{arabic}
\vskip0.5cm
In this chapter we introduce the fundamental tools that are
used in this thesis. In particular, in the first section we recall some
very well-known facts about simple Lie algebras, and we fix the  notations
that will be used throughout. After that, we shall discuss some recent
results obtained by Cowling, De Mari, Koranyi and Reimann on the contact
structures generalized to the boundaries of symmetric spaces of the type
G/P. In the third and last section we illustrate these results in one
example.
\section{Simple Lie algebras}\label{sla}
We shall work with real simple Lie algebras, although
most of what we do holds, {\it mutatis mutandis}, for
semisimple Lie algebras.
For the reader's convenience, we collect some facts about simple
Lie algebras and their decompositions. For more details the reader
can look at~~\cite{B},~\cite{KN02},~\cite{V}.
\par
Let $\Lie{g}$ be a simple Lie
algebra with  Cartan involution $\theta$. Let $\Lie{k} \oplus
\Lie{p}$ be the Cartan decomposition of $\Lie{g}$, where
$\Lie{k} = \{ X \in \Lie{g} : \theta X=X\}$
  and
$\Lie{p} = \{ X \in \Lie{g} : \theta X=-X \}.$
Fix a maximal abelian subspace $\Lie{a}$ of $\Lie{p}$. The
dimension of $\Lie{a}$ is an invariant of $\Lie{g}$ and is
called the real rank of $\Lie{g}$. Denote by $\Lie{a}^\prime$
the dual of $\Lie{a}$. For $\alpha \in \Lie{a}^\prime$, set
\begin{equation}
\Lie{g}_\alpha
=\{ X \in \Lie{g} :[H,X] = \alpha
(H)X,\,\forall H\in\Lie{a}\}.
\label{rs}\end{equation}
When $\alpha \neq 0$ and $\Lie{g}_\alpha$ is not trivial,
$\alpha$ is said to be a restricted\footnote{The reason of
the adjective ``restricted'' resides in the fact that
restricted roots arise as restrictions to $\Lie{a}$
of the roots relative to the Cartan subalgebra
$\Lie{h}=(\Lie{a}+\Lie{t})^c$ of the complexification of
$\Lie{g}$, where $\Lie{t}$ is a maximal abelian subspace of $\Lie{m}$,
as defined in~\eqref{emme}.}  root of
$\Lie{g}$.  Denote by $\Sigma = \Sigma (\Lie{g} , \Lie{a})$
the set of restricted roots of $\Lie{g}$.  It satisfies
all the axioms of a (not necessarily reduced) root system
in the usual sense \cite{B}; we refer to it as the root
system of $\Lie{g}$. The basic decomposition of $\Lie{g}$ is then
the so-called restricted root space decomposition
\begin{equation}
\Lie{g}=
\Lie{m}\oplus\Lie{a}\oplus
\bigoplus_{\alpha\in\Sigma}\Lie{g}_\alpha,
\label{rsd}
\end{equation}
where
$\Lie{m}$ is the centralizer of $\Lie{a}$ in $\Lie{k}$, that is
\begin{equation}
\Lie{m}=\left\{X\in\Lie{k}:[X,H]=0,\,H\in\Lie{a}\right\}.
\label{emme}\end{equation}
Observe that $\Lie{m}\oplus\Lie{a}=\Lie{g}_0$, the space that  corresponds to the
choice $\alpha=0$ in \eqref{rs}. Also, we remark that the direct sums  in~\eqref{rsd}
only concern the vector space structure.
\par
Fix a partial ordering $\succ$ on
$\Sigma$ and denote by $\Sigma _+$ the subset of
$\Lie{a}^\prime$ of positive restricted roots. The space
$\Lie{a}^\prime$ is endowed with the inner product $(\cdot
, \cdot )$ induced by the Killing form $B$ of $\Lie{g}$.
It is defined by $(\alpha ,
\beta) = B(H_\alpha , H_\beta)$, where $H_\alpha$ is the
element of $\Lie{a}$ that represents the functional $\alpha$ via the  Killing
form, that is
$\alpha (H) = B(H, H_\alpha)$ for all $H \in \Lie{a}$.
  From the Jacobi identity one
immediately gets that $[\Lie{g}_\alpha , \Lie{g}_\beta ]
\subseteq \Lie{g}_{\alpha + \beta}$, provided that $\alpha + \beta$
is a root.
Choose a basis for each restricted root space. The above version of the
Jacobi identity allows us to define the structure constants
$c_{\alpha,\beta}$, that is the real numbers that satisfy
$[X_\alpha , X_\beta ] = c_{\alpha, \beta} X_{\alpha + \beta}$.
\par
A positive root $\alpha$ is called simple if it cannot be written as
a sum of positive roots. If
$\Delta = \{\delta_1 , \ldots , \delta_r \}$ denotes the set of simple
roots, then the cardinality $r$ of
$\Delta$ is equal to the real rank of $\Lie{g}$. Every $\alpha \in
\Sigma_+$ is a linear combination of elements of $\Delta$ with
coefficients in $\mathbb{N} \cup \{0\}$.
Thus every positive root $\alpha$ can be written as
$\alpha =\sum_{i=1}^{r} n_i \delta_i$ for uniquely defined non-negative
integers $n_1,\dots,n_r$, and the positive integer
${\rm ht}(\alpha) =\sum_{i=1}^{r} n_i$ is called the height of
$\alpha$. It is well-known that there is exactly one root
$\omega$, called the highest root, that satisfies $\omega \succ \alpha$
(strictly) for every other root~$\alpha$.
\par
Of central importance in the present context are the nilpotent Lie
algebra
\begin{equation}\Lie{n}=
\bigoplus_{\gamma \in \Sigma_+} \Lie{g}_\gamma
\label{enne}\end{equation}
and its counterpart $\overline{\Lie{n}}=\theta(\Lie{n})$. For
instance $\Lie{n}$ appears in one of the most useful features in the  theory of
semisimple Lie algebras, the Iwasawa decomposition of
$\Lie{g}$, namely
$\Lie{g}=\Lie{k} \oplus \Lie{a} \oplus \Lie{n}$. We shall thus refer
to $\Lie{n}$ as the nilpotent Iwasawa subalgebra of $\Lie{g}$.
\par
One of the main ingredients in the theory of multicontact
mappings as developed in \cite{CDKR2} is the nature of the index set
$\Sigma_+$ that labels the direct sum~\eqref{enne} or, more
importantly, of the corresponding direct summands $\Lie{g}_\alpha$
as $\alpha$ runs in $\Sigma_+$. Indeed, it provides what we call a
{\it multistratification}, as we now briefly explain. The notion
of stratified Lie algebra refers properly to the strata
\begin{equation}
\Lie{n}_i=\bigoplus_{{\rm ht}(\gamma)=i} \Lie{g}_\gamma,
\qquad i=1,\dots,{\rm ht}(\omega).
\label{protagonista}\end{equation}
because they satisfy
$$
[\Lie{n}_i,\Lie{n}_j]\subset\Lie{n}_{i+j},
$$
as a further application of the Jacobi identity shows. Since
 the ``ground'' stratum $\Lie{n}_1$ is the direct sum of the root
spaces that are labeled by the simple roots and since each positive
restricted root is the sum of simple ones, it follows that
$\Lie{n}_1$ generates $\Lie{n}$ as a Lie algebra. Thus
$\Lie{n}$ is a stratified Lie
algebra in the usual sense (see~\cite{FS}). What is more important is  the fact
that $\Lie{n}_1$, and every higher stratum, is a direct sum of
finer building blocks, as indicated in~\eqref{protagonista}. The way in  which
these finer blocks (i.e. the restricted root spaces) behave under  bracket is
governed by the root system, which has a highly non-trivial structure.  This is the
multistratification.

\section{Multicontact mappings on G/P}\label{mm}
In this section we recall some results contained in~\cite{CDKR2}.
In particular, we quote a  theorem that plays a fundamental
role in this thesis and then we describe the steps that
lead to the proof of it.
\par
Let $\Lie{g}$ be as above and let $\group{G}$ be a Lie group whose
Lie algebra is $\Lie{g}$. Let
$\group{P}$ be a minimal parabolic subgroup of G. We may assume
that the center of
$\group{G}$ is trivial. Indeed, if $\group{Z}$ is the center of
$\group{G}$, then $\group{Z}\subset \group{P}$, and so G/P and (G/Z)/(P/Z)
may be identified. Moreover, the action of G on G/P factors to an
action of G/Z.
\par
Among all groups with trivial centers and the same
Lie algebra $\Lie{g}$, the largest is the group $\Aut(\Lie{g})$
of all automorphisms of $\Lie{g}$, and the smallest is the group
$\Int(\Lie{g})$ of the inner automorphisms of $\Lie{g}$, the
connected component of the identity of $\Aut(\Lie{g})$. Any group
$\group{G}_1$ such that
$\Int(\Lie{g})\subseteq\group{G}_1\subseteq\Aut(\Lie{g})$,
with corresponding minimal parabolic subgroup $\group{P}_1$,
gives rise to the same space, meaning that $\group{G}_1/\group{P}_1$
may be identified with $\Aut(\Lie{g})/\group{P}$ if P is a minimal
parabolic subgroup of $\Aut(\Lie{g})$. For the purposes of this
thesis the correct assumption is that G is connected and centerless,
and hence we can assume  $\group{G}=\Int(\Lie{g})$
and that P is a  minimal parabolic subgroup of G.
\par
The most natural choice for P is as follows.
According to the notation introduced in Section~\ref{sla}, let
$\Lie{m}$ denote the centralizer of $\Lie{a}$ in $\Lie{k}$ as  previously defined.
As easily verified, $\Lie{m}$ is the Lie algebra of the centralizer
of $\Lie{a}$ in K, the latter being the connected Lie subgroup of G
whose Lie algebra  is $\Lie{k}$. In other words
$$
\group{M}=\left\{m\in\group{K}:\Ad m(H)=H,\,H\in\Lie{a}\right\}.
$$
Let now $\group{A}=\exp\Lie{a}$ and $\group{N}=\exp\Lie{n}$ denote
the connected (and simply connected) Lie subgroups of G with Lie
algebras
$\Lie{a}$ and
$\Lie{n}$, respectively. Finally,  put
$\overline{\group{N}}=\exp \overline{\Lie{n}}$. Then
$\group{P}=\group{MA}\overline{\group{N}}$ is the minimal
parabolic subgroup of G that we shall be concerned with, and the
latter expression is the Langlands decomposition of it.
\par
By means of the
Bruhat decomposition (see
\cite{KN02},Ch.VII, Sec.4) the group $\group{N}$ may be seen as open
and dense in
$\group{G}/\group{P}$. Indeed, if we denote by $b$ the base point in
G/P (that is, the identity coset), the Bruhat lemma states that the
mapping
$\psi:\group{N}\to \group{G/P}$ defined by $\psi(n)=nb$ is injective and its
image is dense and open. The differential $\psi_*$ then maps $\Lie{n}$,
the tangent space to N at the identity $e$, onto $T_b$, the
tangent space to G/P at the
base point. When $\delta$ is a simple restricted root, we denote by
$S_{\delta,b}$ the subspace $\psi_*(\Lie{g}_\delta)$ of $T_b$.
In Lemma~2.2  of~\cite{CDKR2} it is shown that the action of any
element $p\in\group{P}$ on G/P induces an action $p_*$ on the
tangent space $T_b$ which in turn induces an action
$\psi^{-1}_*p_*\psi_*$ on $\Lie{n}$. This last action preserves all
the spaces $\Lie{g}_\delta$ for simple $\delta$. This lemma
allows us to identify  $\Lie{n}$ with the tangent space $T_x$
at any point $x$ in G/P, and to identify the subspaces
$\Lie{g}_\delta$ of $\Lie{n}$ with subspaces $S_{\delta,x}$ of $T_x$.
Indeed we may write $x$ as $gb$, where $g\in\group{G}$; then
the images $g_*\psi_*\Lie{g}_\delta$ are well defined, and
independent of
the representative $g$ of the coset, although the identification
$\Lie{g}_\delta\to S_{\delta,x}$
does depend on the representative. Since we never make use of the
explicit identification, we shall always write $\Lie{g}_\delta$ in
place of $S_{\delta,x}$.
\par
Consider a diffeomorphism $f$ between open
subsets of $\group{G}/\group{P}$, say $\mc{U}$ and $\mc{V}$.
By density, we can assume that $\mc{U}$ and $\mc{V}$ are subsets of
$\group{N}$. Then, $f$ is called a
\emph{multicontact} map if
$f_*$ maps $\Lie{g}_\delta$ in itself, for every simple root
$\delta$. In fact, the original definition given in~\cite{CDKR2}
is slightly weaker, and is designed to cover a wider class of
situations, essentially allowing some disconnectedness of G.
The authors of~\cite{CDKR2} allow $f_*$ to permute the various
$\Lie{g}_\delta$ for simple $\delta$.
In their context, with $\group{G}=\Aut(\Lie{g})$, they prove that  every
multicontact mapping on
$\group{G}/\group{P}$ is the restriction of the action of a uniquely
determined element
$g\in\group{G}$.
The crucial point of their proof is to focus on the infinitesimal  analogue
of the notion of multicontact map, where no permutation comes into play.
In our setting, this latter notion is even more important, and we  actually
take it as the basic notion. For this reason we recall it in
full detail.
\par
The starting point is to consider multicontact vector fields, that is,
vector  fields $V$ on $\mc{U}$ whose local flow
$\{\phi^V_t \}$ consists of multicontact maps.
If $X_\delta \in \Lie{g}_\delta$ and $\delta$ is a simple
root, then
$$
\frac{d}{dt} (\phi^V_t)_{*}(X_\delta)\Bigr|_{t=0} =
-\mathcal{L}_V(X_\delta)=[X_\delta,V],
$$
where $\mathcal{L}$ denotes the Lie derivative. Thus a smooth
vector field $V$ on $\mc{U}$ is a
multicontact vector field if and only if
\begin{equation}\label{mceq}
[V, \Lie{g}_\delta] \subseteq \Lie{g}_\delta
\qquad\text{for every simple root $\delta$.}
\end{equation}
This is to be interpreted  in the sense that, if
$Y\in\Lie{g}_\delta$, then $[V,Y]$ is a section of the subbundle
$\Lie{g}_\delta$ of
$T\mc{U}$, as explained above.
\par
Define a representation
$\tau$ of the Lie algebra $\Lie{g}$
as a set of vector fields on N as follows:
$$
(\tau (X) f) (n) = \frac{d}{dt} f ([\exp(tX)  n])
\Bigr|_{t=0}.
$$
Here $[\exp(tX) n]$ denotes the action of G on G/P, restricted to an
action on N. Hence
$[\exp(tX)  n]$ is the N-component of the product $\exp(tX)\cdot n$ in
the  Bruhat decomposition of G/P.
The infinitesimal analogue of the main result contained in
\cite{CDKR2} is Theorem~4.1, that we
recall here.
\begin{theorem}{\rm (}\cite{CDKR2}{\rm )}\label{cdkr}
Suppose that $\Lie{g}$ has real rank at least two. Then every $C^1$
multicontact vector field is in fact
smooth, and the Lie algebra of multicontact vector fields on $\mc{U}$
consists of the restrictions of
$\tau(\Lie{g})$ to $\mc{U}$.
\end{theorem}
The action of G on N is multicontact by
construction. Hence, $\tau(\Lie{g})$ gives multicontact vector fields
by definition.
Theorem~\ref{cdkr} is thus about the converse implication, which is
proved in several steps. Although in the end the basic assertions that
lead to the proof are (essentially) coordinate-free, they
are best formulated in terms of a chosen canonical basis, as we  describe next.
For every
$\alpha\in\Sigma_+$,  let $m_\alpha$ denote the dimension of
$\Lie{g}_\alpha$ and fix a basis
$\{X_{\alpha,i}:\alpha\in\Sigma_+,i=1,\dots,m_\alpha\}$  of $\Lie{n}$
consisting of
left--invariant vector fields on N. Thus, a smooth vector field $V$ on
$\mc{U}$  can be written as
\begin{equation}
V =
\sum_{\alpha \in \Sigma_+}\sum_{i=1}^{m_\alpha} v_{\alpha,i}  X_{\alpha,i},
\label{canonica}\end{equation}
for some smooth functions $v_{\alpha,i}$.
The main steps in the proof of Theorem~\ref{cdkr} are summarized in the
following statements:
\begin{itemize}
\item[(i)] a multicontact vector field is determined by its
component in the direction of the highest
root, namely $\{v_{\omega,i}:i=1,\dots,m_\omega\}$, and these
functions are determined by a
particular set of differential equations;
\item[(ii)] the differential equations imply that any $v_{\omega,i}$
is a polynomial function in
canonical coordinates;
\item[(iii)] the polynomial nature of the functions $v_{\omega,i}$
implies
Theorem~\ref{cdkr} by general homogeneity arguments.
\end{itemize}

The first two steps tell us that any multicontact vector field
corresponds to a vector
$(v_{\omega,1},\dots,v_{\omega,m_{\omega}})$ of polynomials. Steps (i)  and (ii) is
where most of the hard work goes. It involves a careful analysis of the
system of differential equations via a detailed study of several  properties
of the restricted root system.
\par
We next discuss in detail the homogeneity concept that is used in the  argument of
step (iii).
\par
We select an
element $H_0$ in the Cartan subspace $\Lie{a}$ such that $\delta  (H_0)=-1$
for all simple roots $\delta$ (this is possible because the Cartan  matrix
is non-singular~\cite{KN02}). A function
$v$ on
$\group{N}$ is said to be homogeneous of degree $r$ if it does not  vanish
identically and if it also satisfies
${\tau(H_0)} v=rv$. A vector field $V$ is said to be
homogeneous of degree $s$ if it does not vanish identically and it
satisfies $[{\tau(H_0)},V]=sV$. Hence
\begin{align*}
{\rm deg}(vV) &= {\rm deg}(V) + {\rm deg}(v),\\
{\rm deg}(V(v)) &= {\rm deg}(v) + {\rm deg}(V) \hskip0.5cm \text{(except
when $V(v)=0$)},\\
{\rm deg}([V,W]) &= {\rm deg}(V) + {\rm deg}(V) \hskip0.4cm
\text{(except when $V$ and
$W$ commute)}.
\end{align*}
  All left invariant vector fields $X$, where $X$ is in the stratum
$\Lie{n}_j$, are of degree
$-j$.
Indeed, write $n_s=\exp(-sH_0)n\exp(sH_0)$.
Then
\begin{align*}
[\tau(H_0),X]f(n)&=\frac{d}{ds}\frac{d}{dt} \left(
f([\exp(-sH_0)n]\exp(tE))\right.\\
&\left. \hskip1.5cm -f(\exp(-sH_0)[n\exp(tE)])\right)\Big|_{t=s=0}\\
&=\frac{d}{ds}\frac{d}{dt} \left(
f(\exp(-sH_0)n\exp(sH_0)\exp(tE))\right.\\
&\left. \hskip1.5cm
-f(\exp(-sH_0)n\exp(tE)\exp(sH_0))\right)\Big|_{t=s=0}\\
&=\frac{d}{ds}\frac{d}{dt} \left(
f(n_s\exp(tE))-f(n_s\exp(e^{js}tE))\right)\Big|_{t=s=0}\\
&=\frac{d}{ds}\left(Xf(n_s)-e^{js}
Xf(n_s)\right)\Big|_{s=0}\\ &=-j Xf(n).
\end{align*}
If $X$ lies in a root space $\Lie{g}_\alpha$, we have
$$
[\tau(H_0),\tau(X)]=\tau([H_0,X])=\alpha(H_0)\tau(X)=-k\tau(X),
$$
where $k$ is the height of $\alpha$.
Thus, all homogeneous vector fields in $\tau(\Lie{g})$ have degree
between $-h$ and $h$, where $h$ is
the height of $\Lie{n}$, that is,
the lenght of its stratification.

\par
Let us go back to step (iii) in the proof of Theorem~\ref{cdkr}.
First, one observes that a multicontact vector field $V$ can be written  as the sum of
its homogeneous parts;
togheter with $V$, these also satisfy the differential equations that  define
the multicontact conditions, i.e., they are multicontact (see Section
3.2 in \cite{CDKR2}). Thus,  $V$ can be assumed to be homogeneous. The  proof proceeds
by treating separately negative, positive and zero degrees. We  reproduce below this
final argument.
\par
{\bf Negative degree.} Fix $Y$ in $\Lie{n}_1$ and assume that ${\rm  deg}V<0$. Then
$$
{\rm deg}([Y,V])={\rm deg}(Y) + {\rm deg}(V)=-1 +{\rm deg}(V)<-1,
$$
so $[Y,V]\in\Lie{n}_1$ cannot be a section of the subbundle of $T\mc{U}$
associated to $\Lie{n}_1$ unless $[Y,V]=0$. Therefore $[Y,V]=0$ for all
$Y\in\Lie{n}_1$, and since $\Lie{n}_1$ generates $\Lie{n}$, it follows
that $[Y,V]=0$ for all $Y\in\Lie{n}$. If a vector field commutes with
infinitesimal right translations, then it is an infinitesimal left
translation. Consequently, $V=\tau(X)$ for some $X$ in $\Lie{n}$.
\par
{\bf Negative degree.}
Denote by $\group{M}^\prime$ the normalizer in K of $\Lie{a}$, namely
$$
\group{M}^\prime=\left\{ k\in\group{K}:\Ad k (H)\in \Lie{a},\,
H\in\Lie{a}\right\},
$$
and let $W$ be the Weyl group $\group{M}^\prime /\group{M}$
(see~\cite{B}).
To treat the case where ${\rm deg}(V)>0$, we consider
the inversion map
$s$ on N, induced by the action of $m^\prime$ on G/P, where $m^\prime\in
\group{M}^\prime$ is a representative of the longest Weyl group element.
The induced map $s_*$ has the property that ${\rm deg}(s_* V)=-{\rm
deg}(V)$ for homogeneous vector fields $V$. Consequently, if $V$ is a
homogeneous multicontact vector field of positive degree, then $s_* V$  is
multicontact (because $s$ is multicontact) and hence polynomial. Now  $s_*
V$ is of negative degree, so $s_*V=\tau(X)$ for some $X$ in $\Lie{n}$,
whence $V=\tau(s_*^{-1} X)$.
\par
{\bf Zero degree.}
Finally, if ${\rm deg}(V)=0$, then $\ad V$ preserves both  $\tau(\Lie{n})$
and $\tau(\overline{\Lie{n}})$, that is, $\tau(\Lie{g})$.
As $\ad V$ is a derivation of the semisimple Lie algebra  $\tau(\Lie{g})$,
there exists $Y$ in $\Lie{g}$ such that $V-\tau(Y)$ commutes with
$\tau(\Lie{g})$. If a vector field commutes with $\tau(\Lie{g})$, then  in
particular it commutes with infinitesimal left translations, so it is an
infinitesimal right translation, and since it also commutes with
dilations, it is zero. Hence $V=\tau(Y)$.
This  conludes the proof ot Theorem~\ref{cdkr}.
\vskip0.5truecm

\section{A case study: $\Lie{sl}(3,\mathbb{R})$}\label{sl3}
We report an example contained in~\cite{CDKR1}. We describe this case  study
in order to illustrate the
techniques and the results we presented in the previous section.
Moreover, this basic example is crucial for our aims. Indeed, we
shall use the results we report here
in Chap.~\ref{multi3}.
\par
  Let $\group{G}=\group{SL}(3,\mathbb{R})$,
and let P be the minimal parabolic subgroup of G of lower triangular
matrices. For $x$, $y$ and $u$ in
$\mathbb{R}$, denote by
$\nu(x,y,u)$ the matrix
$$
\bmatrix
0&x&u\\
0&0&y\\
0&0&0
\endbmatrix.
$$
Take $\alpha$ and $\beta$ to be the simple roots relative to the
standard Cartan subalgebra of $\Lie{sl}(3,\mathbb{R})$ of diagonal
matrices: $\alpha({\rm diag}(a,b,c))=(a-b)$ and $\beta(({\rm
diag}(a,b,c))=(b-c)$. Then
\begin{align*}
\Lie{g}_\alpha&=\{\nu(x,0,0):x\in\mathbb{R}\},\\
\Lie{g}_\beta &=\{\nu(0,y,0):y\in\mathbb{R}\},\\
\Lie{g}_{\alpha +\beta}&=\{\nu(0,0,u):u\in\mathbb{R}\}.
\end{align*}
Further, $\Lie{n}=\{\nu(x,y,u):x,y,u\in\mathbb{R}\}$;
the algebra $\Lie{n}$ has the multistratification $\Lie{g}_\alpha
\oplus \Lie{g}_\beta
\oplus\Lie{g}_{\alpha +\beta}$ and the stratification
$\Lie{n}_1\oplus\Lie{n}_2$, where $\Lie{n}_1 =
\Lie{g}_\alpha \oplus \Lie{g}_\beta$ and $\Lie{n}_2 =\Lie{g}_{\alpha  +\beta}$.
\par
We write a vector field $V$ on $\mc{U}$ (open subset of N) as $fX +
gY +hU$, where $f$, $g$ and $h$ are smooth functions on $\mc{U}$ in
the coordinates $x$, $y$, $u$ and where $\{X,Y,U\}$ is the canonical  basis
of $\Lie{n}$. Viewed as left--invariant vector fields,  they are
$$
X=\frac{\partial}{\partial x},\qquad Y=\frac{\partial}{\partial y}
+x\frac{\partial}{\partial u},\qquad U=\frac{\partial}{\partial u}.
$$
Clearly
$$
[X,Y]=U
$$ is the only non-zero bracket.
We ask ourselves when $V$ is a multicontact vector field. The
multicontact vector field
equations~\eqref{mceq} state that this happens if and only if  $[V,X]=\lambda X$ and
$[V,Y]=\mu  Y$, for some smooth functions
$\lambda$ and $\mu$ on $\mc{U}$. These equations imply immediately that
\begin{align*}
\lambda X=-gU-(Xf)X-(Xg)Y-(Xh)U\\
\mu Y=fU-(Yf)X-(Yg)Y-(Yh)U,
\end{align*}
which in turn imply that
\begin{align*}
Xf&=-\lambda  &  Yg&=-\mu\\
Xg&=0 & Yf&=0\\
Xh&=-g & Yh&=f.
\end{align*}
We see at once that $f$ and $g$ are determined by $h$ and that $h$  itself
satisfies the differential equations
\begin{equation}\label{sleq}
X^2h=Y^2h=0.
\end{equation}
The equation $X^2 h=\partial^2 h/\partial x^2=0$ has the general  solution
$$
h(x,y,u)=h_0(y,u) + xh_1(y,u),
$$
for some functions $h_0$ and $h_1$.
The equation $Y^2 h=0$ then becomes
\begin{align*}
0&=\big( \frac{\partial^2}{\partial y^2} +2x
\frac{\partial^2}{\partial y\partial u} +x^2
\frac{\partial^2}{\partial u^2} \big) (h_0 + xh_1)\\
&=x^3\big( \frac{\partial^2 h_1}{\partial u^2} \big)
+x^2\big(\frac{\partial^2 h_0}{\partial u^2} +2\frac{\partial^2
h_1}{\partial y\partial u}\big)
+x\big(\frac{\partial^2 h_1}{\partial y^2} +2\frac{\partial^2
h_0}{\partial y\partial u}\big)
+\big(\frac{\partial^2 h_0}{\partial y^2}\big).
\end{align*}
Since the right hand side vanishes identically in some open set, the
coefficients of the various
powers of $x$ must vanish. Considering the $x^3$ term, we see that
$\partial^2 h_1/\partial u^2 =0$.
Differentiating the coefficient of the $x^2$ term once with respect
to $u$, we deduce that
$\partial^3 h_0/\partial u^3 =0$. Next, considering the constant term
yields $\partial^2 h_0/\partial
y^2 =0$, and then differentiating the coefficient of the $x$ term
once with respect to $y$, we deduce
that $\partial^3 h_1/\partial y^3 =0$. Summarizing, we have shown that
$$
\frac{\partial^3}{\partial u^3} h_0=0,\qquad\frac{\partial^2}{\partial
u^2}h_1=0,\qquad\frac{\partial^2}{\partial y^2} h_0
=0,\qquad\frac{\partial^3}{\partial y^3}h_1 =0.
$$
The first two equations imply that
$$
h_0 (y,u)=u^2 a(y) +ub(y)+c(y),\qquad h_1(y,u)=ud(y)+e(y),
$$
and the second two equations then imply that $a^{\prime
\prime}=b^{\prime \prime}=c^{\prime
\prime}=d^{\prime \prime\prime}=e^{\prime \prime\prime}=0$, so that
both $h_0$ and $h_1$ are
polynomials, whence $h$ is too.
The calculations we did so far correspond to the steps (i) and (ii)
of the proof of
Theorem~\ref{cdkr}, and the techniques here presented can be
generalized. Nevertheless, in this case the differential system can be  integrated
explicitly, and we obtain:
\begin{align}
\label{sl3pol}
h(x,y,u)=&c_0 + c_1 x+c_2y+c_3 u+c_4xy\\
&+c_5x(u-xy)+c_6uy+c_7u(u-xy),\nonumber
\end{align}
where $c_0,\dots,c_7 \in\mathbb{R}$. At this point, one can easily
calculate a basis of
$\tau(\Lie{g})$ given by homogeneous vector fields, and check the
following correspondence (up to
constants) between polynomials and Lie algebra generators, in the
sense that the polynomial $p$
corresponds to the unique multicontact vector field whose $U$
component is $pU$:
\begin{align*}
U&\leftrightarrow 1 & \theta U&\leftrightarrow u(u-xy)\\
Y&\leftrightarrow x & \theta Y&\leftrightarrow yu\\
X&\leftrightarrow y & \theta X&\leftrightarrow x(u-xy)\\
H_\alpha&\leftrightarrow u-2xy & &\\
H_\beta&\leftrightarrow u+xy,  & &
\end{align*}
where $H_\alpha$ and $H_\beta$ are the elements in $\Lie{a}$ that
represent the simple roots by the
Killing form.
\par
We remark that the polynomials that appear in the second column are
products of those in the
first column. Indeed, there exists $H\in\Lie{a}$ such that the
polynomial corresponding to $\tau(H)$
is $u-xy$ (namely $H=(2/3) H_\alpha + (1/2) H_\beta$), and similarly
there exists $H\in\Lie{a}$ such that the
polynomial corresponding to $\tau(H)$
is $u$ (namely $H=(1/3) H_\alpha + (2/3) H_\beta$). We can
interpret this property in a general
setting, by showing that the polynomials that correspond to the root
spaces labeled by the roots in a certain subset generate all the  others. This
remark  motivates the next chapter, where we give explicit  factorization formulas
for these polynomials in the  case that $\Lie{g}$ is a split
simple Lie algebra.

\chapter{Polynomial basis for split simple Lie algebras}~\label{po}
\vskip0.5cm
We investigate the
polynomial nature of the multicontact vector fields on some open subset  of an
Iwasawa nilpotent Lie group N. We consider the vector space $\mc{P}$
consisting of the polynomials that characterize the multicontact vector
fields on N, endowed with the Lie algebra structure induced by the  vector
space isomorphism of $\mc{P}$ with $\Lie{g}$ \cite{CDKR2}. In
particular, in the second section we compute some explicit formulas for  a
basis of
$\mc{P}$, pointing out some factorization properties.
\section{The polynomial algebra $\mc{P}$}\label{pol}
Let us summarize parts of the discussion of the previous chapter.
Theorem~\ref{cdkr} asserts that, under the assumption that the real  rank of
$\Lie{g}$ is greater than one, the Lie algebra of multicontact vector  fields
is $\tau(\Lie{g})$. This latter algebra may be viewed as a Lie algebra
of polynomials, because in the basis~\eqref{canonica}, the components
along the highest restricted root $\omega$ are polynomials and  determine all the
other components. Also, there is a natural notion of homogeneity that  comes into the
picture, which is very useful in order to describe the polynomials.
\par
Now, let us
consider again a basis adapted to the restricted root space
decomposition~\eqref{rsd}. Clearly the Lie algebra of multicontact  vector fields is
generated by the set
$$
\{\tau(X_{\alpha,i}),\alpha\in\Sigma\cup\{0\},i=1,\cdots,m_\alpha\}.
$$
As we just remarked, to any such vector field one associates a vector of
polynomials,  namely the coefficients along the
$\omega$-components. w We shall consider this problem
in the next section.
\section{The split case}\label{split}
Let $\Lie{g}$ be a real split simple Lie algebra, that is, a split real  form of its
complexification. This means that if $\Lie{g}^c=\Lie{g}\oplus i\Lie{g}$  is the
complexification of $\Lie{g}$, then there exists a Cartan subalgebra  $\Lie{h}$  of
$\Lie{g}^c$ such that if
$\Phi$ is the set of roots relative to the pair $(\Lie{g}^c,\Lie{h})$,  then
$\Lie{g}$ contains the real subspace of $\Lie{h}$ on which all the  roots in $\Phi$
are real, namely
$$
\left\{H\in\Lie{h}:\alpha(H)\in\R,\;\alpha\in\Phi\right\}\subset\Lie{g}.
$$
It is well-known that any complex semisimple Cartan subalgebra contains  a split real
form (see e.g. Corollary 6.10 in \cite{KN02}). The most relevant  consequences of
this assumption for our considerations are that $\Lie{m}=\{0\}$ and  that each
restricted root space has real dimension one. In particular, this  implies that
$I(\Lie{g}_\alpha)$ consists of the real multiples of a single  polynomial.
\par
Our decomposition formulas  are relative to a
suitable decomposition of the restricted root system that first
appeared in~\cite{C}.
Given two roots $\alpha$ and
$\beta$, the
$\alpha-$series of
$\beta$ is the set $\{\gamma \in \Sigma \cup \{0\} : \gamma= \beta +
n\alpha \}$. It turns out that the $\alpha-$series is an uninterrupted
string, namely that $n$ takes all the integer values in the interval
$[p,q]$, where the two integers $p\leq0$ and $q\geq0$ satisfy the
equality $p+q = -2 (\beta,
\alpha)/(\alpha, \alpha)$. In particular, for the $\omega-$series of  any root
$\beta \in \Sigma_+$, we have $p\in\{ -1,0 \}$. This means that
either $\omega - \beta \in \Sigma_+$ or $\omega - \beta$
is not a root, according as
$(\omega,\beta)=\frac{1}{2}(\omega,\omega)$ or
$(\omega,\beta) =0$. This gives us the natural decomposition of
$\Sigma_+$ into the disjoint union
$$
\Sigma_+ = \Sigma_{1/2} \cup \Sigma_0 \cup \Sigma_1,
$$
where
\begin{align*}
\Sigma_0&= \{ \beta\in \Sigma_+ : (\omega, \beta)=0 \},\\
\Sigma_{1/2}&= \{ \beta \in \Sigma_+ :
(\omega,\beta)=\frac{1}{2}(\omega,\omega)\},\\
\Sigma_1&=\{\omega\}.
\end{align*}
We shall write $\Delta_{1/2} = \Sigma_{1/2} \cap \Delta$ and
$\Delta_0 = \Sigma_0 \cap \Delta$.
According to the
decomposition of $\Sigma_+$, we put
$$
\Lie{n} =\Lie{n}_{(0)}\oplus\Lie{n}_{(1/2)}\oplus \Lie{n}_{(1)},
$$
with obvious notations.
Since $[\Lie{g}_\alpha , \Lie{g}_\beta] \subseteq
\Lie{g}_{\alpha + \beta}$, and $(\alpha + \beta,\omega) =
(\alpha,\omega) + (\beta, \omega)$, one has that $\Lie{n}_{(0)}$ is
a subalgebra and $ \Lie{n}_{(1/2)} \oplus \Lie{n}_{(1)}$
is an ideal in $\Lie{n}$.
Finally, we recall that the Cartan involution $\theta$ maps each root  space
$\Lie{g}_\alpha$ to $\Lie{g}_{-\alpha}$, so that
$\overline{\Lie{n}}=\theta \Lie{n}=\oplus_{\gamma\in\Sigma_-}  \Lie{g}_\gamma$,
where $\Sigma_-=-\Sigma_+$. According to the notations introduced  above, we
write
$$
\overline{\Lie{n}} =\overline{\Lie{n}}_{(0)}\oplus  \overline{\Lie{n}}_{(1/2)}
\oplus \overline{\Lie{n}}_{(1)},
$$
so that
\begin{equation}
\Lie{g}= \Lie{n}_{(1)} \oplus \Lie{n}_{(1/2)} \oplus
\Lie{n}_{(0)} \oplus
\Lie{a} \oplus \overline{\Lie{n}}_{(0)}\oplus
 \overline{\Lie{n}}_{(1/2)}\oplus \overline{\Lie{n}}_{(1)}.
\label{lunga}\end{equation}
By the linearity of scalar product, it is easy to check the
following commutation rules
\begin{align}
[\Lie{n}_{(r)},\Lie{a}]&\subset\Lie{n}_{(r)},
\hskip2truecm r=0,\frac{1}{2},1\nonumber\\
[\Lie{n}_{(0)},\overline{\Lie{n}}_{(0)}]
&\subset\Lie{n}_{(0)}\oplus\Lie{a} \oplus  \overline{\Lie{n}}_{(0)}\nonumber\\
[\Lie{n}_{(1/2)},\overline{\Lie{n}}_{(0)}]
&\subset\Lie{n}_{(1/2)}\nonumber\\
[\Lie{n}_{(1)},\overline{\Lie{n}}_{(0)}]&=\{0\}\nonumber\\
[\Lie{n}_{(0)},\overline{\Lie{n}}_{(1/2)}]
&\subset\overline{\Lie{n}}_{(1/2)}\label{rules}\\
[\Lie{n}_{(1/2)},\overline{\Lie{n}}_{(1/2)}]
&\subset\Lie{n}_{(0)}\oplus\Lie{a} \oplus  \overline{\Lie{n}}_{(0)}\nonumber\\
[\Lie{n}_{(1)},\overline{\Lie{n}}_{(1/2)}]
&\subset\Lie{n}_{(1/2)}\nonumber\\
[\Lie{n}_{(0)},\overline{\Lie{n}}_{(1)}]
&=\{0\}\nonumber\\
[\Lie{n}_{(1/2)},\overline{\Lie{n}}_{(1)}]
&\subset\overline{\Lie{n}}_{(1/2)}\nonumber\\
[\Lie{n}_{(1)},\overline{\Lie{n}}_{(1)}]
&\subset\Lie{a}.\nonumber
\end{align}
The proof of the next lemma is based on the above rules, and
leads us to a key explicit formula.
Take the following canonical coordinates on $\group{N}$:
$$
n=n_1  n_\frac{1}{2}  n_0 = \exp{(zZ)}
\exp{(\sum_{\alpha \in \Sigma_{1/2}} y_\alpha Y_\alpha)}
\exp{(\sum_{\beta \in \Sigma_0 } x_\beta X_\beta )},
$$
where $\{X_\beta:\beta\in\Sigma_0\}$ and $\{Y_\alpha:
\alpha\in\Sigma_{1/2}\}$ are a basis of
$\Lie{n}_{(0)}$ and $\Lie{n}_{(1/2)}$ respectevely, and
$Z\in\Lie{n}_{(1)}$.
\begin{lemma}\label{rid}
Let $X\in\Lie{g}_{\alpha}$, $\alpha\in\Sigma\cup\{0\}$, and $n$ in  $\group{N}$.
By the Bruhat lemma, for $t$ small enough there exists  $b(t)\in\group{P}$ such that
$\exp(tX)nb(t)\in\group{N}$. Consider the decomposition of  $n^{-1}\exp(tX)nb(t)$
with respect to the chosen coordinates, namely
$$
n^{-1}\exp(tX)
nb(t)=n^X_1(t)n^X_{1/2}(t) n^X_0(t).
$$
Write $n=n_{1}n_{1/2}n_{0}$,
then
\begin{itemize}
\item[(i)]
there exists $A\in\Lie{n}_{(1)}$ and $B\in
\Lie{n}_{(1/ 2)}\oplus\Lie{n}_{(0)}\oplus\Lie{a}\oplus\overline{\Lie{n}}$ such that
$$
n_{1/2}^{-1}n_{1}^{-1} \exp(tX) n_{1}n_{1/2}
=\exp (tA)\exp(tB)\exp(o(t));
$$
\item[(ii)]
$$
\frac{d}{dt}\,\left(n_1^X (t)\right)\Big|_{t=0}
=A.
$$
\end{itemize}
\end{lemma}
\begin{proof}
Write
\begin{align*}
n^{-1}\exp(tX)n=n_0 ^{-1} n_{1/2} ^{-1}
n_1 ^{-1} \exp{(tX)} n_1 n_{1/2} n_0 .
\end{align*}
Observe first that since $n_1=\exp(zZ)$,
$$
n_1 ^{-1} \exp{(tX)} n_1 =\exp(e^{-{\rm ad}(zZ)}tX).
$$
Now, by~\eqref{rules}
$$
[Z,X]\in
\begin{cases}
   \Lie{n}_{(1)}&\text{ if }X\in\Lie{a}\\
\Lie{n}_{(1/2)}&\text{ if }X\in\overline{\Lie{n}}_{(1/2)}\\
\Lie{a}&\text{ if }X\in\overline{\Lie{n}}_{(1)},
\end{cases}
$$
and if $X$ belongs to some other summand in the  decomposition~\eqref{lunga},
then $[Z,X]=0$. Therefore
\begin{align}\label{conuno}
n_1 ^{-1} \exp{(tX)} n_1 &= \exp(tX +t(H_1 +A_{1/2}+A_1) +o(t))\\
&=\exp(tA_1 +o(t))\exp(tX +t(H_1 +A_{1/2}) +o(t)),\nonumber
\end{align}
where $H_1\in\Lie{a}$, $A_{1/2}\in\Lie{n}_{(1/2)}$ and
$A_1\in\Lie{n}_{(1)}$.
\par
Secondly, since $\Lie{n}_{(1)}$ commutes with $\Lie{n}$, we consider
$$
n_{1/2}^{-1} \exp(tX +t(H_1 +A_{1/2}))n_{1/2}
=\exp(e^{{{-\sum_\alpha y_\alpha{\rm ad}   
Y_\alpha}}}(tX+tH_1+tA_{1/2})).
$$
Recall that $n_{1/2}$ is obtained exponentiating some element
in $\Lie{n}_{(1/2)}$. Thus, in the above formula  $\alpha\in\Sigma_{1/2}$, so that if
the commutator
$[Y_\alpha , X]\neq 0$, then by~\eqref{rules} the following  possibilities arise:
$$
[Y_\alpha,X]\in
\begin{cases}
\overline{\Lie{n}}_{1/2}&\text{ if }X\in\overline{\Lie{n}}_{(1)}\\
\Lie{a}&\text{ if }X\in\overline{\Lie{n}}_{(1/2)}\\
\Lie{n}_{(1/2)}&\text{ if }X\in\overline{\Lie{n}}_{(0)} \oplus
\Lie{n}_{(0)} \oplus\Lie{a}\\
\Lie{n}_{(1)}&\text{ if }X\in\Lie{n}_{(1/2)}.
\end{cases}
$$
Moreover,
$$
[Y_\alpha, H_1]\in\Lie{a},\qquad
[Y_\alpha, A_{1/2}]\in\Lie{n}_{(1)}.
$$
Therefore
\begin{align}\label{comezzo}
n_{1/2}^{-1} \exp(&tX +t(H_1 +A_{1/2}))n_{1/2}\\
&\hskip-1.5truecm=\exp(tX + t(B_{1/2}^- +H_2
+B_{1/2}+B_1)+o(t)),\nonumber\\ &\hskip-1.5truecm=\exp(tB_1  +o(t))\exp(tX
+ t(B_{1/2}^- +H_2 +B_{1/2})+o(t)),\nonumber
\end{align}
for some $B_{1/2}^- \in\overline{\Lie{n}}_{(1/2)}$, $H_2\in\Lie{a}$,
$B_{1/2}\in\Lie{n}_{(1/2)}$ and
$B_1\in\Lie{n}_{(1)}$.
Also, observe that by the Baker-Campbell-Hausdorff formula
\begin{equation}\label{ch1}
\exp(tL+o(t))=\exp(tL)\exp(o(t))
\end{equation}
for any $L\in\Lie{g}$.
Thus, by~\eqref{conuno} and~\eqref{comezzo} we deduce that
$$
n_{1/2}^{-1}n_{1}^{-1} \exp(tX) n_{1}n_{1/2}
=\exp (tA)\exp(tB)\exp(o(t)),
$$
with
$$
A =
\begin{cases}
A_1 + B_1 \hskip0.2cm &\text{ if } X\notin\Lie{n}_{(1)}\\
A_1 + B_1 +X \hskip0.2cm &\text{ if }  X\in\Lie{n}_{(1)}
\end{cases}
$$
and
$$
B=
\begin{cases}
 B_{1/2}^- +H_2 +B_{1/2}\hskip0.2cm &\text{ if } X\in\Lie{n}_{(1)}\\
 B_{1/2}^- +H_2 +B_{1/2}+X \hskip0.2cm &\text{ if }  X\notin\Lie{n}_{(1)}.
\end{cases}
$$
This proves (i).
Next, consider
$$
n_0^{-1} \exp(tX + t(B_{1/2}^- +H_2 +B_{1/2})) n_0=\exp(e^{-{\rm
ad}(\sum_\beta x_\beta X_\beta)}(tX +
tB_{1/2}^- +tH_2 +tB_{1/2})).
$$
If $[X_\beta,X]\neq 0$, then by~\eqref{rules}
$$
[X_\beta,X]\in
\begin{cases}
\overline{\Lie{n}}_{(1/2)}&\text{ if }X\in\overline{\Lie{n}}_{(1/2)}\\
\Lie{n}_{(0)} \oplus \overline{\Lie{n}}_{(0)} \oplus \Lie{a}&\text{
if }X\in\Lie{n}_{(0)} \oplus
\overline{\Lie{n}}_{(0)} \oplus \Lie{a}\\
\Lie{n}_{(1/2)}&\text{ if }X\in\Lie{n}_{(1/2)}.
\end{cases}
$$
Furthermore,
$$
[X_\beta , B_{1/2}^- ]\in\overline{\Lie{n}}_{(1/2)},\hskip0.2cm  [X_\beta,
H_2]\in\Lie{n}_{(0)},\hskip0.2cm [X_\beta, B_{1/2}]\in\Lie{n}_{(1/2)}.
$$
Hence
\begin{align}\label{cozero}
n_0^{-1} \exp(tX &+ t(B_{1/2}^- +H_2 +B_{1/2})) n_0\\
&=\exp(tX +t(C_{1/2}^- +C_{0}^- +H_3 +C_0^+ +C_{1/2})+o(t)),\nonumber
\end{align}
for some $C_{1/2}^- \in\overline{\Lie{n}}_{1/2}$, $C_{0}^-\in
\overline{\Lie{n}}_{(0)}$, $H_3\in\Lie{a}$,
$C_0^+\in\Lie{n}_{(0)}$ and $C_{1/2}\in\Lie{n}_{(1/2)}$.
\par
Recall that $\Lie{n}_{(1)}$ commutes with all $\Lie{n}$.
Thus using~\eqref{conuno}, \eqref{comezzo},
\eqref{ch1} and~\eqref{cozero}, we obtain
\begin{eqnarray}
n^{-1}\exp(tX)n&=&\exp(tA_1 +tB_1+o(t))\times\nonumber\\
                 & &\times \exp(tX+t(C_{1/2}^- +C_{0}^- +H_3  +C_0^+
+C_{1/2})+o(t))\nonumber\\
  &=&\exp(tA_1 +tB_1+o(t))\exp(tX+tC_0^+tC_{1/2}+o(t))\times\nonumber\\
& &\times \exp(t(C_{1/2}^- +C_{0}^- +H_3)+o(t))\nonumber\\
&=& \exp(tA_1 +tB_1 +tk_1 (X) )\exp (tC_{1/2}+tk_{1/2}  (X))\times\nonumber\\
& &\times \exp(tC_0 +tk_0 (X)) \exp(tC_{1/2}^- +tC_{0}^- +tH_3 +t
k(X))\times\nonumber\\
 & &\times\exp(o(t))\nonumber\\
&=&\exp(tA)\exp(tC)\exp(tD)\exp(tE)\exp(o(t)),\label{euno}
\end{eqnarray}
where
\begin{align*}
k(X) &= \begin{cases}
X &\text{ if }X\in \Lie{a}\oplus\overline{\Lie{n}}\\
0 &\text{ otherwise },
\end{cases}\\
k_i (X)&= \begin{cases}
X &\text{ if }X\in \Lie{n}_{(i)},\,i=0,1/2,1\\
0 &\text{ otherwise },
\end{cases}
\end{align*}
and
$$
C= C_{1/2}+k_{1/2} (X),\qquad D=C_0 +k_0 (X), \qquad E=C_{1/2}^-  +C_{0}^-
+H_3 +k(X).
$$
On the other hand, by hypothesis
\begin{equation}\label{edue}
n^{-1}\exp(tX)n = n^X_1(t)n^X_{1/2}(t) n^X_0(t) b(t)^{-1}.
\end{equation}
Observe that since $n^{-1}\exp(tX)n$ is the identity for $t=0$, then
necessarly $n_r^X(0)=e$ for every $r=1,1/2,0$, and $b(0)=e$.
Therefore, comparing~\eqref{euno} and~\eqref{edue},
\begin{align*}
\frac{d}{dt} \left( \exp(tA_{})\exp(tC) \exp(tD)\right. &\left.  \exp(tE)\exp(o(t))
\right)
\Big|_{t=0}\\
 &=\frac{d}{dt} \left( n^X_1(t)n^X_{1/2}(t) n^X_0(t) b(t)^{-1} \right)
\Big|_{t=0},
\end{align*}
whence
\begin{align*}
A + C + D + E
=& \, \frac{d}{dt} \left( n^X_1(t) \right)\Big|_{t=0}n^X_{1/2}(0)  n^X_0(0)
b(0)^{-1}\\ &+ n^X_1(0)\frac{d}{dt} \left( n^X_{1/2}(t)  \right)\Big|_{t=0}n^X_0(0)
b(0)^{-1}\\ &+ n^X_1(0)n^X_{1/2}(0)\frac{d}{dt} \left( n^X_0(t)  \right)\Big|_{t=0}
b(0)^{-1}\\ &+ n^X_1(0)n^X_{1/2}(0) n^X_0(0)\frac{d}{dt} \left(  b(t)^{-1}
\right)\Big|_{t=0}.
\end{align*}
This implies
$$
A
= \frac{d}{dt} \left( n^X_1(t) \right)\Big|_{t=0},
$$
because $A$ and $\frac{d}{dt} \left( n^X_1(t) \right)\Big|_{t=0}$ are  the only two
terms in the above sum that lie along $Z$. Thus also (ii) is proved.
\end{proof}
Let $X\in\Lie{g}_\alpha$, $\alpha\in\Sigma\cup\{0\}$. The multicontact
vector field associated to $X$ is defined by
$$
\tau(X)f(n)=\frac{d}{dt} f([\exp(tX) n])\Big|_{t=0},
$$
where $[\exp(tX) n ]$ is the $\group{N}$- component of
$\exp(tX) n$ in the Bruhat decomposition. This is equivalent to
saying that for $t$ small enough there exists $b(t)\in \group{P}$ such  that
$[\exp(tX) n]=\exp(tX)  n b(t)\in\group{N}$.
Hence
\begin{align*}
\tau(X)f(n)&= \frac{d}{dt} f(\exp(tX) nb(t))\Big|_{t=0}\\
&= \frac{d}{dt} f(nn^{-1}\exp(tX) nb(t))\Big|_{t=0}\\
&= \frac{d}{dt} f(n_1 n_{1/2} n_0 n^X_1(t)n^X_{1/2}(t)
n^X_0(t))\Big|_{t=0}.
\end{align*}
In the last part of the proof of Lemma~\ref{rid}, we observed that  $n^X_r (0)=e$,
for every $r=0,1/2,1$.
Recalling that $p$ is the
coefficient of $Z$ in the decomposition of
$\tau(X)$, by the latter assertion we have
\begin{equation}\label{formula}
p(n)=\frac{d}{dt}(n_1^X (t))\Big|_{t=0}.
\end{equation}
In particular, by (ii) of Lemma~\ref{rid} we desume that the
calculation of
$p(n)$ consists in computing the conjugation $ n_{1/2}^{-1} n_1 ^{-1}   
\exp(tX)
n_1 n_{1/2}$ and writing it in the form $\exp(tA)\exp(tB+o(t))$, with
$A\in\Lie{n}_{(1)}$ and $B\in \Lie{g}\setminus\Lie{n}_{(1)}$. In short,  $p(n)=A$.
\par
We shall obtain explicit formulas for the homogeneous polynomials
corresponding to $\Lie{g}$ using~\eqref{formula}. We
consider separetely the cases
with
$\alpha$ that lies respectevely in
$\Sigma_{0}$, $\Sigma_{1/2}$, $\Sigma_1$, $\{0\}$, $-\Sigma_{1}$,
$-\Sigma_{1/2}$, $-\Sigma_0$.
The formulas will point out that the polynomials corresponding to
$\Sigma_0$, $-\Sigma_{1/2}$, $-\Sigma_0$, $-\Sigma_1$
arise as products and suitable linear combinations of those
corresponding to the roots in $\Sigma_{1/2}$ and $\{0\}$.
Before collecting the formulas in a list of propositions, we still
introduce a couple of notations.
\par
We define on the set $\Sigma_{1/2}$ the equivalence
relation
$\sim$ given by
$$
\alpha \sim \beta \Leftrightarrow \alpha +
\beta =\omega,
$$
and we choose one representative for each element of the quotient
$(\Sigma_{1/2}/\sim~)$.
Denote the set of such representatives by $\tilde{\Sigma}_{1/2}$.
  From now until the end of this chapter, we write $p^\alpha$ and
$p^H$ for the polynomial that
corresponds to
$X_\alpha\in\Lie{g}_\alpha$ and to $H\in\Lie{a}$, respectevely.
\par
\begin{prop}\label{me}
{\rm (i)} If $\gamma \in \Sigma_{1/2}$, then
$p^\gamma (n) = c_{\gamma , \omega - \gamma} y_{\omega - \gamma}.$
\par
{\rm (ii)} If $H \in \Lie{a}$, then
$$
p^H (n) = \omega(H) z - \frac{1}{2} \sum_{[\alpha] \in
\Sigma_{1/2} / \sim } y_\alpha y_{\omega - \alpha} \left( (\omega -
\alpha)(H) - \alpha(H) \right) c_{\alpha,\omega-\alpha}.
$$
\par
{\rm (iii)} $p^\omega (n) =1.$
\end{prop}
\begin{proof}

(i)
If $\gamma \in \Sigma_{1/2}$ and $\alpha$ is another root in
$\Sigma_{1/2}$, then $\gamma + \alpha =
\omega$ provided $\gamma +\alpha$ is a root. This implies that
$\alpha=\omega-\gamma$. Furthermore, $[Z,\Lie{n}]=0$.
Therefore
\begin{align*}
n_{1/2}^{-1} n_1^{-1}\exp{(tY_\gamma)}
n_1 n_{1/2}
&=n_{1/2}^{-1}\exp\big( \sum_{n=0}^{+ \infty}(-1)^n  \frac{(\ad{zZ})^n}{n!}
tY_\gamma\big) n_{1/2} \\
&= n_{1/2}^{-1}
\exp(tY_\gamma)n_{1/2}\\
&= \exp \big( \sum_{n=0}^{+ \infty}(-1)^n
\frac{(\ad(\sum_{\alpha\in\Sigma_{1/2}}y_\alpha
Y_\alpha))^n}{n!} tY_\gamma \big)\\
   &= \exp ( tY_\gamma -
ty_{\omega - \gamma} [Y_{\omega - \gamma}, Y_\gamma])\\
&=\exp(tc_{\gamma,\omega-\gamma}y_{\omega-\gamma}
Z)\exp\big( tY_\gamma \big).
\end{align*}
By~\eqref{formula} and the remark thereafter, we have
$$
p^\gamma (n)= c_{\gamma,\omega-\gamma} y_{\omega - \gamma}.
$$
\par
(ii)
Since $[\Lie{n}_{(1/2)},\Lie{n}_{(1/2)}] \subseteq \Lie{n}_{(1)}$,
every bracket involving three or more
vectors in $\Lie{n}_{(1/2)}$ is zero. Then, for
$H
\in
\Lie{a}$, we have
\begin{align*}
n_{1/2}^{-1} n_1 ^{-1}& \exp{(tH)} n_1 n_{1/2}
= n_{1/2} ^{-1} \exp{(tH -tz[Z,H])}n_{1/2}\\
=&\exp(t\omega(H)zZ)\times\\
& \times \exp\big( tH -t\sum_{\alpha \in
\Sigma_{1/2}} y_\alpha
[Y_\alpha, H] + t/2 \sum_{\alpha + \beta = \omega} y_\alpha
y_\beta [Y_\beta,[Y_\alpha,H]] \big) \\
=& \exp{\big( \omega(H)z +\frac{1}{2}\sum_{\alpha + \beta=\omega}
\alpha(H) c_{\alpha,\beta} y_\alpha y_\beta \big) tZ}\ldots\\
\end{align*}
where again the only relevant component is the linear term in $t$ on  $Z$.
Therefore
$$
p^H (n) = \omega(H) z - \frac{1}{2} \sum_{[\alpha] \in
\Sigma_{1/2} / \sim } y_\alpha y_{\omega - \alpha} \left((\omega -
\alpha)(H) - \alpha(H) \right) c_{\alpha,\omega-\alpha},
$$
as required.
\par
(iii) Since $[Z,\Lie{n}]=0$, the conclusion is obvious.
\end{proof}
\begin{prop}\label{ze}
Let $\nu$ be either in $\Sigma_0$ or in $-\Sigma_0$. Let
$\mc{A}=\{\alpha\in\Sigma_{1/2} :\nu+\alpha\in\Sigma\}$, and write  $\mc{A}=\mc{A}_=
\cup \mc{A}_{\neq}$, where
$
\mc{A}_{\neq}  =\{\alpha\in \mc{A} :\alpha\neq\omega-(\nu+\alpha)\}
$
and
$
\mc{A}_= =\{\alpha \in \mc{A} : \alpha \notin
\mc{A}_{\neq} \}.
$
Then
$$
p^\nu (n) =\sum_{\alpha\in
\mc{A}_{\neq}}\frac{c_{\alpha,\nu}}{c_{\alpha,\omega-\alpha}} p^{\nu
+\alpha}(n) p^{\omega-\alpha}(n) + \frac{1}{2} \sum_{\alpha \in
\mc{A}_=}\frac{c_{\alpha,\nu}}{c_{\alpha,\omega-\alpha}}  [p^{\omega-\alpha}(n)]^2 .
$$
\end{prop}
\begin{proof}
First recall that if $\alpha \in \mc{A}$, then
$$
(\nu + \alpha,\omega
)=(\nu,\omega)+(\alpha,\omega)=\frac{1}{2}(\omega,\omega),
$$
whence $ \nu + \alpha \in \Sigma_{1/2}$.
Moreover, by definition $\omega + \nu$ cannot be a root even if $\nu$
is negative. Therefore,
proceeding as in the previous proposition, we obtain
\begin{align*}
& n_{1/2} ^{-1} n_1 ^{-1} \exp{(tX_\nu)} n_1 n_{1/2}
= n_{1/2}^{-1} \exp{(tX_\nu)}n_{1/2}  \\
&=  \exp{(tX_\nu - t\sum_{\alpha\in\Sigma_{1/2}} y_\alpha
[Y_\alpha,X_\nu]+\frac{t}{2}
\sum_{\alpha_1,\alpha_2 \in \Sigma_{1/2}} y_{\alpha_1} y_{\alpha_2}
[Y_{\alpha_2},
[Y_{\alpha_1},X_\nu]])} \\
&=\exp{\big( \frac{t}{2}\sum_{\nu+\alpha_1+\alpha_2=\omega}
c_{\alpha_1,\nu}c_{\alpha_2,\nu + \alpha_1} y_{\alpha_1} y_{\alpha_2}
Z\big) }\exp\big( tX_\nu-t\sum_{\alpha\in
A}c_{\alpha,\nu}y_\alpha Y_{\alpha+\nu}\big)
\end{align*}
Again by~\eqref{formula},
$$
p^\nu (n)= \frac{1}{2} \sum_{\nu+\alpha_1+\alpha_2=\omega}
c_{\alpha_1,\nu}c_{\alpha_2,\nu + \alpha_1}
y_{\alpha_1} y_{\alpha_2},
$$
where $\alpha_1$ and $\alpha_2$ are in $\Sigma_{1/2}$.
\par
We now notice that
   if $\alpha_1 \neq \alpha_2$, then both
$\nu+\alpha_1+\alpha_2$ and $\nu+\alpha_2+\alpha_1$
are $\omega-$chains. Thus the coefficient of the monomial
$y_{\alpha_1} y_{\alpha_2}$ is
$c_{\alpha_1,\nu} c_{\alpha_2,\nu+\alpha_1} +
c_{\alpha_2,\nu}  c_{\alpha_1,\nu+\alpha_2}$.
Furthermore, since $\nu+\alpha_1+\alpha_2 = \omega$,
the root $\alpha_2$ must be equal to $\omega-
(\nu+\alpha_1)$. Using the Jacobi identity:
$$
[X_\alpha,[X_{\omega-(\nu+\alpha)},X_\nu]]=[X_{\omega- (\nu+\alpha)},[X_\alpha,X_\nu]]
+ 0
$$
{\em i.e.}
$$
c_{\omega-(\nu+\alpha),\nu}c_{\alpha,\omega-\alpha}
=c_{\alpha,\nu}c_{\omega-(\nu+\alpha),\nu+\alpha}.
$$
So, by (i) of Proposition~\ref{me}, we can write $p^\nu$ as
follows
$$
p^\nu (n) =\sum_{\alpha\in
\mc{A}_\neq}\frac{c_{\alpha,\nu}}{c_{\alpha,\omega-\alpha}} p^{\nu
+\alpha}(n) p^{\omega-\alpha}(n) + \frac{1}{2} \sum_{\alpha \in
\mc{A}_=}\frac{c_{\alpha,\nu}}{c_{\alpha,\omega-\alpha}}  [p^{\omega-\alpha}(n)]^2.
$$
\end{proof}
\bigskip
Since $[X_\alpha , X_{-\alpha} ] = B(X_\alpha
,X_{-\alpha})H_\alpha$, by suitably normalizing the basis vectors,
one may assume that
$B(X_\alpha, X_{-\alpha}) = 1$. In order to semplify the notations,
from now on we fix such a basis.
Then the following relations for the structure constants hold (see
\cite{KN02}, Sec.1, Chap.VI):
\begin{eqnarray} \label{struttura}
c_{\alpha,\beta}&=&c_{\beta,\gamma}=c_{\gamma,\alpha} ,\nonumber \\
c_{-\alpha,\alpha + \beta} &=& c_{\alpha, \beta},
\end{eqnarray}
for every
$\alpha,\beta,\gamma \in \Sigma$ s.t. $\alpha + \beta +
\gamma =0$.
\begin{prop}\label{ne}
   If $\gamma \in \Sigma_{1/2}$, then
\begin{equation}
p^{-\gamma} (n) = -\frac{1}{c_{\gamma,\omega-\gamma}}
p^{\omega-\gamma} (n) p^{H(\gamma)} (n) +
\frac{1}{3}
\sum_{\alpha\in\Sigma_{1/2}:-\gamma + \alpha \in \pm \Sigma_0}
\frac{c_{\alpha,-\gamma}}{c_{\alpha,\omega-\alpha}} p^{\omega-\alpha}
(n) p^{\alpha - \gamma} (n),\nonumber
\end{equation}
where $H=H(\gamma) \in \Lie{a}$ is the solution of the linear system
\begin{equation} \label{sistema1}
\begin{cases}
\omega(H)=-\omega(H_\gamma)\\
(3\alpha-\omega)(H)=-\alpha(H_\gamma) \hspace{.3in} \forall \alpha
\in \tilde{\Sigma}_{1/2}.
\end{cases}
\end{equation}
Furthermore, let $\gamma$ be a simple root and $H(\gamma)$ the
corresponding solution of
{\rm \eqref{sistema}}. If
$\gamma^\prime =
\gamma +
\delta_1 + \delta_2 + \ldots + \delta_p$, for some simple roots
$\delta_1,\ldots,\delta_p \in
\Sigma_0$, then
\begin{equation} \label{altrigamma}
H(\gamma^\prime) = H(\gamma) -\frac{1}{3} (H_{\delta_1} + \ldots +
H_{\delta_p}).
\end{equation}
\end{prop}
\begin{proof}
As before
\begin{align*}
&n_{1/2} ^{-1} n_1 ^{-1}  \exp{(tY_{-\gamma})}  n_1 n_{1/2}
= n_{1/2} ^{-1}
\exp{(tY_{-\gamma} -t c_{\omega,-\gamma} zY_{\omega-\gamma} )} n_{1/2}  \\
&\hskip0.5cm = \exp\big( tY_{-\gamma} -  tc_{\omega,-\gamma}zY_{\omega-\gamma}
-t\sum_{\alpha\in\Sigma_{1/2}} y_\alpha[Y_\alpha,Y_{-\gamma}]\\
&\hskip2cm +tz\sum_{\alpha\in\Sigma_{1/2}}
c_{\omega,-\gamma} y_\alpha [Y_\alpha,Y_{\omega-\gamma}]\\
&\hskip2cm  +\frac{t}{2}\sum_{\alpha_1,\alpha_2 \in  \Sigma_{1/2}}y_{\alpha_1}
y_{\alpha_2}
[Y_{\alpha_2},[Y_{\alpha_1}, Y_{-\gamma}]\\
&\hskip2cm  - \frac{t}{6} \sum_{\alpha_1,\alpha_2,\alpha_3 \in
\Sigma_{1/2}} y_{\alpha_1} y_{\alpha_2} y_{\alpha_3}
[Y_{\alpha_3},[Y_{\alpha_2},[Y_{\alpha_1},Y_{-\gamma}]]]\big)  .
\end{align*}
Since $(-\gamma + \alpha,\omega)=-(\gamma,\omega) +(\alpha, \omega)=
-\frac{1}{2}(\omega,\omega)+\frac{1}{2}(\omega,\omega)=0$ for every
$\alpha \in \Sigma_{1/2}$, it
follows that $-\gamma+\alpha$ is either in $\pm \Sigma_0$ or $0$ or
not a root. This implies
that the bracket $[Y_{\alpha_1},Y_{-\gamma}]$ is respectively in
$\Lie{n}_{(0)}$, $\Lie{a}$ or zero.
Then, by~\eqref{formula} we have
\begin{align*}
p^{-\gamma} (n) =& -\omega(H_\gamma)y_\gamma z \hspace{.1in} +
\hspace{.1in} \frac{1}{6}
\sum_{\alpha\in\Sigma_{1/2}} \alpha(H_\gamma)
c_{\omega-\alpha,\alpha} y_\gamma y_\alpha
y_{\omega-\alpha}\\
   &-\frac{1}{6} \sum_{-\gamma + \alpha_1 + \alpha_2
+\alpha_3 =\omega}
c_{\alpha_1,-\gamma} c_{\alpha_2,-\gamma+\alpha_1} c_{\alpha_3,
-\gamma+ \alpha_1 + \alpha_2}
y_{\alpha_1}y_{\alpha_2}y_{\alpha_3},
\end{align*}
with $\alpha_1$,$\alpha_2$,$\alpha_3 \in \Sigma_{1/2}$, and where we  used that
$c_{\omega,-\gamma}c_{\gamma,\omega-\gamma}=-\omega(H_\gamma)$
(Jacobi identity). In particular, if
$-\gamma +
\alpha_1
\in
\pm
\Sigma_0$, by Proposition~\ref{ze} we can write the factor
$y_{\alpha_2} y_{\alpha_3} (c_{\alpha_2,
-\gamma +\alpha_1}c_{\alpha_3, -\gamma+\alpha_1+\alpha_2})$ as
$2p^{-\gamma+\alpha_1} (n)$. Moreover,
by (i) of Proposition~\ref{me},
$y_{\alpha_1}=
-\frac{p^{\omega-\alpha_1}(n)}{c_{\alpha_1,\omega-\alpha_1}}$ and
$y_\gamma=-\frac{p^{\omega-\gamma}(n)}{c_{\gamma,\omega-\gamma}}$.  Therefore
\begin{align*}
p^{-\gamma} (n)  &=-\frac{p^{\omega-\gamma}(n)}{c_{\gamma,\omega-\gamma}} \big\{
   -\omega(H_\gamma)z \phantom{\sum_{\Sigma_{1/2}}}    \\
   &\hskip0.5cm +  \frac{1}{6} \sum_{\alpha\in\Sigma_{1/2}/\sim}
c_{\alpha,\omega-\alpha} y_\alpha  y_{\omega-\alpha}((\omega-\alpha)(H_\gamma)
   -\alpha(H_\gamma)) \big\}\\
&\hskip0.5cm
+\frac{1}{3}\sum_{\alpha\in\Sigma_{1/2}/-\gamma+\alpha\in\pm\Sigma_0}
\frac{c_{\alpha,-\gamma}}{c_{\alpha,\omega-\alpha}}
p^{\omega-\alpha}(n)  p^{-\gamma+\alpha} (n).
\end{align*}
Consider the polynomial in curly braces and compare it with (ii) of
Proposition~\ref{me}. We desume
that it has the form of a polynomial corresponding to some element $H
\in \Lie{a}$, provided that $H$
satisfies
$$
\begin{cases}
\omega(H) =-\omega(H_\gamma)\\
-\frac{1}{2}
c_{\alpha,\omega-\alpha}\left(
(\omega-\alpha)(H)-\alpha(H)\right)=\frac{1}{6}c_{\alpha,\omega-\alpha}
\left( (\omega-\alpha)(H_\gamma)-\alpha(H_\gamma)\right) ,
\end{cases}
$$
for every $ \alpha \in\tilde{\Sigma}_{1/2}$.
This system is equivalent to
\begin{equation*}
\begin{cases}
\omega(H)=-\omega(H_\gamma)\\
(3\alpha-\omega)(H)=-\alpha(H_\gamma) \hspace{.3in} \forall \alpha
\in \tilde{\Sigma}_{1/2}.
\end{cases}
\end{equation*}
In order to conclude the proof of the first statement, it is enough
to prove that this linear system
has a solution.
To this end,
take a maximal set of linear independent vectors in
$\tilde{\Sigma}_{1/2}$, and denote it by $\mc{B}
(\tilde{\Sigma}_{1/2})$. Since the restricted roots generate a vector  space of
dimension equal to the rank of $\Lie{g}$, say $l$, there are at most
$l$ elements in $\mc{B}
(\tilde{\Sigma}_{1/2})$.  By inspection of the complete Dinkin
diagrams of split semisimple
Lie algebras (\cite{B}) one sees that the cardinality of $\mc{B}
(\tilde{\Sigma}_{1/2})$ is at
least $l-1$. Hence there are two possible cases.
\par
(a) $\# \mc{B} (\tilde{\Sigma}_{1/2}) =l-1$, say $\mc{B}
(\tilde{\Sigma}_{1/2}) =\{
\alpha_1,\ldots,\alpha_{l-1}\}$. In this case $\beta=
\sum_{i=1}^{l-1} a_i \alpha_i$ for every $\beta
\in
\tilde{\Sigma}_{1/2}$. Since $\beta \in \Sigma_{1/2}$, its inner
product with $\omega$ is half of the
square norm of $\omega$, hence  $\frac{1}{2}(\omega,\omega)=(\sum_{i=1}^{l-1}
a_i\alpha_i,\omega)=\sum_{i=1}^{l-1} a_i (\alpha_i ,
\omega)= \sum_{i=1}^{l-1} a_i \frac{1}{2} (\omega,\omega)$, which  implies $
\sum_{i=1}^{l-1} a_i=1$.
If $H$ is a solution of the subsystem of (\ref{sistema1})
\begin{equation} \label{sottosistema}
\begin{cases}
\omega(H)=-\omega(H_\gamma)\\
(3\alpha_i-\omega)(H)=-\alpha_i(H_\gamma) \hspace{.3in} \forall  i=1\ldots l-1,
\end{cases}
\end{equation}
then $H$ solves also (\ref{sistema}).
Indeed
$$
(3\beta-\omega)(H)=\left(3\sum_{i=1}^{l-1} a_i \alpha_i -
\sum_{i=1}^{l-1} a_i \omega \right)(H) =
-\sum_{i=1}^{l-1} a_i \alpha_i (H_\gamma)=-\beta(H_\gamma).
$$
The linear system (\ref{sottosistema}) has $l$ equations and $l$
variables, and the associated
matrix is diagonal. Therefore it has exactly one solution.
\par
(b) $\# \mc{B} (\tilde{\Sigma}_{1/2}) =l$, so that
$\omega=\sum_{i=1}^{l} b_i \alpha_i$,
that is $(\omega,\omega)=\sum_{i=1}^{l} b_i
(\alpha_i,\omega)=\sum_{i=1}^{l} b_i \frac{1}{2}
(\omega,\omega)$. This implies that $\sum_{i=1}^{l} b_i=2$. It turns
out that any equation in
(\ref{sistema1}) depends linearly on the set of equations
$$
(3\alpha_i-\omega)(H)=-\alpha_i(H_\gamma) \hspace{.3in} \forall  i=1\ldots l.
$$
Indeed
$$
-\omega(H_\gamma)=-\sum_{i=1}^{l} b_i \alpha_i (H_\gamma) = \big(
3\sum_{i=1}^{l} b_i \alpha_i(H) -
\sum_{i=1}^{l} b_i \omega \big) (H) = \omega(H),
$$
and
$$
(3\beta-\omega)(H)=\big( 3\sum_{i=1}^{l} a_i \alpha_i -
\sum_{i=1}^{l} a_i \omega \big) (H) =
-\sum_{i=1}^{l} a_i \alpha_i (H_\gamma)=-\beta(H_\gamma),
$$
where $\beta=\sum_{i=1}^{l}a_i \alpha_i$ is any root in  $\tilde{\Sigma}_{1/2}$.
\par
Let now  $\gamma \in \Delta_{1/2}$, and let $H(\gamma)$ be the
corresponding solution of
(\ref{sistema1}). Let $\gamma^\prime$ be another root in
$\Sigma_{1/2}$ and $\delta_1,\ldots ,
\delta_p$ simple roots such that
$\gamma^\prime = \gamma + \delta_1 + \ldots + \delta_p$. It is an
easy calculation to check that
$H(\gamma^\prime) = H(\gamma) -\frac{1}{3} (H_{\delta_1} + \ldots +
H_{\delta_p})$.
\end{proof}
{\bf Remarks.}
 From the theory of root systems it follows  that every root
$\gamma^\prime$ in
$\Sigma_{1/2}$ can be written as $\gamma + \delta_1
+\ldots+\delta_p$, with $\gamma \in \Delta_{1/2}$
and for some $\delta_1 ,\dots,\delta_p \in
\Delta_0$ (see, e.g., Lemma 3.5 in \cite{CDKR2}).
This fact, toghether with~\eqref{altrigamma}, tells us that we must
solve~\eqref{sistema1} only for $\gamma\in\Delta_{1/2}$.
  Furthermore, by the
classification of root systems,
we know that there exists exactly one simple root belonging to
$\Sigma_{1/2}$, except the case of the root system
$A_n$, for which $\Sigma_{1/2}$ consists of two roots \cite{B}.
\begin{prop}
Fix $H =
\frac{1}{\sqrt{2\omega(H_\omega)}} H_\omega$. Then
$$
p^{-\omega} (n) = -(p^H (n))^2 -\frac{1}{4}
\sum_{\alpha\in\Sigma_{1/2}} p^{\omega-\alpha}(n)
p^{\alpha-\omega}(n).
$$
\end{prop}
\begin{proof}
Notice that in order to complete $-\omega$ to an $\omega-$chain we
need exactly four roots in
$\Sigma_{1/2}$. We obtain
\begin{align*}
 & n_{1/2}^{-1} n_1^{-1} \exp tX_{-\omega} n_1 n_{1/2}\\
 &\hskip0.5cm = n_{1/2}^{-1} \exp
        \big( tX_{-\omega} -tzH_\omega - \frac{t}{2} z^2
        \omega(H_\omega)Z\big) n_{1/2}\\
 &\hskip0.5cm = \exp(-\frac{t}{2} z^2 \omega(H_\omega)Z) \times\\
 &\hskip1cm \times \exp\big(tX_{-\omega} -tz  H_\omega-t\sum_{\alpha\in\Sigma_{1/2}}
          c_{\alpha,-\omega}y_\alpha Y_{-\omega+\alpha}\\
 &\hskip1cm -tz\sum_{\alpha\in\Sigma_{1/2}}\alpha(H_\omega)y_\alpha  Y_\alpha
      +\frac{t}{2}
       \sum_{\alpha_1,\alpha_2 \in\Sigma_{1/2}}
    y_{\alpha_1}y_{\alpha_2}[Y_{\alpha_2},[Y_{\alpha_1},X_{-\omega}]]\\
 &\hskip1cm +\frac{t}{2}z
\sum_{\alpha\in\Sigma_{1/2}}c_{\omega- \alpha,\alpha}\alpha(H_\omega)y_\alpha
        y_{\omega-\alpha}Z\\
 &\hskip1cm -\frac{t}{6}\sum_{\alpha_1,\alpha_2,\alpha_3  \in\Sigma_{1/2}}
       y_{\alpha_1}y_{\alpha_2}y_{\alpha_3}
       [Y_{\alpha_3},[Y_{\alpha_2},[Y_{\alpha_1},X_{-\omega}]]]\\
 &\hskip1cm  +\frac{t}{24}
       \sum_{\alpha_1,\alpha_2,\alpha_3,\alpha_4\in\Sigma_{1/2}}
y_{\alpha_1}y_{\alpha_2}y_{\alpha_3}y_{\alpha_4}[Y_{\alpha_4},[Y_{\alpha _3},
       [Y_{\alpha_2},[Y_{\alpha_1},X_{-\omega}]]]]\big) \\
 &\hskip0.5cm =\exp\big( \big( -\frac{t}{2} z^2 \omega(H_\omega)  +\frac{t}{2}z
\sum_{\alpha\in\Sigma_{1/2}}c_{\omega- \alpha,\alpha}\alpha(H_\omega)y_\alpha
        y_{\omega-\alpha}\\
 &\hskip1cm  +
\frac{t}{24} \sum_{\alpha_1,\alpha_2,\alpha_3,\alpha_4
\in\Sigma_{1/2}}c_{\alpha_1,-\omega}c_{\alpha_2,- \omega+\alpha_1}\times\\
 &\hskip4cm  \times\, c_{\alpha_3,-\omega+\alpha_1+\alpha_2}
             c_{\alpha_4,-\omega+\alpha_1
+\alpha_2+\alpha_3}y_{\alpha_1}y_{\alpha_2}y_{\alpha_3}y_{\alpha_4}\big)
             Z\big)\ldots.
\end{align*}
Writing the corresponding polynomial, we can split the last summand of
the  above formula
   according to the fact that the chains $-\omega + {\alpha}_1 +
{\alpha}_2 + {\alpha}_3 +{\alpha}_4$
are of two kinds. In fact we have either ${\alpha}_2=\omega
-{\alpha}_1$, so that $-\omega + {\alpha}_1
+ {\alpha}_2=0$, or $-\omega + {\alpha}_1 + {\alpha}_2 \in \pm
{\Sigma}_0$. Therefore
\begin{align}
p^{-\omega}(n) = &-\frac{1}{2} z^2 \omega(H_\omega) \nonumber\\
&+\frac{1}{2}z
\sum_{\alpha\in\Sigma_{1/2}}c_{\omega- \alpha,\alpha}\alpha(H_\omega)y_\alpha
y_{\omega-\alpha}\label{pezzo}\\
&-\frac{1}{24}\sum_{\alpha,\beta\in\Sigma_{1/2}}
c_{\alpha,-\omega}c_{\omega-\beta,\beta}
\beta(H_{\omega-\alpha})y_\alpha
y_{\omega-\alpha}y_\beta y_{\omega-\beta}\nonumber\\
&+\frac{1}{24} \sum_{\alpha_i
\in\Sigma_{1/2}:-\omega+\alpha_1+\alpha_2\in\pm\Sigma_0}
c_{\alpha_1,-\omega}c_{\alpha_2,-\omega+\alpha_1}\times\nonumber\\
&\hskip3cm \times c_{\alpha_3,-\omega+\alpha_1+\alpha_2}
c_{\alpha_4,-\omega+\alpha_1+\alpha_2+\alpha_3}y_{\alpha_1}
y_{\alpha_2}y_{\alpha_3}y_{\alpha_4}.\nonumber
\end{align}
Since for every $\alpha\in\Sigma_{1/2}$
$$
B(H_\omega,H)=\omega (H)=(\omega -\alpha) (H) + \alpha
(H)=B(H_{\omega -\alpha},H) +B(H_\alpha
,H),
$$
it follows that $H_\omega = H_{\omega-\alpha}+H_{\alpha}$. Moreover,
\begin{align*}
\sum_{\alpha\in\Sigma_{1/2}} c_{\omega-\alpha,\alpha}
\alpha(H_\omega)y_\alpha y_{\omega-\alpha}
&=\sum_{\alpha\in\Sigma_{1/2}} c_{\omega-\alpha,\alpha}
\frac{\omega(H_\omega)}{2}  y_\alpha
y_{\omega-\alpha}\\
&=\frac{\omega(H_\omega)}{2}\sum_{\alpha\in\Sigma_{1/2}}
c_{\omega-\alpha,\alpha}   y_\alpha
y_{\omega-\alpha}\\
&= 0,
\end{align*}
because
$$
\sum_{\alpha\in\Sigma_{1/2}} c_{\omega -\alpha ,\alpha} y_\alpha
y_{\omega -\alpha}
=\sum_{\alpha\in\tilde{\Sigma}_{1/2}} ( c_{\omega -\alpha ,\alpha} -
c_{\omega -\alpha
,\alpha})y_\alpha  y_{\omega -\alpha}.
$$
Comparing with the polynomial formula of Proposition~\ref{ne}, and
observing that~\eqref{struttura} implies
$$
c_{\omega-\alpha,\alpha}=c_{\alpha,-\omega},
$$
the sum in~\eqref{pezzo} becomes
\begin{align*}
&\frac{1}{2}z
\sum_{\alpha\in\Sigma_{1/2}}c_{\omega- \alpha,\alpha}\alpha(H_\omega)y_\alpha
y_{\omega-\alpha}\\
& - \frac{1}{24}\sum_{\alpha,\beta\in\Sigma_{1/2}}
c_{\alpha,-\omega}c_{\omega-\beta,\beta}\beta(H_{\omega-\alpha})y_\alpha
y_{\omega-\alpha}y_\beta y_{\omega-\beta}\\
& + \frac{1}{24} \sum_{\alpha_i
\in\Sigma_{1/2}:-\omega+\alpha_1+\alpha_2\in\pm\Sigma_0}c_{\alpha_1,- \omega}
   c_{\alpha_2,-\omega+\alpha_1}\times \\
& \hskip3cm \times \, c_{\alpha_3,-\omega+\alpha_1+\alpha_2}
c_{\alpha_4,-\omega+\alpha_1+\alpha_2+\alpha_3}
y_{\alpha_1}y_{\alpha_2}y_{\alpha_3}y_{\alpha_4}\\
&\\
=& -\frac{1}{4} \big\{
\sum_{\alpha\in\Sigma_{1/2}}c_{\alpha,-\omega}y_\alpha\big(
-\omega(H_{\omega-\alpha})zy_{\omega-\alpha}\\
&+\frac{1}{6}y_{\omega-\alpha}\sum_{\beta\in\Sigma_{1/2}}
c_{\omega-\beta,\beta}\beta(H_{\omega-\alpha})y_\beta
y_{\omega-\beta}\\
   &-\frac{1}{6}
\sum_{\alpha_i \in
\Sigma_{1/2}:-\omega+\alpha+\alpha_1\in\pm\Sigma_0}c_{\alpha_1,- \omega+\alpha}
c_{\alpha_2,-\omega+\alpha+\alpha_1}
c_{\alpha_3,-\omega+\alpha+\alpha_1+\alpha_2}y_{\alpha_1}
y_{\alpha_2} y_{\alpha_3}\big)\big\}\\
=&-\frac{1}{4} \sum_{\alpha\in\Sigma_{1/2}}
p^{\omega-\alpha}(n) p^{\alpha-\omega} (n),
\end{align*}
where the last equality follows by (i) of Proposition~\ref{me} and
Proposition~\ref{ne}.
Finally, by Proposition~\ref{ze}, $\frac{z^2}{2} \omega(H_\omega) =
p^H (n) p^H (n)$ if $H$ satisfies
$$
\begin{cases}
\omega(H)-2\alpha(H) =0\\
\omega(H)=\sqrt{\frac{\omega(H_\omega)}{2}},
\end{cases}
$$
{\em i.e.}
$$
\begin{cases}
\alpha(H) = \frac{1}{2} \omega(H)\\
\omega(H) = \sqrt{\frac{\omega(H_\omega)}{2}}.
\end{cases}
$$
It is a simple calculation to verify that
$H=\frac{1}{\sqrt{2\omega(H_\omega)}}H_\omega$ satisfies the
equations above.
We then conclude that
$$
p^{-\omega} (n) = -(p^H (n))^2 -\frac{1}{4}
\sum_{\alpha\in\Sigma_{1/2}} p^{\omega-\alpha}(n)
p^{\alpha-\omega}(n),
$$
as required.
\end{proof}

\chapter{Hessenberg manifolds}\label{hesse}
\vskip0.5cm
The Hessenberg manifolds arise as a natural class of submanifolds of  the spaces
G/P. A crucial point of the present work is to investigate in detail the
stratification of the tangent bundle of these manifolds. We shall  see
that they inherit from G/P a structure that allows us to define the
appropriate version of multicontact mapping. In the first section we
define the classical Hessenberg manifolds, viewed as submanifolds of the
complete flag manifolds. In the second section we define them
in the more general context  of G/P,
where G is a real semisimple Lie group and P is a minimal parabolic
subgroup of G. The results and definitions that are considered in the
second section are taken from \cite{DP}.

\section{The basic context: the Hessenberg flags}
We collect some notions about the Hessenberg flag manifolds. For
further details, see~\cite{FDMth}. The presentation that follows differs
somewhat from the standard version outlined in the introduction. The  main
reason for doing so is that our natural (local) environment is the
Iwasawa nilpotent group that in the standard setting would be
$\overline{\group{N}}$, a lower triangular group. We find it more
natural to be working on N, the unipotent upper triangular group.
\par
 Let
$\group{G}=\group{SL}(n,\mathbb{R})$, and let P
be its minimal parabolic subgroup given by the
lower triangular matrices. The homogeneous space
obtained by the quotient G/P realizes the
complete $flag$ $manifold$, in the following
sense. Define
$$
{\rm
Flag}(n)=\{(S_1,\dots,S_{n-1}):S_1\subset\dots\subset S_{n-1}\},
$$
where each $S_k$ is a subspace of $\mathbb{R}^n$ such that ${\rm  dim}S_k=k$. The set ${\rm
Flag}(n)$ is a smooth and compact manifold, called complete flag  manifold. It is easy to
check  that
$\group{SL}(n,\mathbb{R})$ acts in a transitive way on
${\rm Flag}(n)$ by the natural action
$$g: (S_1,\dots,S_{n-1})\mapsto (gS_1,\dots,gS_{n-1}),$$
for every $g\in\group{SL}(n,\mathbb{R})$.
We can view any flag as a matrix of column vectors
\begin{equation} \label{flag}
\bmatrix
v_n &\cdots&v_1
\endbmatrix
\end{equation}
where $S_1={\rm span}(v_1)$, $S_2={\rm  span}(v_1,v_2),\dots,\mathbb{R}^n={\rm
span}(v_1,v_2,\cdots,v_n)$. Consider the flag corresponding to the  identity matrix. The
isotropy group at this flag is the group P of lower triangular  matrices, so that
$${\rm Flag}(n)=\group{G}/\group{P}.$$
Each matrix $A\in \group{SL}(n,\mathbb{R})$ representing a flag can be  reduced to
a  unipotent upper triangular matrix by changing the representative in  G/P. This is done by
the Gauss algorithm, provided we restrict ourselves to an open and  dense subset of G/P where
we assume that some entries are not zero in order to make the Gauss
 algorithm work.
Summarizing, the nilpotent subgroup of $\group{SL}(n,\mathbb{R})$  defined by the unipotent
upper triangular matrices and denoted by N is identified with an open  and
dense subset of ${\rm Flag}(n)$.
\par
Consider now the Lie algebra $\Lie{sl}(n,\mathbb{R})$ given by the  matrices whose trace is
zero, and take a diagonal matrix $H$ in $\Lie{sl}(n,\mathbb{R})$ with  distinct entries and
non--zero determinant. Furthermore, fix an
integer
$p$ such that $1\leq p <n-1$. We say that a flag
$(S_1,\dots,S_{n-1})$ is a $type$ $p$
$Hessenberg$ $flag$ $for$ $H$ if the matrix $H$
shifts the linear space
$S_i$ within $S_{i+p}$, that is
$$HS_i\subset S_{i+p}.$$
The set of all such flags is a smooth manifold that we denote by
${\rm Hess}_p(H)$ and call {\it  p-th
Hessenberg manifold}. By considering
again the flag manifold as the
homogeneous space G/P, a flag
$X\in\group{G}/\group{P}$ is a type $p$ Hessenberg flag for $H$ if 
and only if it satisfies
\begin{equation} \label{hessflag}
HX=XR,
\end{equation}
where $R=(r_{i,j})$ is a matrix in $\Lie{sl}(n,\mathbb{R})$ such that
$$
r_{i,j}=
\begin{cases}
0 \hskip0.1cm\text{ if }(i,j)=(i,p+i+1)\text{ or }(i,j)=(p+j+1,j)\\
* \hskip0.1cm\text{ otherwise},
\end{cases}
$$
i.e. $R$ is a matrix where all entries above the {\it p+1-th} diagonal  are zero. In
order to see that (\ref{hessflag}) is equivalent to the definition of a  Hessenberg flag it is
enough to visualize $X$ as a matrix of column vectors as in  (\ref{flag}), and to notice that
the product $XR$ maps the last column, that corresponds to $S_1$, in a  linear
combination of the last {\it p+1-th}, and so on.
Since $X$ is invertible, we can rewrite (\ref{hessflag}) as
\begin{equation} \label{char}
X^{-1}HX=R.
\end{equation}
The formula above can be used for computing local algebraic equations  for ${\rm
Hess}_p(H)$. We restrict ourselves to the dense subset N of upper
unipotent  triangular matrices, and we write $X\in\group{N}$ as
$(x_{i,j})_{i,j}$, with $i<j$. Moreover, write
$H=(\lambda_i)_i$. Then we compute the product $X^{-1} H X$ and we put  equal to zero the
entries above the {\it p+1-th} diagonal, obtaining
\begin{align} \label{polhess}
f_{i,j}(X)&=(\lambda_j - \lambda_i)x_{i,j} \\
&\hskip-0.5truecm+ \sum_{t=1}^{j-i-1}
\sum_{i<k_1<\cdots<k_t<j} (-1)^t (\lambda_{k_t} -\lambda_j)
x_{i,k_1}x_{k_1,k_2}\cdot \dots \cdot x_{k_t,j}=0,\nonumber
\end{align}
for every pair $(i,j)$, $i<j$. Thus we have a set of equations defining  the Hessenberg
manifold in a dense subset.
\par
Notice that in the formula above, the coefficients
$\lambda_j -\lambda_i$ and $\lambda_{k_t} -\lambda_j$
are always non--zero,
because the entries of $H$ are all distinct.
This implies that the Jacobian matrix associated
with the set of equations (\ref{polhess}) has
maximal rank, so that ${\rm Hess}_p(\mathbb{R})$
is smooth.
\par
We conclude observing that the formula (\ref{char}) continues to make  sense if the matrix $R$
is taken in a set which is closed under conjugation by elements in P.  Hence we can abstractly
give a more general definition of Hessenberg
manifold by considering matrices of
the form
\\
\vskip4.5truecm
\setlength{\unitlength}{1cm}
\begin{picture}(0,0)\thicklines\put(3,2.5){$R=$}
\put(4.5,0.5){\line(0,1){4}}
\put(4.5,0.5){\line(1,0){4}}
\put(4.5,4.5){\line(1,0){4}}
\put(8.5,0.5){\line(0,1){4}}
\put(5.5,4){\line(0,1){0.5}}
\put(5.5,4){\line(1,0){1.5}}
\put(7,3.5){\line(0,1){0.5}}
\put(7,3.5){\line(1,0){0.5}}
\put(7.5,2){\line(0,1){1.5}}
\put(7.5,2){\line(1,0){1}}
\put(5.15,4.15){$*$}
\put(5.65,3.65){$*$}
\put(6.15,3.65){$*$}
\put(6.65,3.65){$*$}
\put(6.15,3.15){$*$}
\put(6.65,3.15){$*$}
\put(7.15,3.15){$*$}
\put(6.65,2.65){$*$}
\put(7.15,2.65){$*$}
\put(7.15,2.15){$*$}
\put(7.65,1.65){$*$}
\put(8.15,1.65){$*$}
\put(8.15,1.15){$*$}
\put(7.5,3.9){$0$}
\put(8.65,2.5){.}
\end{picture}
\\
By this remark one generalizes the definition of Hessenberg
manifolds to the context of semisimple Lie group, as we show in the
next section.

\section{Real Hessenberg manifolds}
We report the definition and some properties of real Hessenberg  manifolds,
that can be found in~\cite{DP}.
Let G be a connected (real) semisimple noncompact Lie group with finite  center.
Let $\Lie{g}$ be its Lie algebra and $\Sigma$ the
corresponding restricted root system. Choose an ordering in $\Sigma$,
fix the set of positive roots $\Sigma_+$ and  the set of positive
simple roots $\Delta=\{\delta_1,\cdots,\delta_l\}$. Let $\mc{R}$ be
some
proper subset of the set of the positive roots
$\Sigma_+$. We call it of Hessenberg type if it satisfies the
following property:
\vskip0.2cm
\centerline{if $\alpha\in\mc{R}$ and $\beta$ is any negative
root such that $\alpha+\beta\in \Sigma_+$, then $\alpha+\beta\in
\mc{R}$.}
\vskip0.2cm
Denote by $\mc{C}$ the complement in $\Sigma_+$ of $\mc{R}$.
Let $\group{P}$ be the minimal parabolic subgroup of G. We
define the Hessenberg manifold corresponding to
$\mc{R}$ and to some regular element $H$ in the Cartan
subspace
$\Lie{a}$ as the following submanifold of $\group{G}/\group{P}$
:
$${\rm Hess}_\mc{R} (H) =
\{\langle g\rangle_\group{P} \in \group{G}/\group{P}:
\Ad{g}^{-1} H \in \Lie{b}_\mc{R} \} ,$$
where $ \Lie{b}_\mc{R}= \Lie{a} \oplus \Lie{n}
\oplus \bigoplus_{\gamma\in\mc{R}} \Lie{g}_\gamma$.
The Hessenberg manifolds are smooth submanifolds of G/P, and they are
algebraic varietes. The next proposition can be found in [DP], but we
recall the proof bacause it shows the explicit algebraic equations
locally defining ${\rm Hess}_\mc{R} (H)$. The notations are those that  we
introduced in Ch.~\ref{pr}. We just remind that W denotes the
Weyl group and MAN the Langlands decomposition
of P.
\begin{prop}{\rm (}\cite{DP}{\rm )}
${\rm Hess}_\mc{R} (H)$ is a smooth submanifold of $\group{G}/\group{P}$
of dimension $\sum_{\alpha \in \mc{R}} m_\alpha$.
\end{prop}
We prove the smoothness of ${\rm Hess}_\mc{R} (H)$ by writing local
defining equations. The local coordinates are given by the ``Bruhat
charts'' in
$\group{G}/\group{P}$. It is useful in this context to recall in some  details
the Bruhat decomposition of
$\group{G}$ and $\group{G}/\group{P}$, in order to construct
a somewhat canonical open covering of $\group{G}/\group{P}$.
The Bruhat decomposition of $\group{G}$ is the disjoint union
$$\group{G}= \coprod_{w\in W} \group{P}w\group{P}.$$
Observe that
$$
\group{P}w\group{P} =\group{MAN}(wN)\\
     = \group{MA\overline{N}}(w\overline{\group{N}}w^{-1})w\\
     = \group{MA}(\group{\overline{N}\overline{N}}^w)w,
$$
where $\overline{\group{N}}^w := w\overline{\group{N}}w^{-1}$; therefore
$\mc{P}:=\group{P}w\group{P}=\group{MA}(\group{\overline{N}\overline{N}} ^w)w$. In
particular, if $w_0$ is the element of the Weyl group which exchanges
the negative and the positive roots,
then
$\overline{\group{N}}^{w_0} =
{\group{N}}:=
\exp {\Lie{n}}$. Therefore the ``big cell'' is
$$
\mc{P}_{w_0} = \group{MA\overline{N}} {{\group{N}}}w_0\\
              = \group{\overline{N}MA}w_0 \overline{\group{N}}\\
              =w_0 ( {{\group{N}}} \group{MA\overline{N}}).
$$
The Bruhat decomposition of $\group{G}/\group{P}$ is trivially
induced by that of $\group{G}$:
$$\group{G}/\group{P}=\coprod_{w\in W} \langle \mc{P}_w
\rangle_{\group{P}}=\coprod_{w\in W}  \overline{\mc{P}}_w.$$
Notice that
$\overline{\mc{P}}_w=\langle w_0 {{\group{N}}}\rangle_\group{P}$.
Now we set, for $w\in W$,
\begin{align*}
\overline{\group{N}}^{(w)} &= \overline{\group{N}} \cap
w {{\group{N}}}w^{-1}=\exp(\overline{\Lie{n}}\cap \Ad w {\Lie{n}})\\
\overline{\group{N}}_{(w)} &=\overline{\group{N}} \cap  w\overline{\group{N}}w^{-1}=
\exp(\overline{\Lie{n}}\cap\Ad
w\overline{\Lie{n}}).
\end{align*}
Then $\overline{\group{N}}=\overline{\group{N}}^{(w)}  \overline{\group{N}}_{(w)}=
\overline{\group{N}}_{(w)}\overline{\group{N}}^{(w)}$.
\begin{prop}{\rm (}\cite{DP}{\rm )} \label{bruhat}
$\overline{\group{N}}\overline{\group{N}}^w \subset w
{{\group{N}}}\overline{\group{N}}w^{-1}$.
\end{prop}
\begin{proof}
Since  $\overline{\group{N}}_{(w)}\overline{\group{N}}^{w}=(\overline{\group{N} }\cap
w\overline{\group{N}}w^{- 1})\overline{\group{N}}^w=(\overline{\group{N}}\cap
\overline{\group{N}}^w){\group{N}}^w=\group{\overline{N}\overline{N}}^w\ cap\overline{\group{N}}^w $, we have that
\begin{align*}
\overline{\group{N}}^{(w)}  \overline{\group{N}}_{(w)}\overline{\group{N}}^w &=
\overline{\group{N}}^{(w)}(\group{\overline{N}\overline{N}}^w  \cap\overline{\group{N}}^w )\\
&= (\overline{\group{N}}\cap w  {{\group{N}}}  w^{-1})(\group{\overline{N}\overline{N}}^w \cap
\overline{\group{N}}^w )\\ &=  [\overline{\group{N}}(\group{\overline{N}\overline{N}}^w \cap  \overline{\group{N}}^w
)] \cap [w
 {{\group{N}}} w^{-1}(\group{\overline{N}\overline{N}}^w \cap  \overline{\group{N}}^w )]\\
&\subset w  {{\group{N}}} w^{-1} (\group{\overline{N}\overline{N}}^w  \cap \overline{\group{N}}^w )\\
&\subset w  {{\group{N}}} w^{-1} \overline{\group{N}}^w\\
&=w {{\group{N}}} w^{-1} w \overline{\group{N}} w^{-1}\\
&= w {{\group{N}}} \overline{\group{N}} w^{-1}.
\end{align*}
\end{proof}
\begin{cor}{\rm (}\cite{DP}{\rm )}
$\mc{P}_w \subset w {{\group{N}}}  \group{MA\overline{N}}$.
\end{cor}
\begin{proof}
Using Proposition~\ref{bruhat}, we get
\begin{align*}
\mc{P}_w &= \group{MA}(\group{\overline{N}\overline{N}}^w)w\\
&\subset \group{MA}(w {{\group{N}}} \overline{\group{N}} w^{-1})w\\
&= \group{MA}w(w_0 \overline{\group{N}} w_0) \overline{\group{N}}\\
&=w\group{MA} w_0 \overline{\group{N}}w_0 \overline{\group{N}}\\
&=w w_0  \group{MA\overline{N}}w_0 \overline{\group{N}}\\
&=w w_0 \group{\overline{N}MA} w_0 \overline{\group{N}}\\
&=w w_0 \overline{\group{N}} w_0  \group{MA\overline{N}}\\
&= w  {{\group{N}}}  \group{MA\overline{N}}.
\end{align*}
\end{proof}
 From this corollary it follows that the open sets
$$ch(w) :=\langle w  \group{N}\rangle_\mc{P},$$
  give a covering of $\group{G}/\group{P}$, since $ch(w) \supset
 \overline{\mc{P}}_w$. Hence every point in $\group{G}/\group{P}$ can be
written as $\langle w  {n}\rangle_\group{P}$, for some $w\in W$, where
$ {n}\in \group{N}$ is unique (once $w$ is fixed). Let $\langle w
 {n} \rangle_\group{P} \in ch(w)$; then $\langle w  {n}  \rangle_\group{P}
\in {\rm Hess}_\mc{R} (H )$ if and only if $\Ad(w  {n} )^{-1} H \in
\Lie{b}_\mc{R}$. But  $\Ad(w  {n} )^{-1} H = \Ad  {n}^{-1}
(w^{-1}   H )$, and therefore $\langle w  {n} \rangle_\group{P} \in
{\rm Hess}_\mc{R} (H ) \Leftrightarrow \langle {n} \rangle_\group{P} \in
{\rm Hess}_\mc{R} (w^{-1}  H )$. Fix $w=1$. Hence we must impose that
$\Ad  {n}^{-1} H \in \Lie{b}_\mc{R}$, that is,
$$(\Ad  {n}^{-1} H)_\alpha = 0 \hspace{.2in} \forall \alpha \in  \mc{C},$$
where $(X)_\alpha$ denotes the component in $\Lie{g}_\alpha$ of $X$.
We write $ {n} = \exp \nu$, with $\nu = \sum_{\alpha \in \Sigma_+}
\nu_{\alpha} \in  {\Lie{n}}$. Therefore, because of nilpotency,
\begin{align*}
\Ad  {n}^{-1} H &= \Ad \exp (- \nu ) H \\
&=  e^{-\ad \nu} H \\
&= H - [\nu , H ] + \frac{1}{2} [\nu, [\nu,H]]+ \dots
\hspace{.2in}  \text{(finite sum)}
\end{align*}
and, for $\alpha\in\Sigma_+$,
$$(\Ad  {n}^{-1} H )_\alpha = \alpha (H) \nu_\alpha +
\text{(terms containing } x_{\beta , i}\text{, with } {\rm ht} (\beta
) <  {\rm ht}( \alpha)).$$
If the root $\alpha$ has multiplicity $m_\alpha$, we write
$\nu_\alpha = \sum_{j=1}^{m_\alpha} x_{\alpha,j} E_{\alpha,j}$, where
$\{ E_{\alpha,j} \}$ is a basis of $\Lie{g}_\alpha$. This means that
we are using coordinates $\{ x_{\alpha,j}\}$, with $\alpha \in
\Sigma_+$, $j=1, \dots , m_\alpha$, in the chart $ch(1)$.
Consequently, the equations (locally) defining ${\rm Hess}_\mc{R} (H)$  are
\begin{equation}\label{hes}
p_{\alpha,j} (x) = 0,\hspace{.2in} \alpha\in
\mc{C},j=1,\dots,m_\alpha,
\end{equation}
where
$$
p_{\alpha,j} = \alpha(H) x_{\alpha,j} + \text{(terms containing
}x_{\beta,i},\text{ with }{\rm ht} (\beta)  < {\rm ht}
(\alpha)\text{)}.
$$
the components $x_{\gamma,k}$ of vector $x$, with $\gamma\in
\Sigma_+,k=1,\dots,m_\gamma$, are ordered in such a way that
$x_{{\gamma_1}, k_1}$ precedes $x_{{\gamma_2},k_2}$, if ${\rm
ht}({\gamma_1}) < {\rm ht}({\gamma_2})$. Obviously, ${\rm Hess}_\mc{R}
(H)$ is smooth if and only if the matrix
$$J= \left[  \frac{\partial p_{\alpha,j}}{\partial
x_{\gamma,k}}\right]_{\alpha\in\mc{C},\gamma\in\Sigma_+}$$
has maximal rank. On the other hand
$$ \frac{\partial p_{\alpha,j}}{\partial x_{\gamma,k}}=
\begin{cases}
0&\text{if ${\rm ht}(\gamma ) \geq {\rm ht}(\alpha)$ and $\gamma \neq
\alpha$}\\0&\text{if $\gamma=\alpha$ and $k\neq j$}\\
\alpha(H)&\text{if $\gamma=\alpha$ and $k=j$}\\
*&\text{otherwise}.
\end{cases}
$$
The square submatrix
$$\tilde{J} = \left[  \frac{\partial p_{\alpha,j}}{\partial
x_{\gamma,k}} \right]_{\alpha,\gamma\in\mc{C}}$$
is lower triangular, and the blocks on the diagonal are $\alpha (H)
I_{m_\alpha}$. Therefore
$${\rm det} \tilde{J} = \prod_{\alpha\in\mc{C}} \alpha
(H)^{m_\alpha} \neq 0,$$
since $H$ is regular, thereby proving that ${\rm Hess}_\mc{R}(H)$ is
smooth. Finally, we compute the dimension of ${\rm Hess}_\mc{R} (H)$:
$${\rm dim} {\rm Hess}_\mc{R} (H) ={\rm dim} (\group{G}/\group{P}) -  {\rm
rank} J = \sum_{\alpha\in\Sigma_+} m_\alpha - \sum_{\alpha\in\mc{C}}
m_\alpha= \sum_{\alpha\in\mc{R}} m_\alpha.$$
Finally, if $H$ is a regular element, then also $w^{-1}
H$ is so. Hence the conclusions above are true in any chart $ch(w)$.
\\
\\
\begin{prop}\label{ideal}
 The vector space
$\displaystyle{
{\Lie{n}}_\mc{C}=\bigoplus_{\alpha\in\mc{C}} \Lie{g}_\alpha}$
is an ideal in $ {\Lie{n}}$.
\end{prop}
\begin{proof}
Let $\alpha\in\mc{C}$ and
$\beta\in\Sigma_+$. If $\alpha+\beta$ is a root in $\mc{R}$, then
$\alpha+\beta-\alpha\in\mc{R}$, that is false. Therefore, if
$X\in\Lie{g}_\alpha$ and $Y\in\Lie{g}_\beta$, then
$[X,Y]\in\Lie{g}_{\alpha +\beta}$, with $\alpha +\beta\in\mc{C}$.
\end{proof}


\chapter{Multicontact vector fields on Hessenberg  manifolds}\label{multi3}
\vskip0.5cm
This is the main chapter of the thesis. First, we transfer the
multicontact structure from G/P to the Hessenberg manifolds
${\rm Hess}_{\mc{R}}(H)$. Then, in the second section, we lift the  problem
to the Lie algebra level, that is, we define the notion of multicontact
vector field on ${\rm Hess}_{\mc{R}}(H)$, and we consider the Lie  algebra
$MC$(S) of all  multicontact vector fields on a local model S of ${\rm
Hess}_{\mc{R}}(H)$, the central object of our investigation. We prove
that $MC$(S) contains canonically a quotient algebra
$\Lie{q}/\Lie{n}_{\mc{C}}$. In the third
section, we show that this quotient exhausts all multicontact
vector fields, if some natural additional hypotheses on the Hessenberg
manifolds are assumed,  synthetized in the notion of Iwasawa
sub-models. In the last section we consider an example that shows that
$MC$(S) can be larger than
$\Lie{q}/\Lie{n}_{\mc{C}}$, if the additional assumptions do not hold.
\section{A multicontact structure on Hessenberg manifolds}
The Iwasawa Lie algebra ${\Lie{n}}$ has a multistratification
$ {\Lie{n}}=\sum_{\gamma\in\Sigma_+} \Lie{g}_\gamma$ and a
stratification by height $ {\Lie{n}}= {\Lie{n}}_1 \oplus
     {\Lie{n}}_2 \oplus \cdots \oplus  {\Lie{n}}_h$, where
$ {\Lie{n}}_i$ is the direct sum of all root spaces
$\Lie{g}_\gamma$ that are sum of $i$ positive simple roots, that is
${\rm ht}(\gamma)=i$. In particular $[ {\Lie{n}}_i ,
     {\Lie{n}}_j ] \subset  {\Lie{n}}_{i+j}$. The structure of
$ {\Lie{n}}$, viewed as the tangent space to $ {\group{N}}$ at the
identity, allows us to give generalized versions of contact
mappings (see Chap.1).
\par
We ask ourselves how to relate with $ {\Lie{n}}$ the tangent space to  some point of
${\rm Hess}_{\mc{R}} (H)$  because we want to
transfer the  stratification
to the tangent bundle on some open  set of ${\rm Hess}_\mc{R} (H)$.
We have seen that a choice of a Hessenberg structure $\mc{R}$
determines the set of independent local coordinates
$\{x_{\alpha,j}: \alpha\in\mc{R},1\leq j \leq m_\alpha\}$ on the
Hessenberg manifold, whereas the remaining
entries $\{x_{\alpha,j}: \alpha\in\mc{C},1\leq
j \leq m_\alpha\}$ are polynomial functions of
the previous ones (see~\eqref{hes}).
The coefficients of the polynomials depend on $H$ and more is true:
those that are not
zero are in fact given by functions that never vanish on the set of
regular elements in
$\Lie{a}$. Thus, the slice
$\group{S}$ of $ {\group{N}}$ obtained by setting
$x_{\alpha,j}=0$ if $\alpha\in\mc{C}$ is diffeomorphic to ${\rm
Hess}_\mc{R} (H)$ for every regular element $H$. The graph mapping
\begin{equation} \label{graph}
\phi:( \{x_{\beta,k}\}_{\beta\in\mc{R}},0)\longmapsto
(\{x_{\beta,k}\}_{\beta\in\mc{R}} , \{p_{\alpha,j}
(x_{\beta,k})\}_{\alpha\in\mc{C}})
\end{equation}
gives the diffeomorphism.
For this reason we shall use $\group{S}$ as a simplified model for ${\rm
Hess}_{\mc{R}} (H)$.
\par
In $\group{SL}(n,\R)$, where the nilpotent subgroup N is given by the
unipotent upper
triangular matrices, the diffeomorphism can be visualized as follows:
\vskip5truecm
\setlength{\unitlength}{1cm}
\begin{picture}(0,0)\thicklines\put(2,2.5){$\{*_{i,j}:(i,j)\in\mc{R}\}\mapsto$}
\put(6,0.5){\line(0,1){4}}
\put(6,0.5){\line(1,0){4}}
\put(6,4.5){\line(1,0){4}}
\put(10,0.5){\line(0,1){4}}
\put(7,4){\line(1,0){1.5}}
\put(8.5,3.5){\line(0,1){0.5}}
\put(8.5,3.5){\line(1,0){0.5}}
\put(9,2){\line(0,1){1.5}}
\put(9,2){\line(1,0){0.5}}
\put(9.5,1.5){\line(0,1){0.5}}
\put(7.15,3.65){$*$}
\put(7.65,3.65){$*$}
\put(8.15,3.65){$*$}
\put(7.65,3.15){$*$}
\put(8.15,3.15){$*$}
\put(8.65,3.15){$*$}
\put(8.15,2.65){$*$}
\put(8.65,2.65){$*$}
\put(8.65,2.15){$*$}
\put(9.15,1.65){$*$}
\put(9,3.9){$p(*)$}
\put(6.2,4.){$1$}
\put(6.65,3.6){$1$}
\put(7.15,3.1){$1$}
\put(7.65,2.6){$1$}
\put(8.15,2.1){$1$}
\put(8.65,1.6){$1$}
\put(9.15,1.1){$1$}
\put(9.6,0.66){$1$}
\end{picture}

Next we define the contact structure on $\group{S}$ and
 we show how it can be transferred to ${\rm Hess}_\mc{R} (H)$.
First we define the
left--invariant  vector fields
$X_{\alpha,j}$ on
$ {\group{N}}$ that concide with the partial derivative operators at the
origin. Next we obtain the differentiable
structure on $\group{S}$ by considering the projections
$\overline{X}_{\alpha,j}$ on
$\group{S}$ of those vector fields that correspond to $\alpha\in\mc{R}$.
Finally, the push-forwards $\phi_*
( {X}_{\alpha,j})$ will define the differentiable structure on the
Hessenberg manifold. This structure will allow us to give a generalized  version
of contact mapping.
\subsection{Multicontact mappings}
Consider the basis $\{X_{\alpha,j} : \alpha\in\Sigma_+,1\leq j\leq
m_\alpha\}$ of
left-invariant vector fields on
$ {\group{N}}$, where
$$X_{\alpha,j} ( {n}) = (l_{ {n}})_{*e}
\frac{\partial}{\partial x_{\alpha,j}}{\Big|_e},$$
and write
$$X_{\alpha,j} =
\sum_{\gamma\in\Sigma_+}\sum_{k=1}^{m_\gamma} a_{\gamma,k}^{\alpha,j}
\frac{\partial}{\partial x_{\gamma,k}},$$ where
$a_{\gamma,k}^{\alpha,j}$ are some smooth
functions on $ {\group{N}}$.
We want to compute their explicit expressions in exponential
coordinates. In order to do this, we recall a result that can be
found in \cite{CDKR2}. If $\alpha=\sum_{\delta\in\Delta}a_\delta \delta$
and $\beta  =\sum_{\delta\in\Delta}
b_\delta \delta$ are two positive roots, we write $\alpha
\preceq
\beta$ if
$a_\delta
\leq b_\delta$ for all $\delta\in\Delta$.
     We say that ${\alpha}_1+\cdots + \alpha_n$ is a $chain$ if each
$\alpha_j$ and
each partial sum $\alpha_1 +\cdots +\alpha_j$ is a root for all
$j=1,\dots,n$. Ordered pairs of roots can be joined by chains:
\begin{lemma}[\cite{CDKR2}] \label{chain}
Let $\alpha$ and $\beta$ be distinct positive roots and suppose that
$\alpha\succeq
\beta$. Then there exist simple roots $\delta_1,\dots,\delta_p$ such
that $\alpha=\beta+\delta_1 + \dots +\delta_p$ is a chain.
\end{lemma}
\noindent
We can now describe the coefficient functions
$a_{\gamma,k}^{\alpha,j}$.
\begin{lemma} \label{left}
For every root $\alpha\in\Sigma_+$ and $j=1,\dots,m_\alpha$ we have
\begin{equation} \label{componenti}
a_{\gamma,k}^{\alpha,j}=
\begin{cases}
0&\text{if ${\rm ht}(\alpha)\geq{\rm ht}(\gamma)$ and  $\alpha\neq\gamma$}\\
0&\text{if $\alpha=\gamma$ and $k\neq j$}\\
1&\text{if $\alpha=\gamma$ and $k=j$}\\
P&\text{if ${\rm ht}(\alpha)<{\rm ht}(\gamma)$},
\end{cases}
\end{equation}
where $P$ is a polynomial that does not vanish only if
$\alpha\preceq\gamma$. In this case, it depends only on those
variables labeled by those roots
$\alpha_1,\cdots,\alpha_q$ for which
$\alpha +
\alpha_1 +\cdots +\alpha_q =\gamma$ is a chain.
     This implies that
\begin{equation} \label{c}
X_{\alpha,j}=\sum_{\gamma\in\mc{C}}\sum_{k=1}^{m_\gamma}
a_{\gamma,k}^{\alpha,j} \frac{\partial}{\partial x_{\gamma,k}},
\end{equation}
for every $\alpha\in\mc{C}$.
\end{lemma}
\begin{proof}
Consider $\{E_{\beta,i}:\beta\in\Sigma_+,1\leq i\leq m_\beta\}$ as a
basis of $\Lie{n}$ viewed as the tangent space to N at the identity
and write
$$
{n}  = \exp \left( \sum_{\epsilon \in \Sigma_+}
\sum_{s=1}^{m_{\epsilon}} y_{\epsilon,s}
E_{\epsilon,s}\right) \hskip0.5cm  {n}^\prime  =
\exp\left( \sum_{\beta\in\Sigma_+}
\sum_{r=1}^{m_{\beta}} x_{\beta , r} E_{\beta,r}\right).
$$
Let $f$ be a smooth function on N. From the left invariance
$$(X_{\alpha,j} f)( {n})= (l_{ {n}})_{*e}
\frac{\partial}{\partial x_{\alpha ,j}}\Big|_e
f=\frac{\partial}{\partial x_{\alpha, j}}f\circ l_{ {n}} (e)$$
it follows that the component $a_{\gamma,k}^{\alpha,j}$ is given as
the derivative with respect to $x_{\alpha,j}$ of $f\circ l_n$,
where $f$ is
the coordinate function $n^\prime \mapsto
(n^\prime)_{(\gamma,k)}$. That is,
$$
a_{\gamma,k}^{\alpha,j}(n')= \frac{\partial}{\partial
x_{\alpha,j}} ( {n} {n}^\prime)_{(\gamma,k)}.
$$
An explicit calculation using the Campbell-Hausdorff formula gives
\begin{eqnarray}
     {n}  {n}^\prime &=& \exp \left( \sum_{\epsilon \in \Sigma_+}
\sum_{s=1}^{m_{\epsilon}} y_{\epsilon,s} E_{\epsilon,s}\right)
\exp\left(
\sum_{\beta\in\Sigma_+} \sum_{r=1}^{m_{\beta}} x_{\beta , r}
E_{\beta,r}\right) \nonumber \\
&=&\exp\left( \sum_{\epsilon \in
\Sigma_+} \sum_{s=1}^{m_{\epsilon}} y_{\epsilon,s} E_{\epsilon,s} +
\sum_{\beta\in\Sigma_+} \sum_{r=1}^{m_{\beta}}
x_{\beta , r} E_{\beta,r} +  \right. \nonumber \\
     && \hskip1cm + \left. \frac{1}{2}\sum_{\epsilon \in \Sigma_+}
\sum_{s=1}^{m_{\epsilon}}\sum_{\beta\in\Sigma_+}
\sum_{r=1}^{m_{\beta}}y_{\epsilon,s}x_{\beta , r} [E_{\epsilon,s},  E_{\beta,r}]
     +\cdots \right) .\label{ch}
\end{eqnarray}
Since $\Lie{n}$ is nilpotent, the sum in (\ref{ch}) is finite, and the  variable
$x_{\alpha,j}$ appears in the coefficient of an iterated bracket of the  form
$$
[E_{\epsilon_i,s},[\dots,[E_{\epsilon_1,s},[E_{\beta_k,r},\dots,[E_{\beta_1,r},
E_{\alpha,j}]\dots ]]\dots],
$$
     provided that the bracket is not zero.
If it appears with a power greater than or equal to two, then
the derivative of the corresponding monomial with respect to
$x_{\alpha,j}$ evaluated at the identity is zero. Thus, the
only relevant brackets are those of the form
$$
[E_{\epsilon_i,s},[\dots,[E_{\epsilon_1,s},
E_{\alpha,j}]\dots ]],
$$
where $\alpha + \eps_1 +\dots +\eps_i =\gamma$ is a chain. The
coefficient of such an iterated  bracket is the monomial
$y_{\eps_i,s}\dots y_{\eps_1,s}x_{\alpha,j}$. Its derivative
with respect to $x_{\alpha,j}$ evaluated at the identity is
$y_{\eps_i,s}\dots y_{\eps_1,s}$. This proves
\eqref{componenti}.

Finally, by
definition of $\mc{R}$, if $\alpha\in\mc{C}$ then $\alpha +
\delta \in \mc{C}$, for any simple root
$\delta$ such that $\alpha +\delta$ is a root. Hence, if  $\alpha\in\mc{C}$,
there are no chains going from $\alpha$ to a root $\gamma\in\mc{R}$.
Thus also (\ref{c})
follows.
\end{proof}
For every $\alpha\in\Sigma_+$, and $1\leq j\leq m_\alpha$, consider the  vector
field $\overline{X}_{\alpha,j}$ whose $(\gamma,k)$ component is
$$r_{\gamma,k}^{\alpha,j} =
\begin{cases}
a_{\gamma,k}^{\alpha,j} &\text{ if }\gamma\in\mc{R}\text{ and }
k=1,\cdots,m_\gamma\\
0&{\rm otherwise.}
\end{cases}
$$
The $\overline{X}_{\alpha,j}$ are vector fields on $\group{S}$,
and from (\ref{c}) $\overline{X}_{\alpha,j}=0$ for
every
$\alpha\in\mc{C}$. Moreover, (\ref{componenti}) implies that the set
$\{\overline{X}_{\alpha,j}:{\alpha\in\mc{R}},{j=1,\cdots,m_\alpha}\}$
is a basis of the tangent space at
any point of $\group{S}$. Indeed, writing the matrix of the
coefficients of
$\{\overline{X}_{\alpha,j}:{\alpha\in\mc{R}},{j=1,\cdots,m_\alpha}\}$,
ordering the roots according to any lexicographic order, we obtain
a triangular matrix with ones along the diagonal. Hence
$\{\phi_*  (\overline{X}_{\alpha,j}):{\alpha\in\mc{R}},{j=1,\cdots,m_\alpha}\}$
is a basis of the tangent space at all points of an open
set of
${\rm Hess}_{\mc{R}} (H)$.
\par
Next, we introduce some special sub-bundles of the tangent bundle of
$\group{S}$. For
$\delta\in\Delta\cap\mc{R}$, put
$\overline{\Lie{g}}_\delta =
{\rm span}\{\overline{X}_{\delta,i}:{i=1,\cdots,m_\delta}\}$.
   From Proposition~\ref{ristretto} below it follows that the
vector fields in the family
$\{\overline{\Lie{g}}_{\delta}\}_{\delta\in\Delta_{\mc{R}}}$,
$\Delta_\mc{R}=\Delta\cap\mc{R}$, satisfy a H\"ormander-type
condition: their iterated
brackets generate at each point the tangent space of
$\group{S}$. We denote by $\mathfrak{X}(\group{N})$ the Lie
algebra of all smooth vector fields on N.
\begin{prop}\label{ristretto}
Given $X$ and $Y \in \mathfrak{X}( {\group{N}})$, the following
formula holds at every point
$n\in\group{S}$
$$[\overline{X},\overline{Y}]( {n})=\overline{[X,Y]} ( {n})
.$$
\end{prop}
\begin{proof}
Let $X$ and $Y \in \mathfrak{X}( {\group{N}})$ and write
$X=\overline{X} +
\underline{X}$, where
$$
\overline{X} := \sum_{\beta\in\mc{R}}\sum_{i=1}^{m_\beta}
r_{\beta,i} \frac{\partial}{\partial x_{\beta,i}},
\qquad
\underline{X} := \sum_{\gamma\in\mc{C}}\sum_{k=1}^{m_\gamma}
c_{\gamma,k} \frac{\partial}{\partial
x_{\gamma,k}},
$$
and similarly $Y=\overline{Y}+\underline{Y}$.
Then
$$\overline{[X,Y]} ( {n})= \overline{[\overline{X},\overline{Y}]}( {n})  +
\overline{[\overline{X},\underline{Y}]} ( {n}) +
\overline{[\underline{X},
\overline{Y}]} ( {n}) + \overline{[\underline{X}, \underline{Y}]}
( {n}).$$ Clearly
$\overline{[\overline{X},\overline{Y}]}=[\overline{X},\overline{Y}]$.
Moreover
$$
\overline{[\overline{X},\underline{Y}]}=\overline{[\underline{X},\overline{Y}]}=0,
$$
because when expanded in terms of partial derivatives, each of the
above brackets
contains only coefficients of the form
$({\partial}/{\partial x_{\gamma,k}})r_{\beta,i}$, which vanish whenever
$\gamma\in\mc{C}$ and $\beta\in\mc{R}$ because of
(\ref{componenti}). Finally,
$\overline{[\underline{X}, \underline{Y}]}=0$, because in $  [\underline{X},
\underline{Y}]$ only the coefficients of components labeled by
$\mc{C}$ will  appear, but they
become zero once we project them on $\group{S}$.
\end{proof}
Let $\mc{A},\mc{B}$ be some open subsets of ${\rm Hess}_{\mc{R}} (H)$.
Without loss of generality, we can assume
$\mc{A},\mc{B}
\subset
\left({\group{N}}\cap{\rm Hess}_{\mc{R}} (H)\right)$.
Let $f:\mc{A}\rightarrow\mc{B}$ be a diffeomorphism. We say that
$f$ is a $multicontact$ map if
$$
f_* (\phi_*(\overline{\Lie{g}}_\delta))\subseteq
\phi_*(\overline{\Lie{g}}_\delta)
$$
for every simple root
$\delta$ in
$\mc{R}$,
where $\phi$ is the graph mapping defined in (\ref{graph}).
\subsection{Example}\label{ex}
Consider the simple Lie algebra $\Lie{g}=\Lie{sl}(4,\mathbb{R})$. The
Cartan subspace
$\Lie{a}$ is the abelian algebra of diagonal matrices. Denote by
${\rm  diag}(a,b,c,d)$
the diagonal matrix with entries $a,b,c,d$. The standard simple
restricted roots are
$\alpha,\beta,\gamma$, where $\alpha({\rm diag}(a,b,c,d))=(a-b)$,  $\beta({\rm
diag}(a,b,c,d))=(b-c)$ and $\gamma({\rm diag}(a,b,c,d))=c-d$. Let
$\mc{R}$ be the set of all positive roots except the highest
one, that is $\alpha+\beta+\gamma$. A natural regular element in
$\Lie{a}$ is
$$H=
\bmatrix
-1&0&0&0\\ 0&\frac{1}{2}&0&0\\ 0&0&-\frac{1}{2}&0\\ 0&0&0&1
\endbmatrix.
$$
    Consider then
$\group{G}=\group{SL}(4,\mathbb{R})$
and its minimal parabolic subgroup consisting of lower triangular
matrices. Since the
Hessenberg manifold is a submanifold of
$\group{G}/\group{P}$
and the problem of multicontact mappings is local, we restrict
ourselves to the big
cell
$\group{{N}}\subset \group{G}/\group{P}$ and fix coordinates on it:
\begin{equation} \label{coordinate}
{n}=
\bmatrix
1&x&u&z\\ 0&1&y&v\\ 0&0&1&t\\ 0&0&0&1
\endbmatrix.
\end{equation}
Thus
$$
\Ad{{n}^{-1}} H =
\bmatrix
-1&{3}/{2}x&({u-3xy})/{2}&2z-({3vx-ut+3xyt})/{2}\\
0&{1}/{2}&-y&({v+yt})/{2}\\
0&0&-{1}/{2}&{3}/{2}t\\
0&0&0&1
\endbmatrix.
$$
The points $n(x,y,t,u,v,z)\in\group{{N}}$
that lie in the
Hessenberg manifold are those that satisfy
$2z-({3vx-ut+3xyt})/{2}=0$.
Therefore, the slice $\group{S}$ is the algebraic submanifold of
N defined by the equation $z=0$. A basis of
$\Lie{n}$ is given by the following left invariant vector fields
\begin{align*}
X&=\frac{\partial}{\partial x}+y\frac{\partial}{\partial u} +
v\frac{\partial}{\partial z}  &   U&=\frac{\partial}{\partial
u}+t\frac{\partial}{\partial z}   &   Z&=\frac{\partial}{\partial z},\\
Y&=\frac{\partial}{\partial y} + t
\frac{\partial}{\partial v}    &  V&=\frac{\partial}{\partial v},\\
T&=\frac{\partial}{\partial t},
\end{align*}
with nonzero brackets $[X,Y]=-U$, $[Y,T]=-V$, $[U,T]=-Z$
and
$[X,V]=-Z$. By projecting $\partial / \partial z$ to zero, we
restrict the above vector
fields to a pointwise basis for the tangent space to
$\group{S}$:
\begin{align*}
\overline{X}&=\frac{\partial}{\partial
x}+y\frac{\partial}{\partial u},  &
\overline{U}&=\frac{\partial}{\partial u},\\
\overline{Y}&=\frac{\partial}{\partial
y} + t
\frac{\partial}{\partial v},     &
\overline{V}&=\frac{\partial}{\partial v},\\
\overline{T}&=\frac{\partial}{\partial t}.
\end{align*}
The nonzero brackets are $[\overline{X},\overline{Y}]=-\overline{U}$
and
$[\overline{Y},\overline{T}]=-\overline{V}$.
A diffeomorphism $f$ on some open set in $\group{S}$
is a multicontact mapping if its differential $f_*$ preserves each of
the  following subspaces
$$
\overline{\Lie{g}}_\alpha = {\rm span}\{\overline{X}\}\qquad
\overline{\Lie{g}}_\beta =  {\rm span}\{\overline{Y}\}\qquad
\overline{\Lie{g}}_\gamma =  {\rm span}\{\overline{T}\}.
$$
If we know the multicontact mappings on $\group{S}$ we also know
those on the corresponding Hessenberg manifold, because of the
diffeomorphism $\phi$. For completeness, we compute the
stratified structure in the tangent space of the Hessenberg
manifold. Since
$$
\phi(x,y,t,u,v,0)=(x,y,t,u,v,({3vx-ut+3xyt})/{4}),
$$
for
a point $p=\phi^{-1}(q)$ in $\group{S}$ we have
\begin{align*}
(\phi_*(\overline{X}))g(q)&=\overline{X}(g\circ \phi)(p)\\
&=\overline{X}g\left( x,y,t,u,v,\frac{(3vx-ut+3xyt)}{4}\right)\\
&=\frac{\partial}{\partial x}g(q) +y\frac{\partial}{\partial u}g(q)
+\frac{3v+2yt}{4}\frac{\partial}{\partial z}g(q),
\end{align*}
for any smooth function $g$ on ${\rm Hess}_\mc{R}(H)$.
Thus
$$\phi_*(\overline{X})=\frac{\partial}{\partial x}
+y\frac{\partial}{\partial u}
+\frac{3v-4yt}{4}\frac{\partial}{\partial z}
$$
and similarly
\begin{align*}
\phi_*(\overline{Y})&=\frac{\partial}{\partial
y} + t\frac{\partial}{\partial v},\\
\phi_*(\overline{T})&=\frac{\partial}{\partial t}
+\frac{(3xy-u)}{4}\frac{\partial}{\partial z}.
\end{align*}
In \cite{CDKR1} and \cite{CDKR2}, the authors study the group of  multicontact
mappings on the
boundaries G/P. Their approach is to lift the problem to the
Lie algebra level.
In fact, most of the effort consists in investigating the set of
multicontact vector
fields, that is vector fields whose local flow is given by
multicontact mappings.
The set of multicontact vector fields is a Lie algebra with respect to  the Lie
bracket.
In particular, they ask whether this algebra is finite dimensional and
which is the
resulting algebra. The integration step to the group level is
then relatively simple, and it uses classical tools of Lie theory.
\par
In the case we are studying, we use the same approach and we ask
ourselves if the Lie
algebra of multicontact vector fields on $\group{S}$ is finite
dimensional and how we can
characterize this algebra. As we shall see in the next section, the
multicontact condition
for a vector field $F$ is equivalent to the fact that ${\rm ad}F$
preserves each
$\overline{\Lie{g}}_\delta$, $\delta\in\Delta_\mc{R}$.
In the present example
this means
$$
[F,\overline{X}]=\lambda \overline{X}\qquad
[F,\overline{Y}]=\mu \overline{Y} \qquad
[F,\overline{T}]=\nu \overline{T} ,
$$
for some smooth functions $\lambda$, $\mu$ and $\nu$ on $\group{S}$.
Put $F= f_x \overline{X} + f_y \overline{Y} + f_t \overline{T} + f_u
\overline{U}
+f_v
\overline{V}$. The equation $[F,\overline{X}]=\lambda
\overline{X}$ gives
\begin{align*}
[f_x \overline{X} + f_y \overline{Y} + f_t \overline{T} + f_u  \overline{U}
+f_v\overline{V},\overline{X}]&=  -\overline{X}(f_x)\overline{X} +f_y
\overline{U}
-\overline{X}(f_y)\overline{Y}\\
&  \hskip0.5cm      -\overline{X}(f_t)\overline{T}
-\overline{X}(f_u)\overline{U}-\overline{X}(f_v)\overline{V}\\
&= \lambda \overline{X},
\end{align*}
whence
$$
\begin{cases}
\overline{X}(f_x)&\! =-\lambda \\
\overline{X}(f_y)&\! =\overline{X}(f_t)=\overline{X}(f_v)=0 \\
\overline{X}(f_u)&\! =f_y.
\end{cases}
$$
Similarly, $[F,\overline{Y}]=\mu \overline{Y}$ yields
$$
\begin{cases}
\overline{Y}(f_y)&=-\mu  \\
\overline{Y}(f_x)&=\overline{Y}(f_t)=0  \\
\overline{Y}(f_u)&=-f_x    \\
\overline{Y}(f_v)&=f_t
\end{cases}
$$
and $[F,\overline{T}]=\nu \overline{T}$ implies
$$
\begin{cases}
\overline{T}(f_t)&= -\nu \\
\overline{T}(f_x)&=\overline{T}(f_y)=\overline{T}(f_u)=0   \\
\overline{T}(f_v)&=-f_y.
\end{cases}
$$
The equations $\overline{X}(f_u)=f_y$, $\overline{Y}(f_u)=-f_x$ and
$\overline{Y}(f_v)=f_t$ show that the coefficients $f_u$ and $f_v$
determine all the
others. The equations involving $f_u$ alone may be viewed as the pair
of systems
\begin{equation} \label{model1}
\begin{cases}
\overline{X}^2 f_u &=0   \\  \overline{Y}^2 f_u &=0
\end{cases}
\hskip2cm
\begin{cases}
\overline{T}f_u &=0  \\  \overline{T}\overline{Y}f_u &=0,
\end{cases}
\end{equation}
whereas for $f_v$ we have
\begin{equation} \label{model2}
\begin{cases}
\overline{Y}^2 f_v&=0    \\ \overline{T}^2 f_v&=0
\end{cases}
\hskip2cm
\begin{cases}
\overline{X} f_v&=0 \\    \overline{X}\overline{Y}f_v&=0.
\end{cases}
\end{equation}
Finally, $f_u$ and  $f_v$ are linked by the extra
cross-condition
\begin{equation} \label{comp}
\overline{X}f_u=\overline{T}f_v .
\end{equation}
The set of equations above are typical of the problem of multicontact
vector fields,
and they are related to the multicontact-type equations in the case
G/P studied in
\cite{CDKR1} and \cite{CDKR2}, in a sense that we show below.
\par
Look first at the systems (\ref{model1}). The second system shows
that $f_u$ is  independent of the variables $v$ and $t$. Indeed,
$$
\overline{T}f_u
=\frac{\partial}{\partial t}f_u=0
$$
implies $f_u = f_u (x,y,u,v)
$, and
$$\overline{T}\overline{Y} f_u=
\frac{\partial}{\partial t}
\frac{\partial}{\partial y}f_u
+\frac{\partial}{\partial v}f_u
+ t\frac{\partial}{\partial t}
\frac{\partial}{\partial v}f_u
=\frac{\partial}{\partial v}f_u=0
$$
implies $f_u=f_u(x,y,u)$. The systems (\ref{model1}) are then
equivalent to the system~\eqref{sleq} of Chap.~\ref{po} when  investigating multicontact vector fields on
the nilpotent
Iwasawa subgroup of
$\group{SL}(3,\mathbb{R})$.
It can be integrated and yields the following explicit
polynomial solution
$$f_u = a_0 + a_1 x +a_2 y+a_3 u + a_4 xy+ a_5 x(u-xy) + c_6 uy + a_7
u(u-xy).$$
The same argument holds for (\ref{model2}) and one obtains
$$f_v= b_0 + b_1 y + b_2 t +b_3 v + b_4 yt + b_5 y(v-yt) +b_6 tv +
b_7 v(v-yt).$$
Thus, we may view the study of multicontact vector fields on the slice
$$\group{S}=\Bigl\{
\bmatrix
1&x&u&0\\
0&1&y&v\\
0&0&1&t\\
0&0&0&1
\endbmatrix
: x,y,t,u,v\in\mathbb{R} \Bigr\}
$$
as the study of two embedded models corresponding to
$\group{SL}(3,\mathbb{R})$
that must satisfy the additional compatibility condition given by
(\ref{comp}).
Using this last condition, the polynomials
$f_v$ and $f_u$ become
\begin{align}
f_v&=b_0 +b_1 y + b_2 t + b_3 v + b_4 yt + b_5 y(v-yt) \label{fv}\\
f_u&= a_0 - b_2 x + a_2 y - (a_4 + b_4 ) u +a_4 xy +b_5 uy. \label{fu}
\end{align}
     The Lie algebra of multicontact vector fields has dimension nine,
that is the
number of free costants appearing in the expressions of $f_u$ and
$f_v$ because to each  possible pair
$(f_u,f_v)$ there corresponds one and only one multicontact vector
field. This answers to the question of finite dimensionality.
\par
We now want to characterize the Lie
algebra of multicontact vector fields.
In order to identify this Lie algebra, we refer again to \cite{CDKR1}.
Indeed,
the problem we are considering is nothing but the study of multicontact
mappings on N projected to the slice that corresponds to setting
one coordinate  equal to zero.
Recall that the special subbundles defining the multicontact
structure are obtained by
projecting onto $\group{S}$ those vector fields on N that
correspond to  simple roots, and
that a multicontact mapping is by definition a map that preserves
these bundles. Thus, we
consider the multicontact vector fields on N, project them to be  tangent to
$\group{S}$ and finally check whether we obtain multicontact
vector fields on
$\group{S}$.
\par
By \cite{CDKR1}, the multicontact vector fields on N are all of the form
$\tau(X)$ for some $X\in\Lie{g}$,  where
$$
\tau(X)f(n)=\frac{d}{ds} f({\rm exp}(-sX)n)\Big|_{s=0},
$$
and
$n\mapsto {\rm exp}(-sX)n$ is the action of
${\rm exp}(-sX)\in\group{G}$ on G/P, restricted to N, namely ${\rm
exp}(-sX)n$ is the N-component of the product in the Bruhat
decomposition of G/P.
\par
Consider a vector generating the root space $\Lie{g}_\alpha$, namely  $E_{1,2}$,
where
$E_{i,j}$ denotes the matrix in $\Lie{sl}(4,\mathbb{R})$ with $1$ in
the $(i,j)$
position and zero elsewhere. Therefore, using the coordinates  introduced in
(\ref{coordinate})
$${\rm exp}(-sE_{1,2})n=\bmatrix
1&-s&0&0\\
0&1&0&0\\
0&0&1&0\\
0&0&0&1
\endbmatrix
n =
\bmatrix
1&x-s&u&z\\
0&1&y&v\\
0&0&1&t\\
0&0&0&1
\endbmatrix,
$$
whence
$$\tau(E_{1,2})f(n)=-\frac{\partial}{\partial x}f(n)=-Xf(n) + yUf(n)
+ (v-yt)Zf(n).$$
The projected vector field obtained by setting $\partial/\partial z=0$,  namely
$$\overline{\tau(E_{1,2})}=-\overline{X}+y\overline{U},$$
     is tangent to $\group{S}$ at each point by construction.
     Moreover, the components $(f_u,f_v)=(y,0)$ satisfy equations
(\ref{model1}),
(\ref{model2}) and (\ref{comp}), so that
$\overline{\tau(E_{1,2})}$ is a
multicontact vector field on $\group{S}$. If we do the same
calculation  for each element in
the  canonical basis of
$\Lie{sl}(4,\mathbb{R})$, we obtain polynomials of the form (\ref{fv})  and
(\ref{fu}) for all matrices of the form
$$
\bmatrix
*&*&*&*\\
0&*&*&*\\
0&*&*&*\\
0&0&0&*
\endbmatrix.
$$
In particular we see that $\overline{\tau(E_{1,4})}=0$, and no other  matrix in
$\Lie{sl}(4,\mathbb{R})$ of the form described above gives a zero
vector field. Therefore
the Lie algebra of multicontact vector fields on $\group{S}$ is
isomorphic to
$\Lie{q}/\Lie{n}_\mc{C}$, where
$$
\Lie{q}={\rm span}\Bigl\{
\bmatrix
*&*&*&*\\
0&*&*&*\\
0&*&*&*\\
0&0&0&*
\endbmatrix
: *\in\mathbb{R}\Bigr\}\cap\Lie{sl}(4,\mathbb{R})
$$
and
$$
\Lie{n}_\mc{C}={\rm span}\Bigl\{
\bmatrix
0&0&0&*\\
0&0&0&0\\
0&0&0&0\\
0&0&0&0
\endbmatrix
: *\in\mathbb{R}\Bigr\}.
$$
\par
The Lie algebra $\Lie{q}/\Lie{n}_\mc{C}$ can be viewed as the union
of the two
submodels of $\Lie{sl}$-type plus the reflection of their intersection.
In the subsequent sections we
study the differential equations defining the multicontact vector
fields on an arbitrary
Hessenberg manifold.
On the one hand we shall show that $\Lie{q}/\Lie{n}_\mc{C}$ always
defines a set of multicontact vector fields. On the other hand, we
shall see that  the converse of this
statement is not true in general, that is the Lie algebra of
multicontact vector
fields can be bigger than $\Lie{q}/\Lie{n}_\mc{C}$. It does, however,
become true under the additional hypothesis that the Hessenberg
structure encodes a finite number of embedded models (i.e. each
corresponding to an Iwasawa nilpotent group N) that
intersect non trivially. This is exactly what happens in the
case study that we have just discussed.
\section{A set of multicontact vector fields}
In this section we consider the infinitesimal version of the notion
of multicontact mapping. We obtain a Lie algebra of  multicontact
vector fields  and we address the problem of
understanding its structure.

\subsection{Lifting the multicontact conditions to the infinitesimal  level}
Since the composition and the inverse of multicontact maps on a
Hessenberg manifold are multicontact maps,
the set of such maps is a group.
Moreover, as we already noticed, all Hessenberg manifolds
corresponding to different
choices of regular $H$ give rise to the same
slice $\group{S}$, so they are mutually diffeomorphic. Thus
the group of multicontact maps does not depend
on $H$, so that from now on we focus our attention on the slice
$\group{S}$ of
$\group{N}$. Fix an open set $\mc{A}$ of $\group{S}$. In order
to characterize the group of multicontact
maps, we lift the problem to the Lie algebra level, by
considering multicontact vector fields, that is, vector fields $F$ on
$\mc{A}$ whose local flow $\{\psi_t^F\}$
consists of multicontact maps. This means that if
$\delta\in\Delta_\mc{R}$, then
$$
\frac{d}{dt} (\psi_t^F)_*
(\overline{X}_\delta)\Big|_{t=0}=-\mathcal{L}_F
(\overline{X}_\delta)=[
\overline{X}_\delta,F],
$$
where $\mathcal{L}$ denotes the Lie derivative.
Hence a smooth vector field $F$ on $\mc{A}$ is a multicontact vector
field if and only if
\begin{equation} \label{multicontact}
[F,\overline{\Lie{g}}_\delta] \subseteq
\overline{\Lie{g}}_\delta\text{  for every $\delta\in\Delta_\mc{R}$}.
\end{equation}
We write a vector field on $\mc{A}$ as
\begin{equation}
F=\sum_{\gamma\in\mc{R}}\sum_{j=1}^{m_\gamma} f_{\gamma,j}
\overline{X}_{\gamma,j},
\label{standard}
\end{equation}
where $f_{\gamma,j}$ are smooth functions on $\mc{A}$.
Condition~(\ref{multicontact}) becomes
$$
[F, \overline{X}_{\delta,i}]=\sum_{k=1}^{m_\delta}
\lambda_{\delta,k}^i \overline{X}_{\delta,k},
\qquad
\delta\in\Delta_\mc{R},\;
i=1,\dots,m_\delta,
$$ where
$\{\lambda_{\delta,k}^i \}$ is a set of smooth functions.
By  Proposition \ref{ristretto} we have
\begin{align*}
[\overline{X}_{\alpha,i},\overline{X}_{\beta,j}]&=\overline{[X_{\alpha,i },X_{\beta,j}]}\\
&=\sum_{k=1}^{m_{\alpha+\beta}} c_{\alpha,\beta}^{ijk}
\overline{X}_{\alpha+\beta,k}
\hspace{.2in}\alpha,\beta\in\mc{R},\alpha+\beta\in\Sigma,
\end{align*}
where $c_{\alpha,\beta}^{ijk}$ are the structure constants with
respect to the chosen basis,
and $\overline{X}_{\alpha+\beta}=0$ if $\alpha+\beta\in\mc{C}$. We
can  write the multicontact conditions as the
system of equations
$$\sum_{\gamma\in\mc{R}} \sum_{j=1}^{m_\gamma}
\overline{X}_{\delta,i}(f_{\gamma,j}) \overline{X}_{\gamma,j}+
\sum_{\gamma\in\mc{R}}\sum_{j=1}^{m_\gamma}\left(
\sum_{l=1}^{m_{\gamma -\delta}}
c_{\delta,\gamma -\delta}^{ilj} f_{\gamma -\delta ,l}\right)
\overline{X}_{\gamma,j}=-\sum_{j=1}^{m_\delta}
\lambda_{\delta,j}^i
\overline{X}_{\delta,j},$$ as
$\delta$ varies in $\Delta_\mc{R}$ and $i=1,\dots,m_\delta$.
Equivalently, $F$ is a multicontact vector field on
$\mc{A}$ if and only if for all
$\gamma\in\mc{R}$ and some functions
$\{\lambda_{\delta,j}^i\}$ the following equations are satisfied on  $\mc{A}$:
\begin{equation} \label{sistema}
\begin{cases}
\overline{X}_{\delta,i} (f_{\delta,j})=-\lambda_{\delta,i}^j &\\
&\\
\overline{X}_{\delta,i} (f_{\gamma,j})=0 &\text{if }\gamma-\delta\not\in
\Sigma_+ \cup \{0\}\\
&\\
\displaystyle{\overline{X}_{\delta,i} (f_{\gamma,j})+
\sum_{l=1}^{m_{\gamma -\delta}}
c_{\delta,\gamma-\delta}^{ilj} f_{\gamma-\delta,l}=0} &\text{if
}\gamma -\delta\in\Sigma_+
\end{cases}
\end{equation}
for all the simple roots $\delta$ in $\Delta_\mc{R}$ and $1 \leq i,j
\leq  m_\delta$. We may
clearly forget the equation $\overline{X}_{\delta,i}
(f_{\delta,j})=-\lambda_{\delta,i}^j$ because
$\lambda_{\delta,i}^j$ is arbitrary.
\par
We fix some notations.  We write $MC(\group{N})$ and $MC(\group{S})$  for the
Lie algebra of multicontact vector fields on some open subset
of N and S respectevely. Let $\mc{C}$ be the complement in $\Sigma_+$ of
some Hessenberg type set. We say that a function $f$ on N
is $\mc{C}$-$independent$ if it does not depend on the coordinates  labeled by
$\mc{C}$. We record a simple consequence of \eqref{sistema} for later
reference.
\begin{lemma}\label{cindip} Let $F\in MC(\group{S})$ be as in  \eqref{standard}.
Then $f_{\gamma,j}$ are $\mc{C}$-independent for every  $\gamma\in\mc{R}$.
\end{lemma}
\begin{proof} All the  $f_{\gamma,j}$ appearing in \eqref{standard}
are functions
on S.
\end{proof}
Because of Lemma~\ref{cindip}, if $F\in MC(\group{S})$ is as in
\eqref{standard}, then
$\overline{X}_{\delta,i}f_{\gamma,j}=X_{\delta,i}f_{\gamma,j}$. Thus,  
from now
on we shall write $X_{\delta,i}$ in place of $\overline{X}_{\delta,i}$
whenever treating multicontact vector fields, if no ambiguity arises.
\par
  From \eqref{c} of Lemma \ref{left} it follows that
if $\mc{R}$ is a Hessenberg type
set of roots and $\gamma\in\mc{C}$, then a (basis) left invariant  vector field
${X}_{\gamma,k}$ on N does not depend on the partial derivative vector  fields
that are labeled by the positive roots in $\mc{C}$. This implies in
particular that the system of equations
\begin{equation} \label{total}
{X}_{\gamma,k} f=0 \hskip1cm \text{ for every
}\gamma\in\mc{C} \text{ and }
k=1,\dots,m_\gamma
\end{equation}
is equivalent to the $\mc{C}$-independence, namely to
\begin{equation} \label{indipendenzaparziale}
\frac{\partial}{\partial x_{\gamma,k}} f =0 \hskip1cm \text{
for every } \gamma\in\mc{C} \text{ and }
k=1,\dots,m_\gamma.
\end{equation}

\subsection{Dark zones}
We split (\ref{sistema}) into  suitable independent subsystems,
each defining multicontact vector fields on some Hessenberg
manifold of  lower dimension, and
we show that we can focus our attention to only one of them at a time.
\par
Call a positive root $\mu$ in $\mc{R}$ $maximal$ if
$\mu +\alpha\not\in\mc{R}$ for any
other root $\alpha\in\Sigma_+$. Since, by definition of $\mc{R}$,
$\mu +\alpha\notin
\mc{R}$ if
$\alpha\in\mc{C}$, it suffices to check maximality for all  $\alpha\in\mc{R}$.
Denote by
$\mc{R}_M$ the set of maximal  roots. For a fixed $\mu\in\mc{R}_M$, we  call
{\it shadow} of
$\mu$ the set
$$S_\mu = \{ \alpha\in\mc{R}:\alpha \preceq \mu \}.$$
Notice that $S_\mu \neq \emptyset$ because $\mu\in S_\mu$.
The union $\bigcup_{\mu\in\mc{R}_M}S_\mu$ covers $\mc{R}$.
Indeed, let $\alpha\in\mc{R}$. Either $\alpha\in\mc{R}_M$, or there
exists $\beta\in\mc{R}$ such that $\alpha + \beta \in\mc{R}$. Again,
$\alpha+\beta$ can be
maximal or not. If it is maximal, then $\alpha$ lies in its shadow.
If it is not
maximal, then we may continue this process and in a finite number of
steps we reach
a maximal root in whose shadow $\alpha$ lies.
\vskip6.8truecm
\setlength{\unitlength}{1cm}
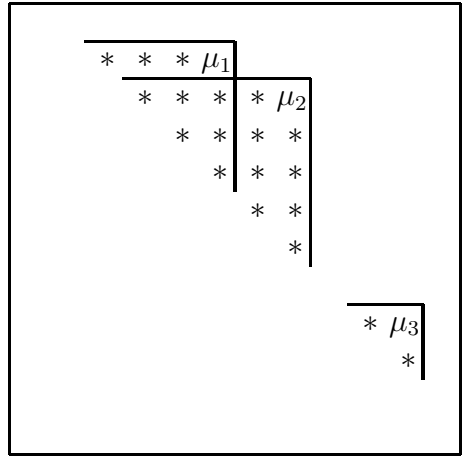
\begin{picture}(0,0)\thicklines
\put(3,0.5){\line(0,1){6}}
\put(3,0.5){\line(1,0){6}}
\put(3,6.5){\line(1,0){6}}
\put(9,0.5){\line(0,1){6}}
\put(4,6){\line(1,0){2}}
\put(6,4){\line(0,1){2}}
\put(4.5,5.5){\line(1,0){2.5}}%
\put(7,3){\line(0,1){2.5}}%
\put(7.5,2.5){\line(1,0){1}}%
\put(8.5,1.5){\line(0,1){1}}%
\put(4.2,5.65){$*$}
\put(4.7,5.65){$*$}
\put(5.2,5.65){$*$}
\put(5.55,5.65){$\mu_1$}
\put(4.7,5.15){$*$}
\put(5.2,5.15){$*$}
\put(5.7,5.15){$*$}
\put(6.2,5.15){$*$}
\put(6.55,5.15){$\mu_2$}
\put(5.2,4.65){$*$}
\put(5.7,4.65){$*$}
\put(6.2,4.65){$*$}
\put(6.7,4.65){$*$}
\put(5.7,4.15){$*$}
\put(6.2,4.15){$*$}
\put(6.7,4.15){$*$}
\put(6.2,3.65){$*$}
\put(6.7,3.65){$*$}
\put(6.7,3.15){$*$}
\put(7.7,2.15){$*$}
\put(8.05,2.15){$\mu_3$}
\put(8.2,1.65){$*$}
\end{picture}
\vskip0.1truecm
\par
In the picture above we see a representation of a
Hessenberg set $\mc{R}$ relative to $\Lie{sl}(n,\R)$ consisting
of three shadows. Here $\mc{R}_M=\{\mu_1,\mu_2,\mu_3\}$. For
simplicity, we label by $\mu_i$ the matrix entry
that corresponds to the root space $\Lie{g}_{\mu_i}$. The direct sum of
all the root spaces associated to roots in the shadow
$\mc{S}_{\mu_i}$ is a subspace whose coordinates belong to the
conical sector that has
$\Lie{g}_{\mu_i}$ as north-east corner. If one sees these corners
as idealized obstacles to a light beam coming from a far  north-east
point, then each cone is the shadow produced by it. This explains the
terminology.

We partition $\mc{R}$ into the disjoint union of $dark$ $zones$, a dark  zone
being a
connected component of $\mc{R}$ in a loose sense, that is, a maximal
union of shadows
$\mathcal{Z}=\cup_{i=1}^k
\mc{S}_{\mu_i}$ with the property that either $k=1$ or any
$\mc{S}_{\mu_i}$ intersects at least another
$\mc{S}_{\mu_j}$ in the same dark zone.
In the picture above, $\mc{R}$ is the union of two dark zones:
one is the union $\mc{S}_{\mu_1}\cup\mc{S}_{\mu_2}$ and the other
consists of the single shadow $\mc{S}_{\mu_3}$.
\par
By their very definition, dark zones are disjoint. This will allows us
to reduce the problem of solving \eqref{sistema} to the problem of  solving
several simpler systems, each naturally associated to a dark zone.  Suppose
that $\mathcal{Z}_1,\dots,\mathcal{Z}_p$ is a numbering of  the
dark zones of $\mathcal{R}$. Given
$F\in \mathfrak{X}(\group{S})$ as in~\eqref{standard},
   we write
$$
F=\sum_{i=1}^pF_i,
$$
where
$$
F_i=
\sum_{\gamma\in\mc{Z}_i}\sum_{j=1}^{m_\gamma}
f_{\gamma,j}\overline{X}_{\gamma,j}.
$$
Clearly, each $F_i$ is itself a vector
field in $\mathfrak{X}(\group{S})$. Since $F_i$ picks the
components of $F$ along the directions labeled by $\mathcal{Z}_i$, it
is natural to consider the sub-slice of $\group{S}$ that
corresponds to
it, as we now explain.
\vskip6.4truecm
\setlength{\unitlength}{1cm}
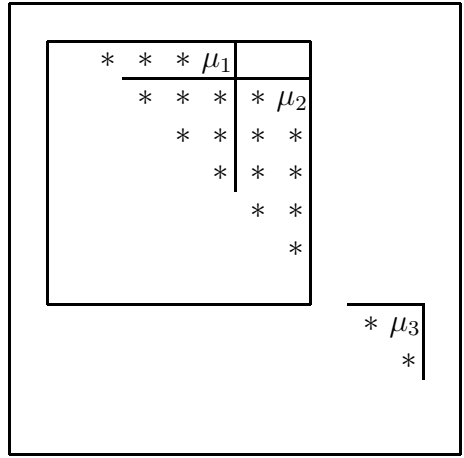
\begin{picture}(0,0)\thicklines
\put(3,0.5){\line(0,1){6}}
\put(3,0.5){\line(1,0){6}}
\put(3,6.5){\line(1,0){6}}
\put(9,0.5){\line(0,1){6}}
\put(3.5,6){\line(1,0){3.5}}
\put(6,4){\line(0,1){2}}
\put(4.5,5.5){\line(1,0){2.5}}%
\put(7,2.5){\line(0,1){3.5}}%
\put(3.5,2.5){\line(0,1){3.5}}%
\put(3.5,2.5){\line(1,0){3.5}}%
\put(7.5,2.5){\line(1,0){1}}%
\put(8.5,1.5){\line(0,1){1}}%
\put(4.2,5.65){$*$}
\put(4.7,5.65){$*$}
\put(5.2,5.65){$*$}
\put(5.55,5.65){$\mu_1$}
\put(4.7,5.15){$*$}
\put(5.2,5.15){$*$}
\put(5.7,5.15){$*$}
\put(6.2,5.15){$*$}
\put(6.55,5.15){$\mu_2$}
\put(5.2,4.65){$*$}
\put(5.7,4.65){$*$}
\put(6.2,4.65){$*$}
\put(6.7,4.65){$*$}
\put(5.7,4.15){$*$}
\put(6.2,4.15){$*$}
\put(6.7,4.15){$*$}
\put(6.2,3.65){$*$}
\put(6.7,3.65){$*$}
\put(6.7,3.15){$*$}
\put(7.7,2.15){$*$}
\put(8.05,2.15){$\mu_3$}
\put(8.2,1.65){$*$}
\end{picture}
\vskip0.1truecm
\par
Fix a dark zone
$\mc{Z}$.
The set of roots contained in $\mathcal{Z}$ generate the
positive set of an irreducible root system,
say
$\Sigma_+(\mathcal{Z})$, and the corresponding Lie algebra
$$
\Lie{n}(\mathcal{Z})=\bigoplus_{\beta\in\Sigma_+
(\mathcal{Z})}\Lie{g}_\beta
$$
is a nilpotent Iwasawa
algebra.  The roots in $\mathcal{Z}$ play, within
$\Sigma_+(\mathcal{Z})$, the r\^ole of a Hessenberg
    set of roots.
Also, $\Lie{n}(\mc{Z})$ is a subalgebra of $\Lie{n}$ and we may
consider the (connected, simply connected, nilpotent) Lie
subgroup $\group{N}(\mc{Z})$ of N whose Lie algebra is
$\Lie{n}(\mc{Z})$.
Thus, if $\mathcal{Z}$ is a dark zone we write
$$
\mc{S}_\mathcal{Z}=\{{n} \in \group{N} :
x_{\gamma,k}=0\;\text{ if }\;\gamma \not\in \mathcal{Z}\}.
$$
Coming back to the decomposition $\mc{R}=\mc{Z}_1 \cup\dots\cup
\mc{Z}_p$, we write for simplicity $\mc{S}_i$ in place of
$\mc{S}_{\mathcal{Z}_i}$. We want to prove the following
reduction result.
\begin{theorem}\label{reduction} If $F\in{MC}(\group{S})$, then
$F_i\in{MC}(\group{S}_i)$ for all $i=1,\dots, p$. Conversely, given
$G_i\in{MC}(\group{S}_i)$ with $i=1,\dots, p$, then
$\sum_iG_i\in{MC}(\group{S})$.
\end{theorem}
The proof requires some remarks, that we state in the
next lemmas.
\begin{lemma} \label{zone}
Let $\mathcal{Z}\subset\mathcal{R}$ be a dark zone and
let $\alpha\in \mathcal{Z}$. The $(\gamma,k)$ component of the vector  field
$X_{\alpha,j}$ is zero for every $\gamma\in\mc{R}\setminus
\mathcal{Z}$.
\end{lemma}
\begin{proof}
Suppose $\alpha\in \mathcal{Z}$, $\gamma\in\mc{R}\setminus
\mathcal{Z}$ and suppose the
$(\gamma,k)$-component of the vector field $X_{\alpha,j}$ is not zero.
By Lemma
\ref{left}, there exist roots
$\alpha_1,\dots,\alpha_q$ such that $\alpha +\alpha_1 +\dots
+\alpha_q =\gamma$ is a
chain, so that in particular
$\gamma - \alpha_q -\dots -\alpha_j$ is also a root for $j=1,\dots ,
q-1$. Now, since
$\gamma\in\mc{R}$, then
$\gamma\in S_\mu$ for some maximal root $\mu$. Therefore
$\alpha=\gamma - \alpha_q -\dots
-\alpha_1 \in S_\mu$.
This implies that
both $\alpha$ and $\gamma$ belong to the same shadow,
    and hence to
the same dark zone, that is a contradiction.
\end{proof}
\begin{lemma} \label{gerarchia}
The coefficients of a
multicontact vector field $F$
are determined by its
$\Lie{g}_\mu$ components, as $\mu$ varies in $\mc{R}_M$.
\end{lemma}
\begin{proof}
The proof of this statement is analogous to the proof of Proposition
3.3 of \cite{CDKR2}. We
give an outline for the reader's convenience.
Let $\beta$ be a positive root distinct from $\mu$. Let $\delta$ be a
simple root such
that $\beta +\delta$ is again a root. The set $\Delta (\beta)$ of all
such simple roots
is not empty by Lemma \ref{chain}. The equations (\ref{sistema})
yield
$$X_{\delta,i} (f_{\beta + \delta,j}) + \sum_{l=1}^{m_{\beta}}
c_{\delta,\beta}^{ilj} f_{\beta,l}=0,\hskip0.2cm\delta\in\Delta
(\beta),\hskip0.2cm,(i,j)\in \mc{I}_\delta (\beta),$$
where
$\mc{I}_\delta (\beta) =\{(i.j):1\leq i\leq d_\delta , 1\leq j\leq  d_{\beta +
\delta}\},\hskip0.2cm\delta\in\Delta (\beta)$.
We are thus led to consider the linear map given by the matrix
$A=(a_{I,l})=(c_{\delta,\beta}^{ilj})$ with row index $I=(i,j)$ and
column index $l$
varying in $\{1,\dots,d_\beta\}$.
    The proof consists then in showing
that $A$ has rank
$d_\beta$, because in this case the $f_{\beta,l}$ are uniquely given
by $X_{\delta,i}
(f_{\beta + \delta,j})$. This may be done using Lemma~\ref{trans} below,
that is proved in \cite{CDKR2} as a consequence of Kostant's
double transitivity theorem.
\end{proof}
\begin{lemma}\cite{CDKR2}\label{trans}
Let $\alpha,\beta\in\Sigma$ such that $\alpha + \beta$ is a root, then
$$\{[X,Y] : X\in\Lie{g}_\alpha, Y\in\Lie{g}_\beta \} =
\Lie{g}_{\alpha + \beta},$$
and $\{Z\in\Lie{g}_\beta : [\Lie{g}_\alpha , Z] =\{0\}\}=\{0\}$.
\end{lemma}
Lemma~\ref{gerarchia} suggests a hierarchic structure of the
equations (\ref{sistema}). In
particular, if $\gamma + \delta_1 + \dots + \delta_s = \alpha$ is a
chain, there exist vector fields $X_1\in\Lie{g}_{\delta_1},\dots ,
X_s\in\Lie{g}_{\delta_s}$ such that
the differential monomial $X_1\cdots X_s$ maps a ${\alpha}$-component  to a
$\gamma$-component of a vector field whose coefficients  solve~(\ref{sistema}).
In the proof the folowing result, we use again Lemma~\ref{trans}.
\begin{lemma}\label{ombra} Let $F\in MC(\group{S})$ be as in  \eqref{standard}.
Then  $Xf_{\gamma,j}=0$ for every $\gamma\in\mc{S}_\mu$, every
$j=1,\dots,m_\gamma$ and every
$X\in\Lie{g}_\alpha$ with
$\alpha\notin\mc{S}_\mu$.
\end{lemma}
\begin{proof}
If $\alpha\notin \mc{S}_\mu$, then it is either out of $\mc{R}$ or it  is in
some other shadow. If
$\alpha\in\mc{C}$, then $Xf_{\gamma,j}=0$ by Lemma~\ref{cindip}.
\par
Assume $\alpha\in\mc{R}$.
It is enough to prove the statement for $\gamma =\mu$. Indeed,
suppose the result true for all $f_{\mu,j}$'s. Then, by the equivalence  of
\eqref{total} and \eqref{indipendenzaparziale}, these functions are
$(\Sigma_+\setminus\mc{S}_\mu)$-independent, because
$\mc{S}_\mu$ is a Hessenberg type subset. If $\gamma +\delta_1+\dots
+\delta_p =\mu$ is a chain, then by
Lemma~\ref{gerarchia} there exist vector fields $X_1,\dots,X_p$ in
$\Lie{g}_{\delta_1},\dots,
\Lie{g}_{\delta_p}$ such that $X_1\cdots X_p f_{\mu,j} =
f_{\gamma,k}$. Each
$X_i$, $i=1,\dots,p$, has the form calculated in Lemma \ref{left}, that  is
$$X_i=\sum_{\alpha\in\Sigma_+} \sum_{j=1}^{m_\alpha} a^i_{\alpha,j}
\frac{\partial}{\partial x_{\alpha,j}},$$
where $a_{\alpha,j}$ is a nonzero polynomial only if there exists a
chain of roots going
from $\delta_i$ to $\alpha$, In this case $a^i_{\alpha,j}$ is a
polynomial in the variables
$\{x_{\beta,l}\}$ with $\beta \prec \alpha$. In particular this holds  for
$i=p$ and we show next that this forces
$X_p f_{\mu ,j}$ to be
$(\Sigma_+\setminus\mc{S}_\mu)$-independent. Indeed,
if $a^p_{\alpha,j}$ depends on some variable in
$(\Sigma_+\setminus\mc{S}_\mu)$,
then $\alpha\in(\Sigma_+\setminus\mc{S}_\mu)$ and therefore
   $\partial f_{\mu ,j}/\partial x_{\alpha,k}=0$ for all  $k=1,\dots,m_\alpha$.
Hence all coefficients $f_{\gamma,j}$ with ${\rm
ht}(\gamma)={\rm ht}(\mu)
-1$ are $(\Sigma_+\setminus\mc{S}_\mu)$-independent. By iterating the  same
argument, the conclusion holds for every possible height, thus for every
$\gamma$.
\par
It remains to be proved that the lemma is true for $f_{\mu,i}$.
If $\alpha$ is simple, then it is clear
by (\ref{sistema}) that $Xf_{\mu,j}=0$. Let now
$\alpha = \delta_1 + \dots + \delta_p$ be a non simple root in
$\mc{R}\setminus\mc{S}_\mu$. Then there exists
$\delta\in\{\delta_1,\dots,\delta_p\}$ such that
$\delta \notin \mc{S}_\mu$, for otherwise $\alpha \succ \mu$ and
$\mu$ would not be maximal.
By Lemma~\ref{trans} there exist vector fields $X_1,\dots,X_p$ in
$\Lie{g}_{\delta_1},\dots,
\Lie{g}_{\delta_p}$, respectively, such that $X=
[X_p,[\dots,[X_{2},X_{1}]]\dots]$. Then there exists a set
$\Lambda$ of permutations of $p$ elements such that
$$[X_p,[\dots,[X_{2},X_{1}]]\dots]f_{\mu,j}=(\sum_{\lambda\in\Lambda}
c_\lambda X_{\lambda(1)} \cdot \dots \cdot X_{\lambda(p)})f_{\mu,j},$$
for some costants $c_\lambda$. Let $h\in\{1,\dots,p\}$ be the largest
index such that
$\delta_{\lambda(h-1)}\notin \mc{S}_\mu$, so that clearly
$\delta_{\lambda(k)}$ is in $\mc{S}_\mu$ for
all $k \geq h$. We show that each differential monomial that appears
in the sum of the right hand side is zero on
$f_{\mu,j}$.
Consider $X_{\lambda(i)}\dots X_{\lambda(p)}$, with $i \geq h$.
Three possible cases arise.
\begin{itemize}
\item[(i)]
$\mu-\delta_{\lambda(p)}-\cdots
-\delta_{\lambda(i)}=0$, so that $\mu=\delta_{\lambda(p)}+\cdots
+\delta_{\lambda(i)}$. In this case $\alpha$ is the sum of
$\mu$ and some other simple roots. Hence $\alpha$ is a root in
$\mc{R}$ greater than
$\mu$, a contradiction.
\item[(ii)]
There exists $i\geq h$ such that $\mu-\delta_{\lambda(p)}-\cdots
-\delta_{\lambda(i+1)}$ is a positive root and  $\mu-\delta_{\lambda(p)}-\cdots
-\delta_{\lambda(i)}$ is not a root. In this case, from
Lemma~\ref{gerarchia} and the remark thereafter, the differential  monomial
$ X_{\lambda(i+1)}
\cdot
\dots
\cdot X_{\lambda(p)}$ maps $f_{\mu,j}$ into a component that belongs to  the
root space associated to
$\mu-\delta_{\lambda(p)}-\cdots
-\delta_{\lambda(i+1)}$, say $g$. Since
$\mu-\delta_{\lambda(p)}-\cdots
--\delta_{\lambda(i+1)}-\delta_{\lambda(i)}$ is not a root,
$X_{\lambda(i)}g=0$ by (\ref{sistema}).
\item[(iii)]
$\mu-\delta_{\lambda(p)}-\cdots
-\delta_{\lambda(i)}$ is a root for all $i\geq h$. Again the
differential monomial
$ X_{\lambda(h)}\dots X_{\lambda(p)}$ maps $f_{\mu,j}$ into a component  along
the root space labeled by $\mu-\delta_{\lambda(p)}-\cdots
-\delta_{\lambda(h)}$. But
$\mu-\delta_{\lambda(p)}-\cdots  -\delta_{\lambda(h)}-\delta_{\lambda(h-1)}$
is not a root, for otherwise
$\delta_{\lambda(h-1)}$ would lie in $\mc{S}_\mu$. Therefore we can
conclude as in the
previous case. Thus $ X_{\lambda(h-1)}
\dots
X_{\lambda(p)}$ maps the function $f_{\mu,j}$ to zero.
\end{itemize}
\end{proof}
{\it Proof of Theorem~\ref{reduction}}.
\par
``$\Rightarrow$''. Lemma \ref{ombra} applies in particular to
each dark zone, in the sense that a coefficient $f_{\gamma,k}$ of a
multicontact vector
field on S is annihilated by those left invariant vector fields
corresponding to the roots
that do not belong to the dark zone where
$\gamma$ lies.
    Since each dark zone plays the r\^{o}le of a Hessenberg set of  roots, its
complement defines an ideal in
$\Lie{n}$, namely
$$
\Lie{n}_{\mc{Z}^c} = \bigoplus_{\alpha\in\Sigma_+ \setminus
\mc{Z}} \Lie{g}_\alpha,
$$
where $\mc{Z}^c = \Sigma_+\setminus \mc{Z}$.
The corresponding nilpotent Lie group admits the set
$\{X_{\alpha,j} : \alpha\in \Sigma_+ \setminus
\mc{Z}\}$ as a basis for its tangent space at each point.
   From \eqref{c} in Lemma~\ref{left}, all these vector fields
depend on the coordinate vector fields labeled by the positive
roots in $\Sigma_+ \setminus\mc{Z}$. Recall in particular that from
\eqref{total} and \eqref{indipendenzaparziale}
$$
{X}_{\gamma,k} f=0
\text{ for all }\gamma\notin\mathcal{Z}
\iff
\frac{\partial}{\partial x_{\gamma,k}} f =0
\text{ for all }\gamma\notin\mathcal{Z} .
$$
This fact, toghether with Lemma \ref{ombra}, tells us
that the coefficients of the vector
field $F_i$ are functions on $\group{S}_i$, that is, they are
$(\mc{R}\setminus\mc{Z})$-independent.
\par
Moreover, by
Lemma \ref{zone}, the projections
$\overline{X}_\delta$ onto the tangent space at each point of S
are in fact projections on the tangent space of $\group{S}_i$.
Therefore
$F_i\in\mathfrak{X}(\group{S}_i)$. Hence $F_i$ is in
$MC(\group{S}_i)$ if and only if
\begin{equation} \label{blocco}
\begin{cases}
\overline{X}_{\delta,i} (f_{\gamma,j})=0
&
\gamma-\delta\not\in \Sigma_+ \cup
\{0\}\\
&\\
\displaystyle{ \overline{X}_{\delta,i} (f_{\gamma,j})+  \sum_{l=1}^{m_{\gamma
-\delta}} c_{\delta,\gamma-\delta}^{ilj}
f_{\gamma-\delta,l}=0}
&\gamma-\delta\in\Sigma_+,
\end{cases}
\end{equation}
with $\delta\in\Delta\cap\mc{Z}_i$ and $\gamma\in\mc{Z}_i$.
We conclude by observing that these equations are satisfied by  assumption.
\par
``$\Leftarrow$''. Each vector field $G_i$ can be naturally viewed
as a vector field on S. Furthermore, since each $G_i$ satisfies
the system of equations \eqref{blocco}, then the vector field
$\sum_i G_i$ satisfies the system \eqref{sistema}. Thus, it
defines a multicontact vector field on S. This concludes the proof of  the
theorem.
\qed
\medskip
\par
Theorem \ref{reduction} allows us to study each dark zone separately.
  From now on we thus assume that the Hessenberg set contains all  simple roots.
\subsection{A set of solutions}
A further step in investigating the system of differential equations
(\ref{sistema})
allows us to find a set of solutions.
\par
In \cite{CDKR2}, the authors determine the
multicontact vector  fields on the Iwasawa group $\group{N}$, by
solving a system of
differential equations similar to (\ref{sistema}). In particular, if
$V=\sum_{\gamma\in\Sigma_+}\sum_{j=1}^{m_\gamma} v_{\gamma,j}
X_{\gamma,j}$ is a vector field on $\group{N}$, then $V$ is of
multicontact type if it satisfies the following system of
equations
\begin{equation} \label{iwasawa}
\begin{cases}
X_{\delta,i}
(v_{\gamma,j})=0 &\text{ if } \gamma-\delta\not\in \Sigma_+ \cup
\{0\}\\
&\\
\displaystyle{X_{\delta,i} (v_{\gamma,j})+ \sum_{l=1}^{m_{\gamma -\delta}}
c_{\delta,\gamma-\delta}^{ilj} v_{\gamma-\delta,l}=0} &\text{ if }
\gamma -\delta\in\Sigma_+,
\end{cases}
\end{equation}
where $\gamma$
varies in $\Sigma_+$, $\delta$ in $\Delta$, and the $v_{\gamma,j}$ are  smooth
functions on
$\group{N}$. Write
$$\overline{V}=\sum_{\gamma\in\Sigma_+} \sum_{j=1}^{m_\gamma}  v_{\gamma,j}
\overline{X}_{\gamma,j}.
$$
If
$V$ solves  (\ref{iwasawa}), then the projection $\overline{V}$  satisfies
(\ref{sistema}). Moreover, if the coefficients $v_{\gamma,j}$ are
$\mc{C}$-independent for every $\gamma\in\mc{R}$, then the vector
field $\overline{V}$ is
tangent at each point to $\group{S}$. Summarizing, in this case
$\overline{V}$ is a
multicontact vector field on $\group{S}$. In \cite{CDKR2} it is proved  that
the  multicontact
vector fields on
$\group{N}$ are all of the form $\tau(E)$ for some $E\in\Lie{g}$, where
\begin{equation} \label{tau}
\tau(E) h({n}) = \frac{d}{dt} h(\exp(-tE)
{n}) \Big|_{t=0}.
\end{equation}
We ask ourselves for which $E\in\Lie{g}$ the coefficients of
$\overline{\tau(E)}$ are
$\mc{C}$-independent. Denote
by $\Lie{q}$ the parabolic subalgebra of $\Lie{g}$ defined as the
normalizer in $\Lie{g}$ of
$\Lie{n}_\mc{C}$
$$\Lie{q}:={N}_{{\Lie{g}}}
\Lie{n}_{\mc{C}}=\{X\in\Lie{g}:[X,Y]\in\Lie{n}_\mc{C}, \forall
Y\in\Lie{n}_\mc{C}\}.$$
Clearly $\Lie{q}\supset \Lie{m}\oplus\Lie{a}\oplus\Lie{n}$, so that
$\Lie{q}$ is a
parabolic subalgebra of $\Lie{g}$.
\begin{theorem} \label{andata}
Let $\mc{R}\subseteq \Sigma_+$ a Hessenberg type set, $\mc{C}$ the
complement of $\mc{R}$, and
$\Lie{q}= {N}_{\Lie{g}}\Lie{n}_{\mc{C}}$.
For
every $E \in \Lie{q}$, $\overline{\tau(E)}$ is a multicontact vector  field on
$\group{S}$. In particular, the map
\begin{equation}
\label{aumorfismo}
\nu: \Lie{q} \longrightarrow
{\mathfrak X}(\group{S})
\end{equation}
defined by
$\nu(E)=\overline{\tau(E)}$ is a Lie algebra homomorphism.
If $\Delta\subset\mc{R}$, then the kernel of $\nu$ is
$\Lie{n}_\mc{C}$. Thus $\nu(\Lie{q})$ is isomorphic to
$\Lie{q}/\Lie{n}_\mc{C}$.
\end{theorem}
\begin{proof}
In order to prove the first claim we show that the coefficients of
     $\overline{\tau (E)}$ are $\mc{C}$-independent for every $E \in  \Lie{q}$.
Let $E^\prime\in\Lie{n}_\mc{C}$. Then
$$
[\tau(E), \tau(E^\prime)]=[\sum_{\alpha\in\mc{R}}
\sum_{i=1}^{m_\alpha} f_{\alpha,i}
X_{\alpha,i}  + \sum_{\beta\in\mc{C}} \sum_{j=1}^{m_\beta}
f_{\beta,j} X_{\beta,j} ,
\sum_{\gamma\in\mc{C}} \sum_{k=1}^{m_\gamma} g_{\gamma,k} X_{\gamma,k}]
$$
must lie in $\tau(\Lie{n}_\mc{C})$. By direct calculation, this
happens if and only if $X_{\gamma,k} (f_{\alpha,i})=0$, or
equivalently if and only if
$$
\frac{\partial}{\partial x_{\gamma,k}} (f_{\alpha,i}) =0
$$
for every $\alpha\in\mc{R}$ and $\gamma\in\mc{C}$.
\par
The map $\nu$
is a homomorphism because $\tau$ and the projection
operator are such. Hence $\nu(\Lie{q})$ is a  Lie algebra of
multicontact vector fields
on
$\group{S}$.
\par
We now investigate the kernel of $\nu$ in the case
$\Delta\subset\mc{R}$. Since
$\tau(E)=\sum_{\gamma\in\mc{C}} \sum_{k=1}^{m_\gamma} g_{\gamma,k}
X_{\gamma,k}$ for every
$E\in\Lie{n}_\mc{C}$, the inclusion
$\Lie{n}_\mc{C} \subseteq {\rm ker}\nu$ follows.
We prove the opposite inclusion by treating
separetely each component of $E\in{\rm ker}\nu$, written
according to
the following vector space direct sum:
$$\Lie{q}=\Lie{m}\oplus\Lie{a}\oplus\Lie{n}\oplus(\overline{\Lie{n}}\cap \Lie{q}).$$
  From now until the end of this proof we write ${n}=\exp
(W)=\exp(\sum_{\alpha\in\Sigma_+}W_\alpha)$, where  $W_\alpha\in\Lie{g}_\alpha$.

\vskip0.5cm
\noindent
If $E\in\Lie{n}\cap{\rm ker}\nu$, then
$\overline{\tau(E)}=0$. Write $E=\sum_{\gamma\in\Sigma_+}  \sum_{k=1}^{m_\gamma}
a_{\gamma,k} E_{\gamma,k}$ and compute
\begin{align*}
\tau(E) f &= \frac{d}{dt} f(\exp (-tE) n)\Big|_{t=0} \\
&=\frac{d}{dt} f(\exp(-tE + W -\frac{t}{2} [E,W] + \dots ))\Big|_{t=0}.
\end{align*}
If $E$ were not in $\Lie{n}_\mc{C}$, there would exist $\beta \in  \mc{R}$ and
$j=1,\dots,m_\beta$ such that $a_{\beta,j} \neq 0$. If $f:n\mapsto
x_{\beta,j}$ then we
have that
$\tau(E) f$ is a polynomial in $\{x_{\alpha,i}\}_{\alpha\in\Sigma_+}$
whose term of
degree zero is
$a_{\beta,j}$. On the other hand
$$
\tau(E)x_{\beta,j}=0\;\forall\;\beta\in\mc{R},
$$
because its
decomposition on the
basis of left invariant vector fields involves only components
corresponding to the roots in
$\mc{C}$. This is a
contradiction.
\vskip0.5cm
\noindent
Let $E\in\Lie{a}\cap{\rm ker}\nu$. Recalling that we view N as a
dense subset of G/P
and that ${\rm exp}(tE)\in \group{P}$, we have
\begin{align*}
\tau(E)f ({n})&=\frac{d}{dt} f(\exp (-tE) {n})\Big|_{t=0}\\
&= \frac{d}{dt}f(\exp (-tE) {n}\exp (tE))\Big|_{t=0}\\
&=\frac{d}{dt} f(\exp (\sum_{\alpha\in\Sigma_+}
e^{-t\alpha(E)}W_\alpha))\Big|_{t=0}.
\end{align*}
Choose now $f:n\mapsto x_{\gamma,j}$, so that
     $$
\tau(E)f({n})=
\frac{d}{dt}(e^{-t\gamma(E)}x_{\gamma,j})\Big|_{t=0}f(n)=- \gamma(E)x_{\gamma,j}.
$$
This is zero for every $\gamma\in\mc{R}$ because $E$ is in the kernel
of $\nu$, so that
$\gamma(E)=0$ for every $\gamma\in\mc{R}$. Since
$\mc{R}\supset\Delta$ and $\Delta$ is a
basis of $\Lie{a}^*$, the dual space of $\Lie{a}$, it follows that  $E=0$.
\par
Let $E\in\Lie{m}\cap {\rm ker}\nu$.
Since $\Lie{m}$ normalizes every root
space, if $f:n\mapsto x_{\gamma,j}$, then
\begin{align*}
\tau(E)f({n})&=\frac{d}{dt}f(\exp (e^{-{\rm ad}tE}W) )\Big|_{t=0}\\
&=\frac{d}{dt}f(\exp (\sum_{\alpha\in\Sigma_+}\sum_{n=1}^\infty
(-1)^n t^n \frac{({\rm
ad} E)^n}{n!} W_\alpha))\Big|_{t=0}\\ &=((-{\rm ad}E)W_\gamma)_j.
\end{align*}
Whenever $\gamma\in\mc{R}$ we have $((-{\rm ad}E)W_\gamma)_j=0$ for
every $j$. Thus
$({\rm ad} E)\Lie{g}_\gamma=0$ for every
$\gamma\in\mc{R}$. In particular $({\rm ad} E)\Lie{g}_\delta=0$ for
every simple root
$\delta$, and Jacobi identity implies $({\rm ad} E)\Lie{n}=0$.
Since $\theta E=E$, it follows that $({\rm ad}
E)\Lie{g}_{-\delta}=({\rm ad} \theta E)\Lie{g}_{\delta}=({\rm
ad} E)\Lie{g}_\delta=0$. Hence $({\rm ad} E)\Lie{g}=0$. Thus
$E\in Z(\Lie{g})=\{0\}$.
\par
Let now $E\in\Lie{g}_\beta\cap\Lie{q}\cap{\rm ker}\nu$ for some  negative root
$\beta$, so that $\overline{\tau(E)}=0$.

For every $E^\prime\in\Lie{n}$ we have
$$
[\tau(E), \tau(E^\prime)]=[\sum_{\alpha\in\mc{C}}
\sum_{i=1}^{m_\alpha} f_{\alpha,i}
X_{\alpha,i},
\sum_{\beta\in\mc{R}} \sum_{j=1}^{m_\beta} g_{\beta,j} X_{\beta,j} +
\sum_{\gamma\in\mc{C}} \sum_{k=1}^{m_\gamma} g_{\gamma,k} X_{\gamma,k}]
$$
All terms of the bracket above lie on $\Lie{n}_\mc{C}$, except for  summands of the form
$$
f_{\alpha,i} X_{\alpha,i} (g_{\beta,j}) X_{\beta,j},
$$
but $X_{\alpha,i} (g_{\beta,j})=0$, for every $\alpha\in\mc{C}$ and
$\beta\in\mc{R}$,
because the coefficients $g_{\beta,j}$ are $\mc{C}$-independent.
It follows in particular that
$$\overline{[\tau(E),\tau(E^\prime)]}=0,$$
thus $[E,E^\prime]\in{\rm ker}\nu$ for every $E^\prime\in\Lie{n}$.
Therefore one can chose $E^\prime$ such that
$[E,E^\prime]\in\Lie{m}\oplus\Lie{a}$. But
this is a contradiction, because no elements of
$\Lie{m}\oplus\Lie{a}$ lie in the kernel
of
$\nu$.
\end{proof}
Notice that in the last step of the above proof we did not use that
the negative root is
in the normalizer, so that there are no multicontact vector fields on
N coming from a
negative root that become zero once projected to a slice representing
a Hessenberg
manifold.


\section{Iwasawa sub-models}
The converse of Theorem~\ref{andata} is true
under the hypothesis $(I)$ of the Theorem~\ref{al} below.
We remind the reader that by Theorem~\ref{reduction} we are assuming
that $\mc{R}$ consists of a single dark zone and that it contains all
the simple restricted roots.
\begin{lemma} \label{niliwa}
If the vector space
$$
\Lie{n}^\mu =\bigoplus_{\alpha\in \mc{S}_\mu}\Lie{g}_\alpha
$$
is a subalgebra of $\Lie{n}$, then in particular it is an Iwasawa
nilpotent Lie algebra.
\end{lemma}
\begin{proof}
The algebra $\Lie{n}^\mu$ coincides with the nilpotent algebra
generated by the root spaces corresponding to the simple roots in
$\mc{S}_\mu$. Hence it is an Iwasawa Lie algebra because it is
the canonical nilpotent algebra associated to a connected Dynkin
diagram, together with admissible multiplicity data.
\end{proof}
The following theorem holds.
\begin{theorem} \label{al}
Let $\Lie{g}$ be a simple Lie algebra of real rank strictly greater
than two and
$\mc{R}\subset\Sigma_+$ a subset of Hessenberg type satisfying
\begin{itemize}
\item[(I)]
each shadow in the
Hessenberg set defines a subalgebra of $\Lie{n}$,
\item[(II)] each shadow contains at least two simple roots.
\end{itemize}
Then the Lie algebra of multicontact vector fields on ${\rm
Hess}_\mc{R}  (H)$ is isomorphic to
$\Lie{q}/\Lie{n}_\mc{C}$, for every regular element
$H\in\Lie{a}$ and where $\Lie{q}=N_\Lie{g} \Lie{n}_\mc{C}$.
\end{theorem}
Hypothesis (II) just avoids the rank one cases. More precisely, if
(II) is not true, then $\mc{R}$ defines a rank one Iwasawa subalgebra.
In this case, however, the finite dimensionality of the Lie algebra
$MC$(S) is no longer guaranteed\footnote{Personal comunication  by
the authors of~\cite{CDKR2}, who intend to clarify this matter in full
detail in a forthcoming paper.}.
\par
We have seen that
if $\mc{R}$ is an arbitrary Hessenberg type set, then
all elements in $\Lie{q}/\Lie{n}_\mc{C}$ define
multicontact vector fields.
In order to show the converse, we look again at the system
of  differential equation
(\ref{sistema}).
If $F\in MC(\group{S})$, then its coefficients solve
(\ref{sistema}), and they solve all the subsystems that we can
extract from (\ref{sistema}). Therefore we obtain necessary
conditions by looking at some special subsystems of
(\ref{sistema}). In particular we consider a subsystem for
each shadow, namely:
\begin{equation} \label{shadow}
\begin{cases}
X_{\delta,i}
(f_{\gamma,j})=0 &\delta\in S_\mu, \text{ if } \gamma-\delta\not\in
\Sigma_+
\cup
\{0\}\\
&\\
\displaystyle{ X_{\delta,i} (f_{\gamma,j})+ \sum_{l=1}^{m_{\gamma -\delta}}
c_{\delta,\gamma-\delta}^{ilj} f_{\gamma-\delta,l}=0} & \text{ if }
\gamma -\delta\in\Sigma_+ \\
&\\
X_{\delta,i}(f_{\gamma,j})=0 &\delta \notin S_\mu,
\end{cases}
\end{equation}
for every root $\gamma$ in $\mathcal{S}_\mu$. Here we stress
that the functions
$f_{\gamma,j}$ are defined on some open
subset of the slice S. Since each shadow defines an Iwasawa
Lie algebra, the
system above looks like the system of differential equations that
defines the multicontact
vector fields on nilpotent Iwasawa Lie groups. We want to interpret
\eqref{shadow} exactly
in this way. Indeed, Lemma \ref{ombra} tells us that the functions
$f_{\gamma,j}$, as
$\gamma$ varies in $\mc{S}_\mu$, are
$(\Sigma_+\setminus\mc{S}_\mu)$-independent.  Hence
$$
X_{\delta,i} (f_{\gamma,j})=X_{\delta,i}^\mu
(f_{\gamma,j}),\hskip0.1cm \text{for every
}\gamma,\delta\in\mc{S}_\mu,
$$
where $X_{\delta,i}^\mu$ is the vector field that is obtained
from $X_{\delta,i}$ by setting to zero all the components that are
labeled by roots that are not in $S_\mu$.
We then consider, in place of \eqref{shadow}, the equivalent system
    \begin{equation} \label{sh1}
\begin{cases}
X_{\delta,i}^\mu
(f_{\gamma,j})=0 &\text{ if } \gamma-\delta\not\in
\Sigma_+
\cup
\{0\}\\
&\\
\displaystyle{X_{\delta,i}^\mu (f_{\gamma,j})+ \sum_{l=1}^{m_{\gamma -\delta}}
c_{\delta,\gamma-\delta}^{ilj} f_{\gamma-\delta,l}=0} & \text{ if }
\gamma -\delta\in\Sigma_+,
\end{cases}
\end{equation}
where $\gamma,\delta \in S_\mu$.
Define $\Lie{n}^\mu$ as in Lemma \ref{niliwa}. Using hypothesis (I) of
Theorem \ref{al}, Lemma \ref{niliwa} implies that the Lie algebra
$\Lie{n}^\mu$ is an Iwasawa nilpotent Lie algebra. The system of
differential equations above coincides with the  multicontact conditions
for a vector field on
$\group{N}^\mu=\exp
\Lie{n}^\mu$, because the vector fields
$ X_{\delta,i}^\mu$ are exactly the left--invariant vector fields on
$\group{N}^\mu$.
This latter assertion is rather obvious, and it is a consequence of a
direct calculation similar to the one involved in Lemma \ref{left}.
Summarizing, we have the following result.
\begin{prop} \label{compatibilita}
Let $F\in \mathfrak{X}(\group{S})$. Then
$F\in MC(\group{S})$ if and only
if its projection
$$
F^\mu=\sum_{\alpha\in
\mc{S}_\mu}\sum_{i=1}^{m_\alpha}f_{\alpha,i}\overline{X}_{\alpha,i},
$$
is a multicontact vector field on $\group{N}^\mu$ for
every maximal root $\mu$.
\end{prop}
\begin{proof}
``$\Rightarrow$''.
By Lemma \ref{ombra}, any multicontact vector field on S can be  naturally
viewed  as a vector field on
$\group{N}^\mu$ for every maximal root $\mu$.
If the coefficients of $F$ solve the system
of differential equations \eqref{sistema}, then in particular they
solve all subsystems \eqref{sh1}, that is any projected vector field
$F^\mu$ is in $MC(\group{N}^\mu)$.
\par
``$\Leftarrow$''. If $F$ has the property that each
    $F^\mu$ solves \eqref{sh1}, then $F$ solves all the equations in
\eqref{sistema}, so that
it is in $MC(\group{S})$.
\end{proof}

The multicontact vector fields on a nilpotent Iwasawa Lie group are
studied in \cite{CDKR2}, where it is showed in particular that the Lie
algebra of such vector fields on $\group{N}^\mu$
is isomorphic to
$\Lie{g}^\mu = \Lie{n}^\mu + \theta \Lie{n}^\mu + \Lie{m}^\mu +
\Lie{a}^\mu$, where
$$\Lie{m}^\mu=\Lie{m}\cap[\Lie{n}^\mu,\theta\Lie{n}^\mu],$$
and
$$\Lie{a}^\mu =\Lie{a} \cap [\Lie{n}^\mu,\theta\Lie{n}^\mu].$$
In fact each element of this algebra defines a vector field on
$\group{N}^\mu$ whose
coefficients solve (\ref{shadow}). These vector fields are
realized as $\tau_\mu (E)$, where
$$\tau_\mu (E)f(n)=\frac{d}{dt} f(\exp(-tE)n)\Big|_{t=0},$$
with $E\in\Lie{g}^\mu$, $n\in\group{N}^\mu$ and some
function $f$ on
$\group{N}^\mu$.
Since $\Lie{g}^\mu$ is a subalgebra of $\Lie{g}$, it is possible to see
the multicontact vector fields on $\group{N}^\mu$ as the projections of
suitable vector fields on N, in the sense that is explained in the
following proposition.
\begin{prop}\label{rappre}
The set of vector fields
$$
\{\tau(E)^\mu, E\in\Lie{g}^\mu\}
$$
generates the Lie algebra $MC(\group{N}^\mu)$, where
$$
\tau(E)^\mu =\sum_{\gamma\in\mc{S}_\mu}\sum_{j=1}^{m_\gamma}  f_{\gamma,j}
\overline{X}_{\gamma,j},
$$
whenever $\tau(E)=\sum_{\gamma\in\Sigma_+}\sum_{j=1}^{m_\gamma}  f_{\gamma,j}
{X}_{\gamma,j}$. In particular, if $E\in\Lie{q}$, it follows that
$\tau(E)^\mu\neq 0$ if and only if $E\in\Lie{g}^\mu\setminus\{0\}$.
\end{prop}
\begin{proof}
Let $E\in\Lie{g}^\mu$.
We show that $E\in\Lie{b}$, the normalizer
in $\Lie{g}$ of the nilpotent ideal consisting of all the
root spaces labeled by $\mc{S}_{\mu}^c=\Sigma_+\setminus
\mc{S}_\mu$, namely
$\Lie{b}=N_\Lie{g}\Lie{n}_{\mc{S}_{\mu}^c}$. Since
$\Lie{g}^\mu= \Lie{m}^\mu +
\Lie{a}^\mu +\Lie{n}^\mu +
\theta
\Lie{n}^\mu $ and $\Lie{m}^\mu +
\Lie{a}^\mu +\Lie{n}^\mu \subseteq \Lie{b}$, we can suppose that
$E\in
\theta\Lie{n}^\mu$. Write $E=\sum E_\beta$. If
$E\notin\Lie{b}$, then there exists $\beta$ such that
$E_\beta\notin\Lie{b}$. In this case, since $\Lie{b}$
normalizes, there
would exist
$\alpha\in{\mc{S}_\mu^c}$ such
that
$
\alpha +\beta \notin \mc{S}_\mu$--a
contradiction,
because the sum of two roots in $\mc{S}_\mu$ is in $\mc{S}_\mu$, as
follows from hypothesis
(I).
Theorem \ref{andata} applied to the Hessenberg set
$\mc{S}_\mu$ implies that $\tau(E)^\mu\in MC(\group{N}^\mu)$.
\par
We now show that the vector fields $\tau(E)^\mu$ are all different
from zero whenever $E\in\Lie{g}^\mu \setminus \{0\}$.
Suppose that there exists $E\in\Lie{g}^\mu$ such that $\tau(E)^\mu =0$.
Write $E=H+K+\sum E_\alpha$, with $H\in\Lie{a}^\mu$ and
$K\in\Lie{m}^\mu$. Since
$Y\mapsto Y^\mu$ preserves (homomorphic images of) root spaces, the
hypothesis
$\tau(E)^\mu =0$ is equivalent to assuming $\tau(H)^\mu=\tau(K)^\mu=0$
and
$\tau(E_\alpha)^\mu =0$ for every $\alpha$.
We show first that $H=0$. Indeed, writing
$n=\exp (\sum_{\alpha\in\Sigma_+} W_\alpha)$, we have
\begin{align*}
\tau(H)f(n)&=\frac{d}{dt} f(\exp (-tH)n)\Big|_{t=0}\\
&=\frac{d}{dt}
f(\exp(\sum_{\alpha\in\Sigma_+}e^{-t\alpha(H)}W_\alpha)\Big|_{t=0}.
\end{align*}
Hence, if $f:n\mapsto x_{\gamma ,j}$, then
$$\tau(H)f(n)=\frac{d}{dt}(e^{- t\gamma(H)}x_{\gamma,j})f(n)\Big|_{t=0}=-\gamma(H)x_{\gamma,j}.$$
This is zero for every $\gamma$ in $S_\mu$, and in particular
$\delta(H)=0$ for every simple roots $\delta$ in $S_\mu$. By duality
this implies that $H=0$.
\par
Suppose now that $\tau(K)^\mu=0$.
Since $\Lie{m}^\mu$ normalizes every root
space, if $f:n\mapsto x_{\gamma,j}$, then
\begin{align*}
\tau(K)f({n})&=\frac{d}{dt}f(\exp (e^{-{\rm ad}tK}W) )\Big|_{t=0}\\
&=\frac{d}{dt}f(\exp (\sum_{\alpha\in\Sigma_+}\sum_{n=1}^\infty
(-1)^n t^n \frac{({\rm
ad} K)^n}{n!} W_\alpha))\Big|_{t=0}\\ &=((-{\rm ad}K)W_\gamma)_j.
\end{align*}
Whenever $\gamma\in\mc{S}_\mu$ we have $((-{\rm ad}K)W_\gamma)_j=0$ for
every $j$. Thus
$({\rm ad} K)\Lie{g}_\gamma=0$ for every
$\gamma\in\mc{S}_\mu$. In particular $({\rm ad} K)\Lie{g}_\delta=0$ for
every simple root
$\delta\in\mc{S}_\mu$, and Jacobi identity implies $({\rm ad}
K)\Lie{n}^\mu=0$. Since $\theta K=K$, it follows that $({\rm ad}
K)\Lie{g}_{-\delta}=({\rm ad} \theta K)\Lie{g}_{\delta}=({\rm
ad} K)\Lie{g}_\delta=0$. Hence $({\rm ad} K)\Lie{g}^\mu=0$. Thus
$K\in Z(\Lie{g}^\mu)=\{0\}$.
\par
Next, suppose that $\tau(E_\alpha)^\mu =0$. Then
$$
0=[\tau(E_\alpha)^\mu ,\tau(\theta E_\alpha)^\mu ]
=\tau([E_\alpha,\theta E_\alpha])^\mu.
$$
Recall that if $E_\alpha\not=0$, then
\begin{equation}
[E_\alpha,\theta E_\alpha]=B(E_\alpha,\theta
E_\alpha)H_\alpha
\label{nondegenerate}\end{equation}
(see e.g. Prop 6.52 in \cite{KN02}),
where $H_\alpha\in\Lie{a}$ represents $\alpha$ via the Killing form
and $B(E_\alpha,\theta E_\alpha)<0$.
But $0=\tau([E_\alpha,\theta E_\alpha])^\mu=\tau(H_\alpha)$ implies
$H_\alpha=0$ by the previous case,
contradicting~\eqref{nondegenerate}.
\par
Therefore the set
$$
\{\tau(E)^\mu, E\in\Lie{g}^\mu\}
$$
generates the Lie algebra of multicontact vector fields on
$\group{N}^\mu$, as required.
\par
Finally, let $E\in\Lie{q}$. We proved above that if $E\in\Lie{g}^\mu$  then
$\tau(E)^\mu\neq 0$. On the other hand, $E\in\Lie{q}$ implies that
$\overline{\tau(E)}\in MC(\group{S})$, so that in particular
$\tau(E)^\mu\in MC(\group{N}^\mu)$. If $E\notin\Lie{g}^\mu$, then the
latter assertion is possible only if $\tau(E)^\mu =0$.
\end{proof}
We prove a result that comes as a direct consequence of
Proposition~\ref{rappre}. It essentially says that the intersection of  two
or more shadows still defines an Iwasawa nilpotent subalgebra, and that  we
can represent the multicontact vector fields on the corresponding  subgroup
again by projecting some multicontact vector field of the form  $\tau(E)$.
This result will be used later, and it can be viewed as a  generalization of
the extra-cross condition~\eqref{comp} of the Example~\ref{ex}.
\begin{cor}\label{inters}
Let $\mc{I}=\bigcap_{\mu\in \mc{E}} \mc{S}_\mu$ with $\mc{E}$ a subset  of
maximal roots in $\mc{R}$. Then:
\begin{itemize}
\item[(i)] the nilpotent Lie algebra
$\Lie{n}^\mc{I}=\bigoplus_{\alpha\in \mc{I}} \Lie{g}_\alpha$
is an Iwasawa Lie algebra.
\item[(ii)] Let
$\Lie{g}^\mc{I}$ denote the Lie subalgebra of $\Lie{g}$ generated by  $\Lie{n}^\mc{I}$ and
 $\theta\Lie{n}^\mc{I}$, and let $\group{N}^\mc{I}=\exp \Lie{n}^\mc{I}$.
The vector fields of the type
$$
\tau(E)^\mc{I}=\sum_{\alpha\in \mc{I}} \sum_{j=1}^{m_\gamma}
f_{\gamma,j}\overline{X}_{\gamma,j},
$$
with $E\in\Lie{g}^\mc{I}$, are in $MC(\group{N}^\mc{I})$.
\item[(iii)] If $E\in\Lie{q}$, then
  $E\in\Lie{g}^\mc{I}\setminus\{0\}$ implies that $\tau(E)^\mc{I}\neq  0$.
\end{itemize}
\end{cor}
\begin{proof} (i) Let $\alpha$ and $\beta$ two roots in $\mc{I}$ such  that $\alpha + \beta$
is a root. Then $\alpha + \beta \in \mc{S}_\mu$ for every  $\mu\in\mc{E}$,
because each shadow defines a subalgebra. This implies that $\alpha +
\beta\in\mc{I}$, so that $\Lie{n}^\mc{I}$ is a subalgebra in $\Lie{n}$.
By Lemma~\ref{niliwa}, $\Lie{n}^\mc{I}$ is an Iwasawa nilpotent Lie  algebra.
\par
(ii) Because of (i), we can
use the results in \cite{CDKR2} in order to describe the multicontact  vector
fileds on $\group{N}^\mc{I}$. Thus, the same argument as in the
proof of  Proposition~\ref{rappre} with
$\mc{I}$ in place of $\mc{S}_\mu$ shows that $\tau(E)^\mc{I}\in MC(  \group{N}^\mc{I})$.
\par
(iii) See the proof of  Proposition~\ref{rappre}.
\end{proof}
Notice that (ii) of the previous corollary asserts  $\tau(\Lie{g}^{\mc{I}})\subset MC( \group{N}^\mc{I})$ and
does not claim equality. This is because the intersection $\mc{I}$ may  well give rise to  rank-one algebras, so that
the results of \cite{CDKR2} do not apply to prove the reverse inclusion.
 Similarly, (iii) gives only one implication.
\par
We shall conclude the proof of Theorem~\ref{al} by setting up a diagram
of the following type
\begin{equation}
\xymatrix{MC(\group{S})\ar[r]^{\!\!\!\!\!\!\!\!\!k}&\Lie{g}^{1}\oplus \dots \oplus\Lie{g}^{p}\\
&\Lie{q}/\Lie{n}_{\mc{C}}\ar@{-->}[ul]\ar[u]_{\bar\ell}}
\label{diagram}\end{equation}
\noindent
and then showing that both $k$ and $\bar\ell$ are injective linear maps
with equal images. Thus $MC(\group{S})$ and $\Lie{q}/\Lie{n}_{\mc{C}}$
will be seen to be isomorphic vector spaces. Moreover, the induced left
arrow  coincides with the quotient map of the map
$E\mapsto\overline{\tau(E)}$ defined on
$\Lie{q}$, which is a Lie algebra homomorphism. This latter fact is
a straightforward consequence of the definition of the various maps,
that we illustrate below.
\par
Fix a numbering
$\mu_1,\dots,\mu_p$ of the maximal roots and write $\Lie{g}^i$ for
$\Lie{g}^{\mu_i}$. By Proposition~\ref{compatibilita}, we can
associate to each
$F\in MC(\group{S})$ a vector $(F^{1},\dots,F^{p})$, where each
$F^i=F^{\mu_i}\in MC(\group{N}^{\mu_i})$ is the natural projection.
Moreover, by Proposition~\ref{rappre}, we know that
$F^i=\tau(E)^{i}$ for some
$E\in
\Lie{g}^i$. These observations allow us to define a map
$$
k:MC(\group{S})\longrightarrow\Lie{g}^{1}\oplus\dots\oplus\Lie{g}^{p}
$$
by setting $k(F)=(F^{1},\dots,F^{p})$. Here a few comments are in order.
First of all, the direct sum on the right is viewed merely as a vector
space. Secondly,  as already observed, we identify each
$\Lie{g}^\mu$ with $\{\tau(E)^\mu,E\in\Lie{g}^\mu\}$, by means of
Proposition~\ref{rappre}.
Finally, the map $k$ is injective because $(F^{1},\dots,F^{p})$ encodes  all
components of $F$.
\par
Next, observe that the assignment
$$
E\mapsto(\tau(E)^1,\dots,\tau(E)^p)
$$
gives a well-defined map
$$
\ell:\Lie{q} \longmapsto \Lie{g}^{1}\oplus\dots\oplus\Lie{g}^{p}
$$
because by Proposition~\ref{rappre} $\tau(E)^i\not=0$ only if
$E\in\Lie{g}^{i}$. The kernel of $\ell$ is $\Lie{n}_{\mc{C}}$
because if all $\tau(E)^i$ vanish, then so does
$\overline{\tau(E)}$, and Theorem~\ref{andata} says
$E\in\Lie{n}_{\mc{C}}$. Therefore $\ell$ projects to a map $\bar\ell$
defined on $\Lie{q}/\Lie{n}_{\mc{C}}$.
\par
All the ingredients appearing in the
diagram~\eqref{diagram} have been defined
and all we need to show is that
$k(MC(\group{S}))=\bar\ell(\Lie{q}/\Lie{n}_{\mc{C}})$.
This will be formalized in
Proposition~\ref{image} below. Its proof needs a technical lemma that
characterizes $\Lie{q}$ in terms of roots. We underscore
that this characterization (i.e. Lemma~\ref{norma}) holds {\it only}  under
the hypothesis~(I) in Theorem~\ref{al}. Before we set up all the needed
machinery, we point out the main issue. Recall that $\Lie{q}$ is
the normalizer of the nilpotent ideal $\Lie{n}_{\mc{C}}$ in $\Lie{g}$,  so
that, as we already observed, it clearly contains
$\Lie{m}\oplus\Lie{a}\oplus\Lie{n}$. Breaking it according to the
root space decompostion of $\Lie{g}$, for short
$\Lie{g}=\overline{\Lie{n}}\oplus\Lie{m}\oplus\Lie{a}\oplus\Lie{n}$,
we may write
$\Lie{q}=(\Lie{q}\cap\overline{\Lie{n}})
\oplus\left(\Lie{m}\oplus\Lie{a}\oplus\Lie{n}\right)$.
Since $\Lie{q}$ is a parabolic subalgebra of $\Lie{g}$, then
$\Lie{q}\cap\overline{\Lie{n}}=\sum_{\alpha\in\mc{D}}\Lie{g}_\alpha$,  for
some $\mc{D}\subset \Sigma_-$ (see \cite{KN02}, Sec.7, Ch.VII ).
The point addressed by Lemma~\ref{norma} is
an  adequate description of $\Lie{q}\cap\overline{\Lie{n}}$.
\par
Given a root $\alpha = \sum_{\delta\in\Delta}n_\delta
(\alpha) \delta$ we denote by $\mc{Y}(\alpha)$ the subset of $\Delta$
consisting of those
$\delta$ for which
$n_{\delta} (\alpha)\neq 0$, and  we call it the {\it simple support} of
$\alpha$. For later use, we recall
the following result (see Corollary~3,~Ch.~VI,~§~1,~pag.~160
in [Bourbaki]).
\begin{prop}[Bourbaki]\label{bour}
\begin{itemize}
\item[(i)]
Let $\alpha\in\Sigma$. Then $\mc{Y}(\alpha)$ is a connected subset of  the
Dynkin diagram associated to
$\Sigma$.
\item[(ii)]
Let $\mc{Y}$ be any connected non empty subset of a Dynkin diagram. Then
$\sum_{\beta\in \mc{Y}} \beta$ is a root.
\end{itemize}
\end{prop}
Yet another piece of notation. We say that a simple root $\delta$ is a  {\it
boundary} simple root if there exists a maximal root $\nu$ in $\mc{R}$
whose simple support is a connected diagram that does not contain  $\delta$
but to which $\delta$ is adjacent, i.e. such that there exists
$\delta^\prime\in\mc{Y}(\nu)$ with the property that
$\delta+\delta^\prime$ is a root. The set of all the boundary simple  roots
will be denoted by $\mc{B}$.
\vskip6.1truecm
\setlength{\unitlength}{1cm}
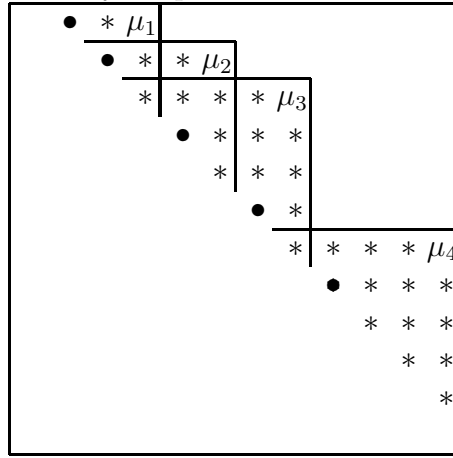
\begin{picture}(0,0)\thicklines
\put(3,0.5){\line(0,1){6}}
\put(3,0.5){\line(1,0){6}}
\put(3,6.5){\line(1,0){6}}
\put(9,0.5){\line(0,1){6}}
\put(5,5){\line(0,1){1.5}}
\put(4,6){\line(1,0){2}}
\put(6,4){\line(0,1){2}}
\put(4.5,5.5){\line(1,0){2.5}}%
\put(7,3){\line(0,1){2.5}}%
\put(6.5,3.5){\line(1,0){2.5}}%
\put(3.7,6.15){$\bullet$}
\put(4.2,6.15){$*$}
\put(4.55,6.15){$\mu_1$}
\put(4.2,5.65){$\bullet$}
\put(4.7,5.65){$*$}
\put(5.2,5.65){$*$}
\put(5.55,5.65){$\mu_2$}
\put(4.7,5.15){$*$}
\put(5.2,5.15){$*$}
\put(5.7,5.15){$*$}
\put(6.2,5.15){$*$}
\put(6.55,5.15){$\mu_3$}
\put(5.2,4.65){$\bullet$}
\put(5.7,4.65){$*$}
\put(6.2,4.65){$*$}
\put(6.7,4.65){$*$}
\put(5.7,4.15){$*$}
\put(6.2,4.15){$*$}
\put(6.7,4.15){$*$}
\put(6.2,3.65){$\bullet$}
\put(6.7,3.65){$*$}
\put(6.7,3.15){$*$}
\put(7.2,3.15){$*$}
\put(7.7,3.15){$*$}
\put(8.2,3.15){$*$}
\put(8.55,3.15){$\mu_4$}
\put(7.2,2.65){$\bullet$}
\put(7.2,2.65){$*$}
\put(7.7,2.65){$*$}
\put(8.2,2.65){$*$}
\put(8.7,2.65){$*$}
\put(7.7,2.15){$*$}
\put(8.2,2.15){$*$}
\put(8.7,2.15){$*$}
\put(8.2,1.65){$*$}
\put(8.7,1.65){$*$}
\put(8.7,1.15){$*$}
\end{picture}
\par
In the picture above the dark dots are the boundary roots.
Each shadow determines a connected line in the $A_{n}$-type Dynkin  diagram,
adjacent to which lie exactly two boundary roots if the line is inside  the diagram, otherwise just one.
(Recall that this refers to $\Lie{sl}(n,\mathbb{R})$).
\begin{lemma}\label{norma}
Let
$\Lie{q}\cap\overline{\Lie{n}}
=\sum_{\alpha\in\mc{D}}\Lie{g}_\alpha$.
\begin{itemize}
\item[(i)] If $\delta$ is a simple root, then
$-\delta \notin \mc{D}$ if and only if $\delta \in \mc{B}$.
\item[(ii)] If $\alpha$ is any positive root, then
$-\alpha\notin\mc{D}$ if and only if the simple support of $\alpha$
contains a simple root in
$\mc{B}$.
\end{itemize}
\end{lemma}
\begin{proof} We prove (i) first.
\par
``$\Leftarrow$''. Let $\delta\in\mc{B}$ and let $\nu$ be a
maximal root to whose shadow $\delta$ is adjacent.
Proposition~\ref{bour} implies that $\sum_{\epsilon\in\mc{Y}(\nu)}\eps
+\delta =\sigma+\delta$ is a root. Moreover, it does not lie in  $\mc{R}$.
Indeed, if $\sigma+\delta\in\mc{R}$, then it would belong to a
shadow containing $\mc{S}_\nu$, contradicting the maximality of $\nu$.
On the other hand, $\sigma$ itself is a root, again by
Proposition~\ref{bour}, and it lies in $\mc{R}$, because it is sum of
simple roots in a same shadow $\mc{S}_\nu$. Thus, we found a root
in $\mc{C}$, namely $\sigma + \delta$, such that $(\sigma +
\delta)-\delta\notin\mc{C}$. Therefore $-\delta\notin\mc{D}$.
\par
``$\Rightarrow$''. Suppose
$\delta\notin\mc{B}$. Let $\alpha\in\mc{C}$ with $\delta\prec\alpha$ and
consider its simple support $\mc{Y}(\alpha)$. We shall show that
$\alpha-\delta\in\mc{C}$ whenever $\alpha-\delta\in\Sigma$. Take a  maximal
connected set $\mc{F}$ of simple roots in
$\mc{Y}(\alpha)$ with the following properties:
\begin{itemize}
\item[$\diamond$] $\delta\in\mc{F}$;
\item[$\diamond$] there exists a shadow containing $\mc{F}$.
\end{itemize}
This
means that
$\delta\in\mc{F}\subset\mc{S}_\nu$ for some $\nu$, but no
larger connected subset of $\mc{Y}(\alpha)$ containing $\delta$  is
contained in any other single shadow. Necessarly $\mc{F}$ is a proper
subset of
$\mc{Y}(\alpha)$, for otherwise $\alpha$ would lie in $\mc{R}$. Take
$\eps\in\mc{Y}(\alpha)$ adjacent to $\mc{F}$. Then two cases arise.
\begin{itemize}
\item[(a)] $\mc{Y}(\alpha -\delta)$ does contain $\delta$. In this case
$\mc{Y}(\alpha -\delta)$ contains both $\mc{F}$ and $\eps$. Thus $\alpha
-\delta\notin
\mc{R}$, for otherwise $\mc{F}\cup\{\eps\}$ would be a connected set
contained in a single shadow (namely any shadow containing
$\alpha -\delta$) and
it would be larger than $\mc{F}$.
\item[(b)] $\mc{Y}(\alpha -\delta)$ does not contain $\delta$. Then
$\mc{Y}(\alpha -\delta)$ is connected and $\delta$ is adjacent to it. If
$\alpha -
\delta\in\mc{R}$ then $\delta$ would be a boundary root because
$\mc{Y}(\alpha -\delta)\subset \mc{S}_\nu$ for some maximal root
$\nu$, and $\delta\notin\mc{S}_\nu$ (for otherwise $\alpha =(\alpha
-\delta)+\delta \in \mc{S}_\nu$, which is impossible). Hence $\delta$  would
be adjacent to the simple support of
$\mc{S}_\nu$, contradicting $\delta\not\in\mc{B}$. Therefore $\alpha
-\delta \notin \mc{R}$.
\end{itemize}
We have seen that in all cases $-\delta\in\mc{D}$. This concludes the
proof of (i).
\par
As for (ii), take a non simple root $-\alpha\notin\mc{D}$. Then
$\mc{Y}(\alpha)$ contains at least one simple root $\delta\notin
-\mc{D}$. Indeed, since
$\Lie{q}$ is a subalgebra, if  $\mc{Y}(\alpha)$ were contained in  $-\mc{D}$, then
$\alpha$ itself would lie in $\Lie{q}$. Thus $\mc{Y}(\alpha)$
contains a boundary simple root. Conversely, if
$\alpha\in\Sigma_+$ is such that $\mc{Y}(\alpha)$ contains a simple root
in $\mc{B}$, then it contains a simple that is not $-\mc{D}$, so that
$-\alpha$ is not in $\mc{D}$.
\end{proof}

\begin{prop}\label{image} In the notations of diagram~\eqref{diagram}
$$
(F^1,\dots,F^p)\in k(MC(\group{S})) \Leftrightarrow
(F^1,\dots,F^p)=\ell(E)\hskip0.2cm\text{ for some }E\in\Lie{q}.
$$
\end{prop}
\begin{proof}
``$\Leftarrow$''. Let $E\in\Lie{q}$. Then
$\ell(E)=(\tau(E)^1,\dots,\tau(E)^p)$. By Theorem~\ref{andata}, it
follows that $\overline{\tau(E)}\in MC(\group{S})$. Therefore
$$
k(\overline{\tau(E)})=(\tau(E)^1,\dots,\tau(E)^p)=\ell(E).
$$
\par
``$\Rightarrow$''. Let $F\in MC(\group{S})$ and $k(F)=(F^1,\dots,F^p)$.
 Proposition~\ref{rappre} implies that for
all
$i=1,\dots,p$
there exists $E^i\in\Lie{g}^i$ such that
$F^i=\tau(E^i)^i$.
Write $E^i=\sum_{\alpha\in\Sigma^i\cup\{0\}}
E^i_\alpha$, with $\Sigma^i=\mc{S}_{\mu_i}\cup(-\mc{S}_{\mu_i})$.
By definition,  $\tau(E^i)^i\in MC(\group{N}^{\mu_i})$ if and only if
$\tau(E^i_\alpha)^i \in MC(\group{N}^{\mu_i})$ for every
$\alpha\in\Sigma^i\cup\{0\}$.
\par
Recall that  $\Lie{q}=\Lie{m}\oplus\Lie{a}\oplus\Lie{n}\oplus(\overline{\Lie{n}}\cap\Lie{q})$, and
write $\Lie{g}_0=\Lie{m}\oplus\Lie{a}$,  $\Lie{n}=\oplus_{\gamma\in\Sigma_+}\Lie{g}_\gamma$ and
$\overline{\Lie{n}}\cap\Lie{q}=\oplus_{\gamma\in\mc{D}}  \Lie{g}_\gamma$. Therefore, the normalizer can be
written as follows:
$$
\Lie{q}=\bigoplus_{\alpha\in\mc{G}}\Lie{g}_\alpha,
$$
where $\mc{G}=\Sigma_+ \cup \{0\} \cup \mc{D}$. Using these notations,  we shall prove the following two claims:
\begin{itemize}
\item[(a)] $\alpha\in\mc{G}\Rightarrow E_\alpha^i = E_\alpha^j$, for  every $i,j$;
\item[(b)] $\alpha\notin\mc{G}\Rightarrow E^i_\alpha =0$.
\end{itemize}
These two facts allows us to define an element  $E=\sum_{\alpha\in\Sigma\cup\{0\}}E_\alpha$ by
$$
E_\alpha=
\begin{cases}
E_\alpha^i&\text{ if }\alpha\in\mc{G}\\
0&\text{ if }\alpha\notin\mc{G},
\end{cases}
$$
for all $i=1,\dots,p$.
In particular, $E\in\Lie{q}$ and $\ell(E)=(F^1,\dots,F^p)$, that proves  the proposition.
\par
$(a)$ If $\alpha\in\mc{G}$, then $E_\alpha^i\in \Lie{q}$ for every  $i=1,\dots,p$. By Theorem~\ref{andata},
$\overline{\tau(E_\alpha^i)}\in MC(\group{S})$ and, by  Proposition~\ref{rappre}, $\tau(E_\alpha^i)\in
MC(\group{N}^\mu)$ for every maximal root $\mu$. Moreover,  Proposition~\ref{rappre} also implies that
$\tau(E_\alpha^i)^j\neq 0$ if and only if  $E_\alpha^i\in\Lie{g}^{\mu_j}$. Suppose that $E_\alpha^i$ belongs to
$\Lie{g}^{\mu_j}$ with $j\neq i$ and let  $\mc{I}=\mc{S}_{\mu_i}\cap\mc{S}_{\mu_j}$ ($\mc{I}$ is not empty,
otherwise $\Lie{g}^{\mu_i}$ and $\Lie{g}^{\mu_j}$ would not have a  common element). Then statement (iii) of
Corollary~\ref{inters} implies that the components of
$\tau(E^i_\alpha)^i$ labeled by
$\mc{I}$ do not vanish identically. This forces $F^j\not=0$,
because
$$
\tau(E^j_\alpha)^\mc{I} = \tau(E_\alpha^i)^\mc{I}\neq 0.
$$
Moreover, since
$\Lie{g}_\beta\subset\Lie{q}$, the identity
$\tau(E^j_\alpha-E_\alpha^i)^\mc{I}= 0$
holds only if
$E_\alpha^i=E^j_\alpha$, again by (iii) in
Corollary~\ref{inters}. This proves (a).
\par
(b). Let $\alpha\notin\mc{G}$, and suppose that $E_\alpha^i\neq 0$.
We show that this hypothesis takes us to a contradiction. In  particular, we shall
show that in the vector $(F^1,\dots,F^p)$ appears one component that is  not of multicontact type, there implying
that $F$ itself is not a multicontact vector field.
\par
By definition of $\mc{G}$, the root $\alpha$
must be negative. Furthermore, by (ii) of Lemma~\ref{norma}, there  exists
$\delta\in\mc{B}$ such that $\delta + \delta_1 +\dots+\delta_q=-\alpha$.
Let $\mc{S}_{\mu_j}$ be a shadow to which $\delta$ is adjacent.
Then there
exists at least a shadow to which $\delta$ belongs that
intersects
$\mc{S}_{\mu_j}$. Indeed, if this does not happen, then $\delta$ would
belong to a dark zone disjoint from $\mc{S}_{\mu_j}$, which is  impossible.
Call $\mc{S}_{\mu_k}$ such a shadow and
$\mc{J}=\mc{S}_{\mu_j}\cap\mc{S}_{\mu_k}\neq\emptyset$.
\vskip6.8truecm
\setlength{\unitlength}{1cm}
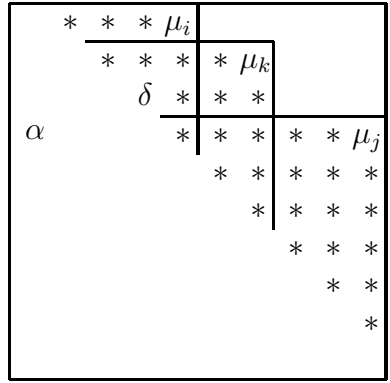
\begin{picture}(0,0)\thicklines
\put(3.5,1){\line(0,1){5}}
\put(3.5,1){\line(1,0){5}}
\put(3.5,6){\line(1,0){5}}
\put(8.5,1){\line(0,1){5}}
\put(6,4){\line(0,1){2}}
\put(4.5,5.5){\line(1,0){2.5}}%
\put(7,3){\line(0,1){2.5}}%
\put(5.5,4.5){\line(1,0){3}}%
\put(4.2,5.65){$*$}
\put(4.7,5.65){$*$}
\put(5.2,5.65){$*$}
\put(5.55,5.65){$\mu_i$}
\put(4.7,5.15){$*$}
\put(5.2,5.15){$*$}
\put(5.7,5.15){$*$}
\put(6.2,5.15){$*$}
\put(3.7,4.2){$\alpha$}
\put(6.55,5.15){$\mu_k$}
\put(5.2,4.65){$\delta$}
\put(5.7,4.65){$*$}
\put(6.2,4.65){$*$}
\put(6.7,4.65){$*$}
\put(5.7,4.15){$*$}
\put(6.2,4.15){$*$}
\put(6.7,4.15){$*$}
\put(7.2,4.15){$*$}
\put(7.7,4.15){$*$}
\put(8.05,4.15){$\mu_j$}
\put(6.2,3.65){$*$}
\put(6.7,3.65){$*$}
\put(7.2,3.65){$*$}
\put(7.7,3.65){$*$}
\put(8.2,3.65){$*$}
\put(6.7,3.15){$*$}
\put(7.2,3.15){$*$}
\put(7.7,3.15){$*$}
\put(8.2,3.15){$*$}
\put(7.2,2.65){$*$}
\put(7.7,2.65){$*$}
\put(8.2,2.65){$*$}
\put(7.7,2.15){$*$}
\put(8.2,2.15){$*$}
\put(8.2,1.65){$*$}
\end{picture}
\vskip0.1truecm
In the picture above we illustrate, in the case of  $\Lie{sl}(n,\mathbb{R})$, a situation that
corresponds to
what we just described. We are going to show that a multicontact vector  field corresponding
to the root $\alpha$ cannot be identically zero in its components  labeled by the intersection $\mc{J}$ ($
\mc{S}_{\mu_j}\cap\mc{S}_{\mu_k}$). The reason of this lies in the fact  that
$\mc{Y}(\alpha)$ contains a boundary simple root, namely $\delta$, and  that $\delta$ is adjacent
to $\mc{J}$. Roughly speaking, the root $\alpha$ is  close enough to  $\mc{J}$, and this allows a
multicontact vector field corresponding to $\alpha$ not to be killed  before its coefficients arrive to
the $\mc{J}$-positions.
\par
In short, we prove that
\begin{equation}\label{zero}
\tau(E^i_\alpha)^\mc{J}\neq 0.
\end{equation}
If the equation above holds, then the relation
$$
\tau(E^i_\alpha)^\mc{J}=\tau(E^j_\alpha)^\mc{J}
$$
forces $F^j=\tau(E^j)^{j}$ to be non-zero
because $\tau(E_\alpha^j)^{j}\not=0$. On the other hand
$-\alpha\notin\mc{S}_{\mu_j}$, for otherwise
$\delta$ would lie in $\mc{S}_{\mu_j}$. This implies that  $\tau(E_\alpha^j)^j$ is not in
$MC(\group{N}^{{\mu_j}})$ by Proposition~\ref{rappre}. This, in turn,  implies that $F^{\mu_j}$, hence $F$, is
not a multicontact vector field, that is the contradiction we expected.
\par
It remains to prove equation~\eqref{zero}.
Suppose $\tau(E^i_\alpha)^\mc{J}=0$. This will give that  $\tau(E^i_{-\delta})^\mc{J}=0$ that, in turn, implies
that $\delta$ is not a boundary root, a contradiction.
First, by $\tau(E^i_\alpha)^\mc{J}=0$, it follows that
for every
$E^\prime\in\Lie{n}$ is
$$
[\tau(E^i_\beta), \tau(E^\prime)]=[\sum_{\gamma_1 \in\mc{J}^c}
\sum_{i=1}^{m_{\gamma_1}} f_{\gamma_1,i}
X_{\gamma_1,i},
\sum_{\gamma_2\in \mc{J}} \sum_{j=1}^{m_{\gamma_2}} g_{\gamma_2,j}
X_{\gamma_2,j} +
\sum_{\gamma_3\in\mc{J}^c} \sum_{k=1}^{m_{\gamma_3}} g_{\gamma_3,k}
X_{\gamma_3,k}].
$$
All terms of the bracket above lie in  $\mathfrak{X}(\group{N}^{\mc{J}^c})$, except
$$
f_{\gamma_1,i} X_{\gamma_1,i} (g_{\gamma_2,j}) X_{\gamma_2,j},
$$
but $X_{\gamma_1,i} (g_{\gamma_2,j})=0$, for every
$\gamma_1\in\mc{J}^c$ and $\gamma_2\in\mc{J}$,
because the coefficients $g_{\gamma_2,j}$ are
$(\Sigma_+\setminus\mc{J})$-independent. Indeed, $\mc{J}$ is a  Hessenberg set and its complement
$\mc{J}^c$ defines an ideal $\Lie{n}_{\mc{J}^c}$ in $\Lie{n}$ whose  normalizer contains $\Lie{n}$.
Therefore, since $E^\prime \in \Lie{n}$, it also lies in $N_{\Lie{g}}  \Lie{n}_{\mc{J}^c}$, so that the
coefficients of $\tau(E^\prime)^\mc{J}$ are  $(\Sigma_+\setminus\mc{J})$-independent by Lemma~\ref{cindip}.
Hence  $[\tau(E^i_\alpha),\tau(E^\prime)]\in\mathfrak{X}(\group{N}^{\mc{J}^c})$ , that implies
$$
\tau([E^i_\alpha,E^\prime])^\mc{J}={[\tau(E^i_\alpha),\tau(E^\prime)]}^\mc{J} =0
$$
for every $E^\prime\in\Lie{n}$. The same argument can be iterated for  showing that
\begin{equation}\label{restr}
\tau([[E^i_\alpha,E^\prime],\dots,E^{{\it
(n)}}])^\mc{J}={[[\tau(E^i_\alpha),\tau(E^\prime)],\dots,\tau(E^{{\it  (n)}})]}^\mc{J} =0
\end{equation}
for every collection of elements $E^\prime,\dots,E^{{\it (n)}}$ in  $\Lie{n}$.
\par
Let $\delta_1,\dots,\delta_q$ simple roots such that $\alpha  +\delta_1+\dots+\delta_q=-\delta$ is a chain, and
 $E_1\in\Lie{g}_{\delta_1},\dots,E_q\in\Lie{g}_{\delta_q}$ such that
$$
[[E^i_\alpha,E_1],\dots,E_q]=E_{-\delta}\in\Lie{g}_{- \delta}\setminus\{0\}.
$$
We apply~\eqref{restr} to the bracket above, there obtaining
$$
\tau(E_{-\delta})^\mc{J}=0.
$$
Since $\theta E_{-\delta}\in\Lie{n}$, again the formula~\eqref{restr}  toghether with Prop 6.52 in \cite{KN02}
gives
\begin{equation}\label{uno}
0=[\tau(E_{-\delta}),\tau(\theta E_{-\delta})]^\mc{J}=
B(E_{-\delta},\theta E_{-\delta})\tau(H_\delta)^\mc{J}.
\end{equation}
By Lemma \ref{norma}, since $\delta\in\mc{B}$, there exists a simple  root
$\delta^\prime\in \mc{S}_{\mu_j}$ such that $\delta + \delta^\prime$ is  a
root. This implies that
$\langle \delta,\delta^\prime\rangle \neq 0$, because
$\delta -\delta^\prime$ is never a root. Hence $\delta^\prime
(H_\delta)\neq 0$, so that  $H_\delta\in\Lie{g}^{\mu_j}\cap\Lie{g}^{\mu_k}$. By
Corollary~\ref{inters}, it follows that $\tau(H_\delta)^\mc{J}\neq 0$,
contradicting~\eqref{uno}.
This concludes our proof.
\end{proof}
\subsection{A remark on the group}
In this section we  show that the group $Q=\Int (\Lie{q})$ acts on S via
multicontact mappings. In order to do that, we define an alternative
model for the Hessenberg manifolds, which is compatible with the
stratified structure introduced in the first section.
\par
Consider the group $\group{N}_\mc{C}=\exp\Lie{n}_\mc{C}$ . Since  $\Lie{n}_\mc{C}$ is an
ideal in $\Lie{n}$, its exponential group is a normal subgroup of N.  Therefore the quotient
$\group{N}/\group{N}_\mc{C}$ is a nilpotent Lie group.
We identify the Lie algebra of $\group{N}/\group{N}_\mc{C}$ with  $\Lie{n}/\Lie{n}_\mc{C}$, and we define
a natural multicontact structure on this quotient simply considering  the subbundles $\{\langle \Lie{g}_\delta
\rangle_{\Lie{n}_\mc{C}}: \delta \in \Delta\}$, where $\langle E  \rangle_{\Lie{n}_\mc{C}}$ denotes the coset
of $E$ in $\Lie{n}/\Lie{n}_\mc{C}$.
Let $f$ be a diffeomorphism between open subsets of  $\group{N}/\group{N}_\mc{C}$. Then $f$
is a multicontact mapping if for every simple root $\delta$
$$
f_* (\langle \Lie{g}_\delta \rangle_{\Lie{n}_\mc{C}})\subset \langle  \Lie{g}_\delta
\rangle_{\Lie{n}_\mc{C}}.
$$
\par
The coordinates system on
the slice S define the analytic structure on  $\group{N}/\group{N}_\mc{C}$.
Thus,
$\group{N}/\group{N}_\mc{C}$ and S are diffeomorphic by the assignment
$$
\chi :
\langle (\{x_{\alpha,i}\}_{\alpha\in\Sigma_+})  \rangle_{\group{N}_\mc{C}} \mapsto
(\{x_{\alpha,i}\}_{\alpha\in\mc{R}},0),
$$
where $\langle n \rangle_{\group{N}_\mc{C}}$ denotes the coset of  $n\in\group{N}$ in the
quotient group. The differential $\chi_*$  maps the
left--invariant vector field
$\langle X_{\alpha,i} \rangle_{\Lie{n}_\mc{C}}$ to
$\overline{X}_{\alpha,i}$, therefore the multicontact structure on  $\group{N}/\group{N}_\mc{C}$ is mapped onto
the multicontact structure on S. The diffeomorphism $\chi$ allows us to  view S, and hence locally a Hessenberg
manifold, as a nilpotent Lie group. From the point of view of  multicontact mappings, we identify S with
$\group{N}/\group{N}_\mc{C}$.
\par
Let $\group{Q}=
{\rm Int}(\Lie{q})$. We have $ {\rm Int}(\Lie{q})\subset {\rm  Int}(\Lie{g})$,
because $ {\rm Int}(\Lie{q})=e^{{\rm ad}\Lie{q}}$, ${\rm  Int}(\Lie{g})=e^{{\rm ad}\Lie{g}}$
and $\Lie{q}\subset\Lie{g}$.
\begin{lemma}\label{action}
The action of every element $q\in\group{Q}$ on $\group{N}$ induces a well-posed  action on the quotient
$\group{N}/\group{N}_\mc{C}$, namely
$$
\hat{q}(\langle n \rangle_{\group{N}_\mc{C}}) =\langle [qn]  \rangle_{\group{N}_\mc{C}},
$$
where $[qn]$ is the $\group{N}$-component of $qn$ in the Bruhat decomposition.
\end{lemma}
\begin{proof}
Let $n\in\group{N}$ and
$n_\mc{C}\in\group{N}_\mc{C}$.
Then $n$ and $nn_\mc{C}$ both represent
$\langle n \rangle_{\group{N}_\mc{C}}\in\group{N}/\group{N}_\mc{C}$. We  show that
$[qn]$ and
$[qnn_\mc{C}]$ represent the same element in
$\group{N}/\group{N}_\mc{C}$, that is  $[qn][qnn_\mc{C}]^{-1}\in\group{N}_\mc{C}$.
Let $p\in\group{P}$ such that $[qn]=qnp$. Since
$\group{N}_\mc{C}$ is a
normal subgroup of Q, there exists ${n}^\prime_\mc{C} \in  \group{N}_\mc{C}$ such that
\begin{align*}
[qnn_\mc{C}]&= [n_\mc{C}^\prime
qn]\\
&=n_\mc{C}^\prime [qn] \\
&=n_\mc{C}^\prime qnp.
\end{align*}
Then $[qn][qnn_\mc{C}]^{-1}=qnp (n_\mc{C}^\prime
qnp)^{-1}=(n_\mc{C}^\prime)^{-1}\in\group{N}_\mc{C}$, as required.
\end{proof}
 We prove the following
proposition.
\begin{prop} \label{mcg}
Let Q be as above, and $\mc{A}$ an open subset of
$\group{N}/\group{N}_\mc{C}$. For
every $q\in\group{Q}$, the map
$$
\hat{q}: \mc{A}\subset \group{N}/\group{N}_\mc{C}\rightarrow  \group{N}/\group{N}_\mc{C}
$$
is a multicontact mapping on $\mc{A}$.
Furthermore $\hat{q}=\id_\mc{A}$ for every
$q\in\group{N}_\mc{C}$.
\end{prop}
\begin{proof}
Since $q\in\group{G}=\Int(\Lie{g})$, it is a multicontact mapping on  G/P. Thus, $q_*
(\Lie{g}_\delta)\subseteq
\Lie{g}_\delta$ for every simple
root $\delta$ (see Ch.~\ref{pr}, Sec.~\ref{mm}).
Let $E\in\Lie{g}_\delta$, for some $\delta\in\Delta$, and consider a  representative
in $\Lie{n}/\Lie{n}_\mc{C}$ of $\langle E \rangle_{\Lie{n}_\mc{C}}$,  say $E + E^\prime$,
with $E^\prime\in\Lie{n}_\mc{C}$.
Then
$$
\hat{q}_* (\langle E \rangle_{\Lie{n}_\mc{C}})=  \langle(l_q)_{*}(E+E^\prime)\rangle_{\Lie{n}_\mc{C}}.
$$
By definition
$$
(l_q)_{*}(E^\prime)=\frac{d}{dt} (q\exp(tE^\prime))\Big|_{t=0}.
$$
Since $[\Lie{q},\Lie{n}_\mc{C}]\subset\Lie{n}_\mc{C}$, a  straightforward calculation implies
that $(l_q)_{*}(E^\prime)\in\Lie{n}_\mc{C}$.
Therefore there exists $E^{\prime\prime}\in\Lie{n}_\mc{C}$ such that
$$
\hat{q}_* (\langle E \rangle_{\Lie{n}_\mc{C}})
= \langle(l_q)_{*e}(E+E^\prime)\rangle_{\Lie{n}_\mc{C}}=\langle q_*(E) +
E^{\prime\prime}\rangle_{\Lie{n}_\mc{C}}
\subset \langle\Lie{g}_\delta +\Lie{n}_\mc{C}  \rangle_{\Lie{n}_\mc{C}}\subset \langle \Lie{g}_\delta
\rangle_{\Lie{n}_\mc{C}}.
$$
Since $\langle n_\mc{C} n \rangle_{\group{N}_\mc{C}}=\langle n  \rangle_{\group{N}_\mc{C}}$,
it follows that $\hat{n}_\mc{C}$ maps $\langle n  \rangle_{\group{N}_\mc{C}}$ in itself, for every
$n\in\group{N}$. Hence the proposition holds.
\end{proof}
The above proposition, toghether with the diffeomorphism $\chi$,
tell us that $\group{Q}/\group{N}_{\mc{C}}$ is a group of multicontact  mappings
on S. Indeed, for every $q\in \group{Q}$, the map  $\chi\hat{q}\chi^{-1}$ is multicontact
on S, and it coincides with the identity whenever  $q\in\group{N}_\mc{C}$.
This fact implies in particular that $\Lie{q}/\Lie{n}_\mc{C}$ is a Lie  algebra of multicontact
vector fields, so that we obtain another proof of Theorem~\ref{andata}.

\subsection{The case of $A_l$}

We conclude this section observing that hypothesis (I) of  Theorem~\ref{al} always holds if
we consider Hessenberg subsets in the root
system $A_l$.
The underlying vector space of $A_l$ is $V=\{ v\in \mathbb{R}^{l+1} :
\langle v, e_1
+\dots +e_{l+1} \rangle=0\}$, and the roots are $A_l=\{e_1 - e_j :
i\neq j\}$, where we
choose the positive ones fixing $\Sigma_+=\{e_i -e_j :i<j\}$. A basis
of simple roots is
given by
$\Delta=\{\delta_1=e_1-e_2 ,\delta_2 =e_2 -e_3,\dots,\delta_l =e_l -
e_{l+1} \}$, and the
highest root with respect to this basis is $\omega =\delta_1 + \dots
+ \delta_l$.
In~\cite{C1},~\cite{C2},~\cite{C3} the author gives a classification of simple Lie
algebras. In particular, by \cite{C2} we desume that the restricted root spaces
associated to a simple Lie
algebra with root
system $A_l$ have the
same dimension, that is all roots have the same
multiplicity, and
the
classification is the following:
\begin{itemize}
\item[(1)] $m_\omega
=1$, that corresponds to
$(\Lie{sl}(l+1,\mathbb{R})$,
$\Lie{so}(l+1))$;
\item[(2)] $m_\omega =2$, that corresponds
to
$(\Lie{sl}(l+1,\mathbb{C})$, $\Lie{su}(l+1))$;
\item[(3)]
$m_\omega =4$, that corresponds to
$(\Lie{sl}(l+1,\mathbb{H})$,
$\Lie{sp}(l+1))$;
\item[(4)] $m_\omega =8$, possible only for $l=2$,
it
corresponds to
$(\Lie{e}_{(6,-26)}$,
$\Lie{f}_4)$.
\end{itemize}
If $l=2$, the Hessenberg proper subsets
in $A_l$ are $\mc{R}=\{\delta_1\}$,
$\mc{R}=\{\delta_2\}$,
$\mc{R}=\{\delta_1,\delta_2\}$. In all these
cases the study
of
multicontact mappings reduces to the case of rank one
nilpotent
Iwasawa Lie algebras. Because of hypothesis (II) of Theorem~\ref{al}
and the remark thereafter,
 we restrict ourselves
to the case $l>2$.
We have a rather obvious
consequence.
\begin{prop}
Let $\Lie{g}$ a simple Lie algebra with
root system $A_l$, $l>2$.
Let $\mc{R}\subseteq \Sigma_+$ a subset of
Hessenberg type. Then each
shadow $\mc{S}_\mu$
define an Iwasawa
nilpotent subalgebra of $\Lie{n}$ of the following
form
$$
\Lie{n}^\prime=\bigoplus_{\alpha\in \mc{S}_\mu} \Lie{g}_\alpha.
$$
\end{prop}
\begin{proof}
Let $\mu\in\mc{R}_M$. Then $\mu=\delta_i +\delta_{i+1}  +\dots+\delta_{i+h}$
for some $i\in\{1,\dots,l-1\}$, and it is the highest root of the root  system
generated by $\mc{Y}(\mu)=\{\delta_i,\dots,\delta_{i+h}\}$. Indeed, any  root
chain in $A_l$ allows a simple root to appear at most once.
Therefore, the Iwasawa subalgebra generated by $\mc{Y}(\mu)$ coincides  with
$$
\Lie{n}^\prime=\bigoplus_{\alpha\in \mc{S}_\mu} \Lie{g}_\alpha,
$$
as required.
\end{proof}


\section{A counter-example}\label{counterex}
In this section we show with an example that the converse of
Theorem~\ref{al} fails when we remove hypothesis (I). More
precisely, this example suggests a general conjecture for the Lie
algebra of multicontact vector fields on a Hessenberg manifold, that
reduces to coincide with $\Lie{q}/\Lie{n}_\mc{C}$ in the case that (I)
holds.
\par
Consider the simple Lie group
$$
\group{Sp}(2,\mathbb{R})= \{A\in\group{SL}(2n,\mathbb{R}):A^t J
A=J\},
$$
where
$$
J=
\bmatrix
0&I_2\\
-I_2&0
\endbmatrix,
$$
and $I_2$ is the identity in $\group{GL}(2,\mathbb{R})$.
Take its Lie algebra
$$
\Lie{sp}(2,\mathbb{R})=\{A\in\Lie{gl}(2n,\mathbb{R}):A^t J+JX=0\}.
$$
Denote by $H_{s,t}={\rm diag}(s,t,-s,-t)$. Then the Cartan space  $\Lie{a}$ is
defined by $\{H_{s,t}:s,t\in\mathbb{R}\}$. The standard restricted  simple roots are $\alpha$ and $\beta$, where
$\alpha(H_{s,t})=t-s$ and
$\beta(H_{s,t})=-2t$.
Moreover, the set of the positive roots is  $\Sigma_+=\{\alpha,\beta,\alpha +\beta,2\alpha +\beta\}$.
The Iwasawa subalgebra $\Lie{n}$ of $\Lie{sp}(2,\mathbb{R})$ is the
direct sum of the restricted root spaces corresponding to the roots in  $\Sigma_+$, namely
$\Lie{n}={\rm span}\{E_U,E_X,E_Y,E_Z\}$, where
\begin{align*}
E_U &=
\bmatrix
0&0& & \\
-\frac{1}{2}&0&&\\
  & &0&\frac{1}{2}\\
  & &0&0
\endbmatrix    &
E_X &=
\bmatrix
  & &0&0\\
  & &0&0\\
0&0& & \\
0&2& &
\endbmatrix,\\
E_Y &=
\bmatrix
  & &0&0\\
  & &0&0\\
0&1& & \\
1&0& &
\endbmatrix    &
E_Z &=
\bmatrix
  & &0&0\\
  & &0&0\\
1&0& & \\
0&0& &
\endbmatrix.
\end{align*}
We fix coordinates on $\group{N}=\exp\Lie{n}$, the nilpotent Iwasawa  subgroup of
$\group{Sp}(2,\mathbb{R})$. Any element $n$ in N
 consists of a $4\times 4$ matrix that we write as
$$ n(u,x,y,z)=
\bmatrix
1&-\frac{1}{2}u&z-\frac{1}{2}uy&y-ux\\
0&1&y&2x\\
0&0&1&0\\
0&0&\frac{1}{2}u&1\endbmatrix.$$
In other words, we consider $\R^4$ with the group law:
$$n(u,x,y,z) n(u',x',y',z')=
\bar n(u+u',x+x',y+y'+ux',z+z'+uy'+\sfrac{1}{2}u^2x').$$
A basis of left--invariant vector fields for the Lie algebra $\Lie{n}$
is given by
\begin{align*}
U&=\frac{\partial}{\partial u};\\
X&=\frac{\partial}{\partial x}+u\frac{\partial}{\partial  y}+\frac{u^2}{2}\frac{\partial}{\partial z};\\
Y&=\frac{\partial}{\partial y}+u\frac{\partial}{\partial z};\\
Z&=\frac{\partial}{\partial z},
\end{align*}
with brackets
$$[U,X]=Y,\quad
[U,Y]=Z,\quad
[X,Y]=[X,Z]=[Y,Z]=0.$$
Take a regular element $H_0$ in $\Lie{a}$ and fix $\mc{R}=\{\alpha,
\beta,
\alpha +\beta\}$.
In order to study the multicontact vector fields on the Hessenberg  manifold that corresponds
to these data, we consider
the slice S, defined by
the equation $z=0$.
 A generating set of vector fields on S is
\begin{align*}
\overline{U}&=\frac{\partial}{\partial u};\\
\overline{X}&=\frac{\partial}{\partial x}+u\frac{\partial}{\partial
y};\\
\overline{Y}&=\frac{\partial}{\partial y},
\end{align*}
where $[\overline{U},\overline{X}]=\overline{Y}$ is the only non-zero
bracket.
\par
Consider a vector field on an open set $\mc{A}\subseteq \group{S}$,  namely
$F=f_u\overline{U} +f_x\overline{X}+f_y\overline{Y}$, where $f_u$,  $f_x$ and $f_y$
are smooth functions on $\mc{A}$.
Then $F$ is a multicontact
vector field on $\mc{A}$ if and only if
\begin{align*}
[F,\overline{U}]=\alpha\overline{U},\\
[F,\overline{X}]=\beta\overline{X},
\end{align*}
for some functions $\alpha$ and $\beta$ on $\mc{A}$.
The two conditions give rise to the differential equations
$$
\begin{cases} f_x+\overline{U}f_y=0\\
\overline{U}f_x=0\end{cases}
\qquad
\begin{cases} -f_u+\overline{X}f_y=0\\ \overline{X}f_u=0\end{cases}
$$
that imply
\begin{equation}\label{sp}
\begin{cases} \overline{X}^2 f_y=0 \\ \overline{U}^2 f_y =0.\end{cases}
\end{equation}
Notice that, since $\mc{R}$ is defined as a single shadow, we obtain  only one model-system, so
that no cross-conditions appear. Moreover, the system above coincides  with the
multicontact conditions for the coefficients of a vector field on the  nilpotent Iwasawa subgroup of
$\group{SL}(3,\mathbb{R})$ (see~\eqref{sleq}, Ch.~\ref{pr}). Therefore,
$$
MC(\group{S}) \cong \Lie{sl}(3,\mathbb{R}).
$$
We
show that this result does not coincide with the claim of  Theorem~\ref{al}, namely $\Lie{sl}(3,\mathbb{R})$
is not isomorphic to $\Lie{q} /\Lie{n}_\mc{C}$. Indeed, in the present  example
$\mc{C}=\{\omega=2\alpha +\beta\}$, and  $\Lie{n}_\mc{C}=\Lie{g}_\omega$. Hence, the normalizer
in $\Lie{sp}(2,\mathbb{R})$ of $\Lie{n}_\mc{C}$ is
\begin{equation}\label{qu}
\Lie{q}=\Lie{g}_{\alpha}\oplus \Lie{g}_{\beta}\oplus
\Lie{g}_{\alpha + \beta}\oplus \Lie{g}_{\omega}\oplus \Lie{a} \oplus
\Lie{g}_{-\beta}.
\end{equation}
Thus
$$
{\rm dim}(\Lie{q} /\Lie{n}_\mc{C})={\rm dim}(\Lie{q})-{\rm
dim}(\Lie{n}_\mc{C})=6,
$$
whereas the dimension of $\Lie{sl}(3,\mathbb{R})$ is $8$.
\par
We want to interpret the result that we obtained for $MC(\group{S})$ in  a more general context. To this end,
first we identify
$\Lie{sl}(3,\mathbb{R})$ with a subspace of $\Lie{sp}(2,\mathbb{R})$,  and secondly
we interpret this subspace in terms of the Hessenberg data.
\par
We consider  the multicontact vector fields on N, and we see how to  interpret them as  multicontact vector fields
on S. If $V=v_u U + v_x X + v_y Y + v_z Z$ is a multicontact vector  field on N, then it has the form $\tau(E)$
for some
$E\in\Lie{sp}(2,\mathbb{R})$, and its coefficients are polynomials.  Decompose $\Lie{sp}(2,\mathbb{R})$ according
to its root space decomposition, and choose a basis corresponding to  the roots.  By~\cite{CDKR1}, any root space
determines a unique (up to a real factor) multicontact vector field  that, in turn, defines a polynomial $v_z$,
namely
\begin{align*}
&\Lie{g}_\omega \leftrightarrow 1 &  &\Lie{g}_{-\omega} \leftrightarrow   (uy-2z)^2\\
&\Lie{g}_{\alpha+\beta} \leftrightarrow u &  &\Lie{g}_{-\alpha-\beta}  \leftrightarrow (y-ux)(uy-2z)\\
&\Lie{g}_{\beta} \leftrightarrow u^2/2 & &\Lie{g}_{-\beta}  \leftrightarrow (y-ux)^2/2\\
&\Lie{g}_{\alpha}\leftrightarrow (y-ux) &  &\Lie{g}_{-\alpha}\leftrightarrow u(uy-2z)\\
& \Lie{g}_0=\Lie{a} \leftrightarrow \{(uy-2z),u(y-ux)\}. & &
\end{align*}
We now define a linear map $D$ between the algebra of polynomials  associated
with $MC(\group{N})$ and the algebra of polynomials associated with  $MC(\group{S})$.
We determine $D$ by extending linearly the assignment
$$
D: v_z(u,x,y,z) \mapsto (Uv_z)(u,x,y,uy/2).
$$
Consider the image of a basis of polynomials corresponding to the root  spaces:
\begin{align*}
&D(1)= 0 & &D((uy-2z)^2)= 0\\
&D(u)= 1 & &D((y-ux)(uy-2z))= y(y-ux)\\
&D(u^2/2)= u & &D((y-ux)^2)= -x(y-ux)\\
&D((y-ux))= -x & &D(u(uy-2z))= uy\\
&D((uy-2z))= y &  &D(u(y-ux))= y-2ux.
\end{align*}
 The image of $D$ gives a set of linear independent polynomials all  solving~\eqref{sp}. Therefore, by dimensional
consideration, they generate
$MC(\group{S})$.
Recalling that $MC(\group{S})\cong \Lie{sl}(3,\mathbb{R})$,
it follows that $D$ establishes in turn a vector space isomorphism  between the
subspace
\begin{equation}\label{vec}
\Lie{a}\oplus\Lie{g}_\alpha\oplus\Lie{g}_\beta\oplus\Lie{g}_{\alpha
+\beta}\oplus\Lie{g}_{-\alpha}
\oplus\Lie{g}_{-\beta}\oplus\Lie{g}_{-\alpha -\beta}
\end{equation}
 of $\Lie{sp}(2,\mathbb{R})$ and $\Lie{sl}(3,\mathbb{R})$.
\par
We still make a further step, that allows us to define~\eqref{vec}, and  hence $MC(\group{S})$, by the Hessenberg
data, in a form that can be viewed as a generalization of the  characterization given in Theorem~\ref{al}.
Let $\mc{D}=\{\gamma\in\Sigma_-:\Lie{g}_\gamma\subset\Lie{q}\}$ and  $\mc{I}=\cap_{\mu\in\mc{R}_M}\mc{S}_\mu$.
Write
\begin{align*}
&\Lie{q}\cap\overline{\Lie{n}}=\sum_{\gamma\in\mc{D}}\Lie{g}_\gamma,\\
&\overline{\Lie{n}}^\mc{I}=\sum_{\gamma\in -\mc{I}}\Lie{g}_{\gamma},\\
&\Lie{b}=\sum_{\gamma\in\mc{D}\cup(-\mc{I})} \Lie{g}_\gamma.
\end{align*}
Notice that $(\overline{\Lie{n}}^\mc{I}\setminus  (\Lie{q}\cap\overline{\Lie{n}}))$ is
isomorphic as a vector space to $\overline{\Lie{n}}/\Lie{b}$.
In the case study $\mc{D}=\{-\beta\}$ and  $\mc{I}=\{\alpha,\beta,\alpha+\beta\}$.
We have the following vector space identifications:
$$
\Lie{a}\oplus\Lie{g}_\alpha\oplus\Lie{g}_\beta\oplus\Lie{g}_{\alpha
+\beta}\oplus\Lie{g}_{-\beta} \cong \Lie{q}/\Lie{n}_\mc{C},
$$
and
$$
\Lie{g}_{-\alpha}\oplus\Lie{g}_{-\alpha -\beta}\cong  \overline{\Lie{n}}/\Lie{b}.
$$
Hence
\begin{equation}\label{conj}
MC(\group{S})\cong \Lie{q}/\Lie{n}_\mc{C}\oplus  \overline{\Lie{n}}/\Lie{b}.
\end{equation}
We conjecture that identification~\eqref{conj}
is true in general settings.
\par
We still make a last remark. If (I) of Theorem~\ref{al} holds, the  space $\overline{\Lie{n}}^\mc{I}$ is a
subset of
$\Lie{q}$ (this follows from Lemma~\ref{norma}), so that the second  term on the right hand in the direct
sum~\eqref{conj} vanishes, there obtaining Theorem~\ref{al}.
In the case study, $\Lie{n}/\Lie{n}_\mc{C}$ is the
Iwasawa algebra whose underlying root system is of type $A_2$.
Hence, $\Lie{n}/\Lie{n}_\mc{C}$ is isomorphic to the nilpotent Iwasawa  subalgebra
of $\Lie{sl}(3,\mathbb{R})$.
 Generalizing, we can say that each shadow $\mc{S}_\mu$ generates a  vector
space that, viewed as the quotient
$\Lie{n}/\Lie{n}_{\mc{S}_\mu^c}$, with $\mc{S}_\mu^c
=\Sigma_+\setminus\mc{S}_\mu$, inherits from
$\Lie{n}$ the stratified structure of an Iwasawa nilpotent Lie algebra.
In this context, we believe that the proof of~\eqref{conj} is an  adaptation of the
proof given for Theorem~\ref{al}, and the most effort seems to be the  characterization
of the normalizer $\Lie{q}$ in terms of roots, that is an analogous of  Lemma~\ref{norma}.



\end{document}